\documentclass{article}

\usepackage{amsmath, amssymb, anysize, datetime, chemarr, latexsym, stmaryrd, xspace}
\usepackage{stackrel, color}
\usepackage[usenames,dvipsnames]{xcolor}
\usepackage{amsfonts,epsfig,colordvi,graphics,textcomp}
\usepackage{amsthm}
\usepackage{longtable}
\usepackage{array}
\usepackage{verbatim}
\usepackage{afterpage}
\usepackage{pdfpages}
\usepackage[font=small,justification=justified]{caption}
\usepackage{bbm}
\usepackage{listings}
\usepackage{hyperref}
\usepackage[nopostdot,acronym,nomain,toc,nonumberlist]{glossaries}
\usepackage{cases}
\usepackage{titlesec}
\usepackage{chemfig}
\usepackage[font=small,justification=centering]{subcaption}
\usepackage{mathrsfs,mathtools}
\usepackage{enumerate}
\usepackage{makeidx}
\usepackage[super]{cite}
\usepackage[document]{ragged2e}

\newtheorem{thm}{Theorem}[section]
\newtheorem{lem}{Lemma}[section]
\newtheorem{prop}{Proposition}[section]
\newtheorem{define}{Definition}[section]
\newtheorem{cor}{Corollary}[section]
\newtheorem{rem}{Remark}[section]

\newtheorem{ppts*}{Properties}[section]
\newtheorem{ansatz}{Ansatz}

\newtheorem{eg}{Example}[section]
\newtheorem{assmp}{Assumption}

\let\oldproofname=\proofname
\renewcommand{\proofname}{\rm\bf{\oldproofname}}

\newcommand{\grad}{\nabla}

\renewcommand{\,}{,\:}

\newcommand{\E}{\mathbb{E}}

\renewcommand{\P}{\mathbb{P}}
\newcommand{\R}{\mathbb{R}}
\newcommand{\C}{\mathbb{C}}

\newcommand{\A}{\mathbf{A}}
\newcommand{\F}{\mathbf{F}}
\newcommand{\Om}{\Omega }
\newcommand{\om}{\omega }

\newcommand{\N}{\mathbb{N}}
\newcommand{\Ra}{\Rightarrow}

\newcommand{\dom}{\mathfrak{D}}

\newcommand{\Lp}{L^{p}(\dom)}

\newcommand{\Lu}{\mathcal{L}}

\newcommand{\lam}{\lambda}

\newcommand{\lap}{\Delta}

\newcommand{\const}{{\mbox{\large{\rm k}}}}
\newcommand{\arbconst}{\const_{_{_*}}}

\newcommand{\lbd}{{\rm m}}
\newcommand{\ubd}{{\rm M}}

\renewcommand{\u}{{\bf u}}

\newcommand{\B}{{\bf B}}

\newcommand{\dto}{\downarrow}

\newcommand{\la}{\langle}
\newcommand{\ra}{\rangle}

\renewcommand{\S}{\mathcal S}

\newcommand{\FT}{\mathfrak{F}}
\newcommand{\IFT}{\mathfrak{F}^{-1}}
\newcommand{\SW}{\mathcal{S}}

\newcommand{\secAng}[1]{\kappa_{_{#1}}}
\newcommand{\secOpSp}{\mathbb S}

\newcommand{\piT}{{\pi \over 2}}

\newcommand{\alp}{\alpha}

\renewcommand{\F}{\mathcal{F}}
\renewcommand{\B}{\mathcal{B}}

\newcommand{\evolOp}{(U_{t,s})_{\underset{t,s\in[0,T]}{t\ge s}}}
 
\newcommand{\Ug}{\mathfrak{U}}

\newcommand{\LP}{L\'{e}vy process }
\newcommand{\LM}{L\'{e}vy measure }

\newcommand{\CSG}{C_0\text{-SG}}
\newcommand{\Csg}{\rm{C_0\text{-semigroup}}}
\newcommand{\hLt}{\hat L^2(\R^d)}

\newcommand{\delInt}{\Upsilon}

\title{From Microscopic SDE to Stochastic Macroscopic Equation: A Framework for Modeling Diversity in Cancer and Other Complex Living Systems.} 

\author{Sandesh Athni Hiremath  (sandesh.hiremath@mv.uni-kl.de) \footnote{
Department of Mechanical and Process Engineering, TU Kaiserlautern, Gottlieb-Daimler-Straße 42, 67663 Kaiserslautern, Germany.}}

\date{}

\usepackage[a4paper, total={7in, 9in}]{geometry}

\begin{document}


\maketitle

\justify

\begin{abstract}
Biological living systems in general exhibit complex and diverse dynamics. The latter, in particular, is essential, since diversification increases the odds of survival of an organism while reducing the risk of extinction of the population. Primarily, diversification is a consequence of the randomness in the replication process of a biological cell, which eventually manifests into a motley set of macroscopic features of an individual. These heterogeneous features of individuals constitutes for diversity in population. Cancer is a prime example of such a complex system where the transformed cells exhibit plethora of disparate features, which in turn makes modeling their dynamics quite challenging. In this paper we consider cancer as a prototype of a complex living system and provide two contrasting perspective for studying and modeling it. Based on this we illicit a deeper role of diversification in the evolution of cancer. Following this, we ask ourselves how can model these diverse dynamics in a multiscale setting.  
The current multiscale modeling techniques lack the ability to capture this diversity at the macroscopic level of experimental observation, thus are inappropriate for describing the behavior of a selfish organism like cancer. 
We address this shortcoming by providing an abstract but mathematically rigorous framework for deducing stochastic evolution equations at the macroscopic level starting from a microscopic description of the involved dynamics. We achieve this by making use of the connection between stochastic process and the semigroup operator generated by them. In particular, we look at the semigroups generated by Levy processes and their connection with the characteristic functions and Levy symbols. The latter turns out to represent pseudo-differential operators using which we eventually provide a mechanism for constructing stochastic evolution equations. Altogether, this provides the framework for modeling diverse dynamics at the macroscale starting from the microscale.
\end{abstract}



\section{Introduction \label{sec:intro}}
\begin{figure}
\centering
	\includegraphics[scale=.3]{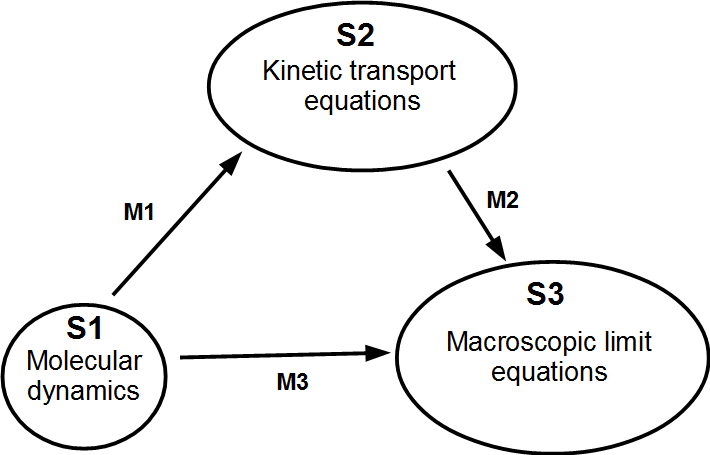} 
\caption{Relation between representation of dynamics at different scales. \label{fig:micMacScaleDiag}}
\end{figure}

The field of multiscale modeling is an active field of research with an aim of accurately and rigorously deducing an equation or a system of equations for the observed phenomenon starting from first principles.  The fundamental physical laws of molecular interactions i.e. the basic Newtonian laws of motion serve as the first principles. Such accurate and mathematically rigorous description of the observed phenomenon is indeed helpful to characterize some important parameters such as cell motility or heat conductivity or viscosity, of the respective phenomenon under investigation such as cancer migration or temperature distribution or fluid flow respectively. 
\noindent
Figure \ref{fig:micMacScaleDiag} illustrates the typical methods of multiscale modeling. Contributions in the direction of S1 to S2 dates back to the pioneering works of J.C. Maxwell and L. Boltzmann on the kinetic theory of gases (fluids in general). A first rigorous proof for the Boltzmann equation starting from Newtonian laws is due to O. Lanford [\citenum{Lanford1975}]. The result was slightly unsatisfactory due to its validity for short time dynamics i.e. just local existence, thereby leaving open the important questions about long time dynamics, i.e. questions relating to the global existence. Some of these gaps were overcome by the works of Illner, Pulvirenti and Shinbrot [\citenum{Illner1984, Pulvirenti1987, Illner1986}]. It was not until 1989 that a more general global result was made available due to DiPerna and Lions [\citenum{DiPerna1989}].
A significant amount of work (both theoretical and applied) has also been carried out in the direction of $S2$ to $S3$ i.e. in the direction of obtaining different types of hydrodynamic limit equations starting from the Boltzmann equation. Hilbert also worked on this problem [\citenum{Hilbert1916}] where he used asymptotic expansions (now known as Hilbert expansions, a slight variant of the Chapman-Enskog expansion, thus sometimes also referred to as Chapman-Enskog-Hilbert expansion) to show that, under appropriate physical setting, the $0$th order term of the expansion solving the Boltzmann equation is a Maxwellian (i.e. the Maxwell-Boltzmann distribution)[\citenum{Delale1982}]. These works not only provide a strong foundation but also a starting point for the deduction of most macroscopic equations. 

There has also been some progress in the direction of $S1$ to $S3$ due to the advances in stochastic calculus. A first attempt in this direction was nearly 60 years ago and is due to C.B. Morrey [\citenum{Morrey54}]. Though the work could be considered to be lacking some mathematical rigor, it did introduce in a subtle way the notions of  \emph{local equilibrium} which is nothing but the pointwise quasi-steady sate of the microscopic dynamics in the appropriate limit and after appropriate space time rescaling. We also see here that the author implicitly uses the notion of ergodicity to justify the existence of equilibrium distribution for which the entropy is supposed to be maximized. 
The latter approach was used by S.R.S Varadhan in [\citenum{Varadhan1993}] and S.Olla in [\citenum{Olla1993}] where they perturbed the N-body Hamiltonian system with white noise in order to ensure ergodicity of the system.  An implicit assumption in these techniques is that the microscopic dynamics is random or in other words random dynamics is incorporated in the description of microscopic evolution. A more intuitive explanation towards this ideology of modeling is provided in [\citenum{DeMasi1991,Kipnis1999}].

One of the major drawback in the above multiscale modeling framework is their inablity to capture complex dynamics of a living system. The main reason being that the framework is mainly developed and for innately passive particles. Many adaptations to the above framework have been used to model complex systems. In [\citenum{Iannini2016}], authors have adopted kinetic-theory based Boltzmann type equation to model active vehicular traffic wherein the notion of active behavior of a vehicle (modeled as a particle) was introduced. More recently, in [\citenum{Aristov2019}] kinetic approach was used to model biological system as a dissipative system arising from non-equilibrium boundary interaction. More precisely, they based their model on the approach of Schrödinger wherein the bio-system is feeds on the negative entropy. On the other hand there are considerable works in the direction of adapting the classical kinetic theory for passive particles to a kinetic theory for active particles (KTAP) with the objective of modeling complex dynamics of living systems. For example in [\citenum{Bellomo2013, Bellomo2017, Bellomo2021}] phenomenological equations at the mesoscopic level have been proposed for describing complex cell-cell and cell-tissue interactions. Since these equations are devised take into consideration multi-scale as well as non-conservative interactions, it allows for modeling a wide range of complex microscopic and macroscopic dynamics [\citenum{Bellomo2004,Bellomo2008,Bellomo08,Bellomo17}]. Although very powerful tool, there are two major shortcomings of KTAP: (i) since it mainly devised from phenomenological arguments it lacks rigorous derivation from first principles, (ii) the obtained equations are mainly of deterministic type and thus lack the ability to model diversity and heterogeneity of observed macroscopic dynamics. 
One way to remedy that is to by incorporating non-determinism at the macroscopic level. 
This is exemplified in the works of [\citenum{SMAMCI15, HiremathZhigun2016, SAIG16}], where stochastic models were proposed to describe the complex dynamics of cancer cell migration. These models illustrate how randomness play an essential role in bringing forth some of the rare, transient, and interesting events  while the averaged dynamics fails to do so. Even more astonishing is that the averaged dynamics tell a very different (or rather an incomplete) story. This is the main drawback of using deterministic models to explain a very diverse, complex, and generic problem such as cancer invasion. A more elaborate comparison between stochastic and deterministic models in acid mediated cancer invasion is provided in [\citenum{hs16b}. Motivated by the interesting dynamics induced by stochastic models, we ask ourselves the following question: how can one deduce a stochastic macroscopic equation starting from  microscopic equations? 



In this paper we answer this question by providing an abstract mathematical framework for deducing stochastic evolution equations at the macroscopic level starting from a microscopic description of the involved dynamics. We achieve this by first observing that a certain class of L\'evy processes generates transport equations at the macroscopic level. Motivated by this we ask ourselves what different types of macroscopic equations would be possible for other types of L\'evy processes. The L\'evy-Ito decomposition formula and L\'evy-Khintchine formula gives a hint for answering this question. Since we want to deduce stochastic evolution equations, we plan to investigate the situation where microscopic equations are driven by different types of L\'evy processes, where the type of the L\'evy process i.e. the law of the L\'evy process, is randomly changing with respect to time. Intuitively speaking, because different L\'evy processes are characterized by different L\'evy symbols, a randomly varying L\'evy noise at the microscopic level should result in stochastic evolution equations at the macroscopic level. To make this mathematically rigorous we invoke the theory of semigroups and pseudo-differential operators which enable use to generate a specific class of random operators which in turn allows us to construct a two parameter semigroups. This eventually lead us to a random non-autonomous Cauchy problem at the macroscopic level. 

To this end, we first start with a short discussion in Section \ref{sec:CanSelfOrg} about some fundamental properties of a complex living system and their implications by considering cancer as a prototypical example. Motivated by these discussions, in Section \ref{sec:M2MACP} we introduce the necessary concepts required for modeling diversity in living systems. This entails discussion about the properties of L\'evy processes, L\'evy-Ito decomposition and the L\'evy-Khintchine formula. The latter introduces a class of pseudo-differential operators which in turn generates a class of operators. The characterization of the domains of these operators is discussed in Section \ref{sec:ndfFunSp}.  Based on these connections we provide a framework for constructing stochastic evolution equations in Section \ref{sec:conRandOp}. Finally we apply this approach for modeling acid mediated cancer invasion in Section \ref{sec:fracSPDE}

\begin{figure}
\centering
	\includegraphics[scale=.5]{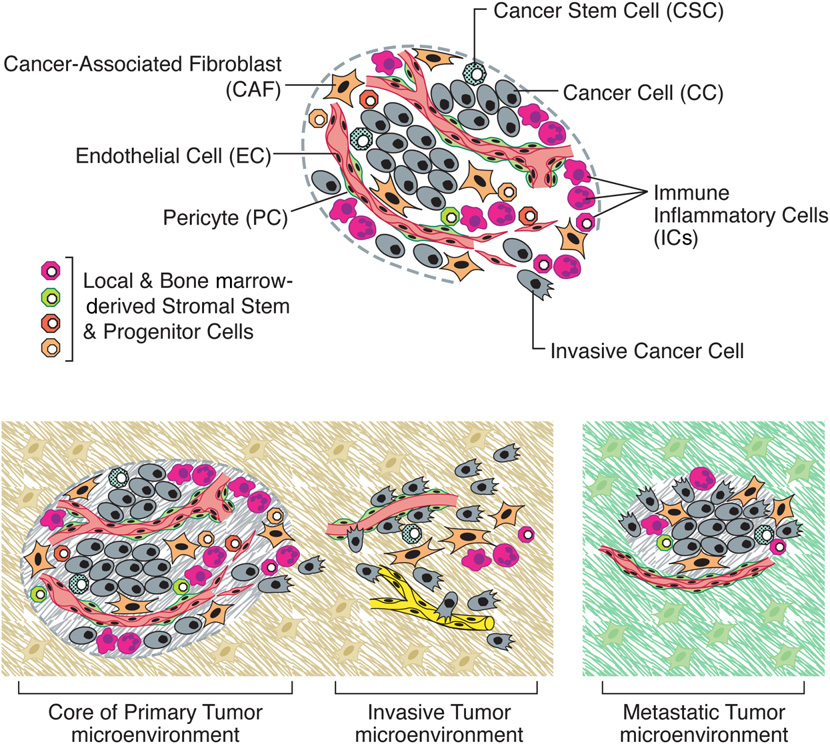} 
\caption{Illustration of a tumor microenvironment [\citenum{HAWEI11}].\label{fig:tumMicEnv}}
\end{figure}
\vspace*{-.5cm}
\section{Two contrasting perspectives of cancer \label{sec:CanSelfOrg}}
There are two main different philosophies for the cause of cancer. The classical viewpoint is that cancer is a neo-transformation, i.e. cells are transformed to some new (probably non-existent in the evolutionary history) type of cells, while a new emerging philosophy being that cancer is a paleo-transformation i.e. cells are transformed to primitive unicellular type of life forms. These ideologies are in fact two faces of the same coin and provide two distinct microscopic views for modeling the problem of cancer. The following discussion is motivated by  [\citenum{Merlo2006, Davies2013, CasasSelves2011, Niklas2014, MYECODEV15, EVOLBIOCONAPP08, EVOCAN15, Chen2015, GCTFT07, WDTHAG76, Yates2012}]. 
\subsection{Cancer: An ecosystem of cells}
\begin{figure}
\centering
\subfloat[Micro interactions{[\citenum{HAWEI11}]}. \label{fig:tumInterSigMicEnv}]{\includegraphics[scale=.15]{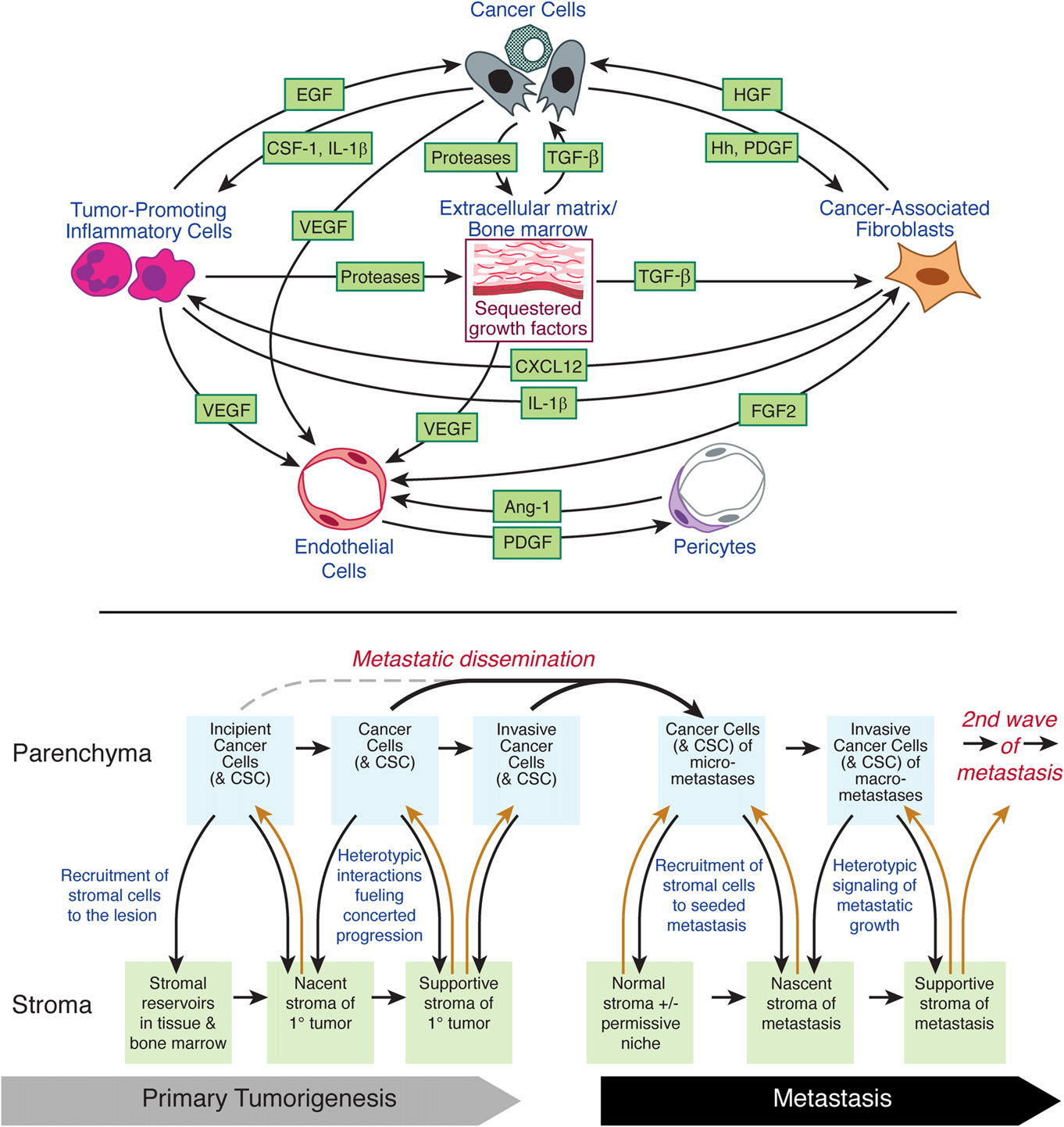} } \qquad 
\subfloat[Mechanisms of cell movement  {[\citenum{Friedl2003}]}.\label{fig:canMigStrat}]{\includegraphics[scale=.3]{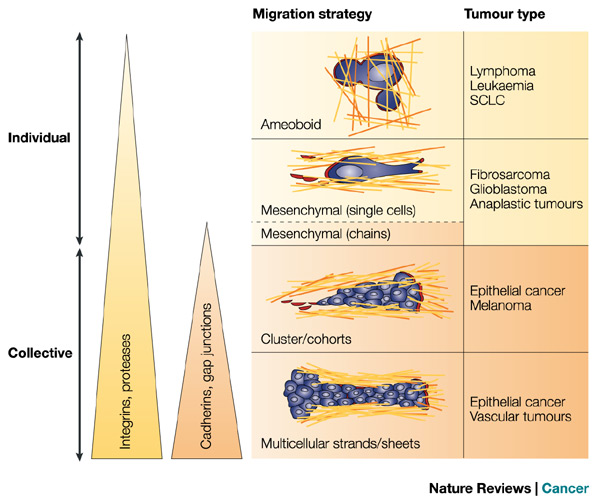}}
\caption{Illustration of plausible interactions within a tumor microenvironment and different mechanisms of cancer cell movement}
\end{figure}

\begin{figure}
\centering
	\includegraphics[scale=.3]{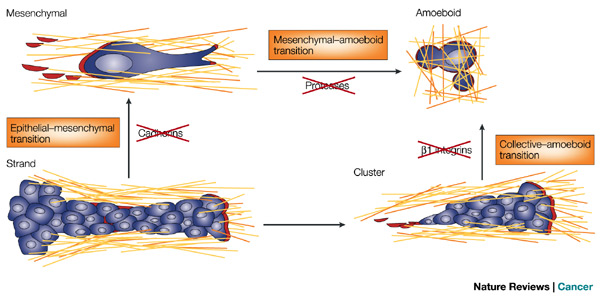} 
\caption{Illustration of plasticity in cell movement, where depending on the molecular interactions cells can switch between different mechanisms of movement. Loss of proteases leads to a switch from EMT to amoeboid type of movement, loss of cell-cell adhesion due inhibition of cadherins leads to a switch from collective epithelial sheet movement to EMT type single cell movement, loss of cell-ECM adhesion due to the inhibition of $\beta_1$-integrins leads to a switch from cluster type collective movement to amoeboid type single cell movement [\citenum{Friedl2003}]. \label{fig:canMigPlast}}
\end{figure}

\begin{figure}
\centering
	\includegraphics[scale=.5]{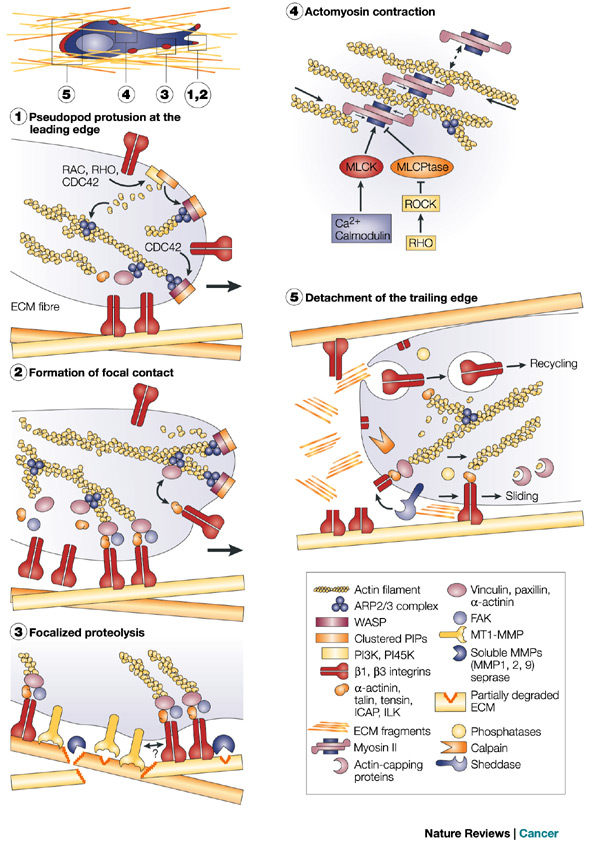} 
\caption{Illustration of the different bio-mechanical processes or steps involved in the movement of a cell [\citenum{Friedl2003}].\label{fig:canMigSteps}}
\end{figure}

Cancer may be interpreted as an ecosystem of neoplastic and stromal cells inhabiting the tissue containing vital resources for cell sustenance. The ecosystem is composed of different types of neoplastic cells such as cancer cells, cancer stem cells, cancer associated fibroblasts, etc., and normal healthy cells such as immune cells, endothelial cells, pericytes, etc. Apart from different types of cells, the micro-environment is also composed of extracellular substances such as interstitial fluid, blood vessels, ECM, acid byproducts, and other various vital biomolecules. Figure \ref{fig:tumMicEnv} illustrates the composition of a generic tumor environment. All these elements together form an interactive ecological system which manifests itself as a tumor. Figure \ref{fig:tumInterSigMicEnv} illustrates some of the possible sets of interactions between different inhabitants of the tumor ecosystem, such as interactions between cancer cells and endothelial cells, endothelial cells and pericytes, immune cells and cancer associated fibroblasts, etc. Various such complex interactions is what orchestrates a metastatic progression of cancer. An obvious implication of this is that the development and progression of a tumor is very much dependent on its  environment. Due to the heterogeneity in the tissue architecture and composition, the way an incipient tumor develops would be different therefore leading to diverse and heterogeneous cancer types. The former (i.e. diversity of cancer types) accounts for the differences between various cancers such as prostate cancer, breast cancer, melanoma, glioma, cervical cancer, etc., while the latter (i.e. heterogeneity of cancer types) accounts for differences in two or more samples of the same type of cancer across different individuals or same individuals. This motivates the study of cancer from an ecological perspective. \\
Recent advances in gene sequencing techniques has resulted into  comprehensive listing of somatic mutations responsible for cancer development. These studies indicate that the evolution of cancer genome is not only diverse, i.e. it varies across different tumor types but is also heterogeneous, i.e. it varies  within same tumor arising across different individuals and within a single individual [\citenum{Yates2012}]. Because of these variations and because of ecological constraints such as competition, predation, parasitism, mutualism, etc., clones with varying genome experience diverse selection forces, thereby producing significant variation in the evolutionary path of a tumor [\citenum{Merlo2006}]. Here the selection forces could be both natural and/or artificial selection. The former is induced due to genetic mutations that increase the fitness condition and are heritable after successive reproduction. The latter is induced via external interferences such as drugs and therapy which choose for more resistive clones of the population. Also, these selection forces act at various levels, namely- at the level of individual cells, at the level of population and at the level of the tumor itself [\citenum{Merlo2006}]. Altogether, they work towards conferring malignancy to the tumor. One of the key malignant feature is the ability to invade and metastasize, which is majorly influenced by the instability and heterogeneity at the level of genes, clones, population, micro-environmental ecosystem and the organism itself. Cancer cells are endowed with diverse cellular mechanisms that enable them to disassociate from the primary tumor and to start wandering through the local tissue and eventually enter the circulatory system. A comprehensive review on some of the different mechanics of movement employed by cancer cells can be found in [\citenum{Friedl2003}]. These mechanisms can be categorized as single cell movement and collective movement. The former involves EMT and amoeboid type movement, while the latter involves clustered or cohort movement  and movement of sheets or strands of cells.  Figure \ref{fig:canMigStrat} illustrates these different categories of cell movement, while Figure \ref{fig:canMigPlast} depicts how cells switch between different movement mechanisms as a response to the changes in their intracellular and extracellular molecular pathways. For the sake of completeness we also include Figure \ref{fig:canMigSteps} depicting the different bio-mechanical steps involved in the process of cellular motion, which includes (i) the formation of protrusions such as lamilapodia, filopodia, (ii) formation of focal contacts i.e. attachment of cells to ECM fibers, (iii) cleaving of ECM via proteolysis near the focal points, (iv) contraction of the cell due to the binding of myosin II to actin filaments, (v) finally detachment of the rear end which results in a displacement of the cell. \\
To summarize, the genetic instability within cancer cells generates a wide fitness landscape and the micro-evnironmental ecosystem selects for cells that are fit enough to survive and reproduce further. Many iterations of such Darwinian dynamics leads to the development of a diverse, heterogeneous, malignant, and invasive tumor.
Next we shall look at a contrasting viewpoint for the process of carcinogenesis and for the behavior of cancer cells. 

%

\subsection{Cancer: A selfish organism}
A typical perspective is that cancer is a transformation of normal healthy cells into a new kind of cells that possess the hallmark of neoplasms. In this philosophy, individual cells, that result after neoplastic transformation have no individuality, instead they are seen as a mere chemical system and are assumed to be lacking the ability to act for their own self interest i.e. they are, in some sense, assumed to be lacking any motive or intent. The whole process of cancer development and progression is viewed from a mechanistic perspective, where the focus is more on the bio-chemical mechanisms that make the cells behave in a certain way. This is simply a physicists perspective, where all matter, including living cells, are considered to be purely mechanistic. Such a philosophy fails to explain the motive or the purpose for the cells to behave the way they do and is too narrow to model a complex and diverse disease such as cancer. The notion of individuality and self interest is a vital factor that contribute towards the complexity of the disease and could in turn be helpful for developing robust models for the cancer invasion.
In contrast to this, a fairly new philosophy is to look at cancer as a consequence of a de-evolution process (i.e. reverse evolution, evolution towards primitive life form). In this perspective cancer is seen as a transformation from collaborative multicellular cells to a selfish free unicellular type of cells. This neither contradicts nor denies the existence and importance of the bio-chemical mechanisms that are responsible for inducing various behavior of cancer cells. Instead, it assigns a notion of individuality to the transformed cells. Since the DNA has been passed on via millions of years of evolutionary processes, from unicellular  free living organisms to multicellular collaborative organisms [\citenum{Davies2013}], it is not totally surprising to hypothesize that mutations can reactivate or enable certain parts of the DNA, phenotypically making them more like primitive life forms. The analogy becomes clear if we relate the behavior of cancer cells and other unicellular organisms like bacteria. Cancer cells grow uncontrollably, they scavenge for nutrients, they move and migrate to other more fertile space, all of which are characteristic of a primitive unicellular organisms. Not only that, some bio-chemical processes are similar to those of primitive unicellular cells, e.g. cancer cells rely on primitive metabolic pathways like glycolysis, they avoid DNA damage by maintaining the length of telomere strand  even after successive replication, just like germ cells and unicellular organisms [\citenum{CasasSelves2011}]. Moreover, looking at the evolutionary pattern of multicellular tissues, organs and so forth, the analogy becomes even more obvious. According to [\citenum{Niklas2014}], the origin of multicellularity is attributed to the genetic variations that led to the acquisition of cell-cell adhesion, communication, cooperation, and specialization. The acquisition of such features can be possible if neighboring cells   show (i) alignment-of-fitness, i.e. they have no inter-cellular conflicts and (ii) export-of-fitness i.e. are able to independently collaborate towards a sustained reproductive goal. Based on this there are two possible paths to multicellularity: (i) unicellular organism $\Ra$ clonal cells $\Ra$ multicellular organism or a more direct path such as (ii) siphonous (a giant unicellular algae cell) $\Ra$ multicellular organism. Now, if we ask ourselves as to what would happen if mutations in a multicellular organism lead to the loss of its vital functionalities, then it is clear that we either get a bunch of clonal cells or just independent individual unicellular organisms which is exactly the scenario in a tumorous tissue. This philosophy of de-evolution and its connection to cancer dynamics is presented next.

\subsubsection{Evolution and de-evolution of cells:}

According to Darwin's theory of evolution, multicellular organisms have evolved, over a period of millions of years, from unicellular ones through random processes such as genetic alterations and natural selection. The tree of evolution clearly indicates a bottom up process where the unicellular organisms are at the root while multicellular organisms occupy their place at the trunk and leaves of the tree. This is the case because genetic/chemical alterations and environmental forces propelled unicellular organisms to aggregate and start to collaborate which eventually gave rise to functionally sophisticated life units such as tissues, organs, and in turn organisms. There are also cases where complex multicelluar organisms have evolved (or de-evolved) to simpler forms of life like viruses. For example, certain species of Cnidaria have de-evolved from multicellular free-living autonomous animals into parasitic organisms comprising of only a handful of cells. Such parasites constitute the class of organisms called Myxozoans [\citenum{MYECODEV15}].  This unique reverse evolution is again due to genetic alterations and environmental selection. \\
Cancer is a disease in which some of the cells within a multicellular life unit (e.g. tissue or organ) have functionally turned into a more primitive life form. What this means is that tumor cells are no more adhering to the rules of their pre-destined functionality but instead they become rogue cells acting in a primitive and selfish way which is often deleterious to the multicellular life unit that it is still a part of. A simple example is the case of melanoma (a type of skin cancer): the pigment producing cells called melanocytes undergo transformations (via mutations induced mainly through exposure to UV radiations) and are no more just responsible for melanine production but they start exhibiting the features of an independent unicellular organism. This means that melanocytes after being carcinogenized start replicating profusely, boost their nutrition supply (e.g. via angiogenesis), may even start to migrate, etc. All these are features of an autonomous living organism (say like protozoans). Thus carcinogenesis is, in a very vague sense, similar to the de-evolution of Cnidarians to Myxozoans, with the main difference that carcinogenesis induces a switch from a stable, cooperative (altruistic), and functionally pre-determined behavior to an instable, selfish, and functionally free behavior while the de-evolution process of Cnidarians to Myxozoans involves a rather drastic physical transformation. This novel analogy can be backed up by the biological finding that many protozoan life cycles features have been observed in mammalian cancers undergoing polyploidy [\citenum{EVOLBIOCONAPP08}].\\
Another aspect of cancer is that it is not a modern disease, but an ancient one, whose roots go back to the time of the beginning of the cellular life itself.  In the early stages of unicellular life, natural selection favored organisms that were fast growing, space invading, and had low mortality rate etc. However, for the development of multicellularity there had to be a paradigm shift, i.e. natural selection had to favor cooperative cells over selfish ones. According to [\citenum{EVOCAN15, EVOLBIOCONAPP08}] this can be a reason for the rarity of cancers in multicellular organisms and also a reason why the occurrence of cancer represents a reversal from altruistic multicellular cells to selfish unicellular cells. Moreover, according to [\citenum{CasasSelves2011}] cancer occurs especially during the later stages of life mainly due to loss of tumor suppressor functions. Evolutionarily, this makes sense because there is less chance of reproduction at  an older age, thus less chance of genetic contribution to the progeny; as a result there is no evolutionary pressure to retain tumor suppressor functionalities. Most of such functionalities are active during the reproductive age, when the chances of creating offsprings are quite high. A more recent study [\citenum{Chen2015}] indicates that metastatic tumors develop not only due to gain-of-function type mutations affecting the oncogenes, but also due to loss-of-function type mutations affecting the tumor-suppressor genes. Such loss-of-function type mutations may lead the loss of genetic constraints responsible for maintaining multicellularity and tissue homogeneity. \\
Based on these considerations, cancer may be viewed as a population of diverse unicellular organisms which are striving to survive within the dynamic micro-ecosystem generated by the host tissue. These unicellular organisms switch between altruistic and selfish behavior in a such a way as to improve their or their community's survival and longevity. \\
Motivated by these intuitive arguments, one is led to the task of finding the fundamental reason for the complexity and diversity observed in cancer. More importantly, one is also led to the question of how to model such diverse and complex living systems. In order to answer such fundamental questions we have to resort to an abstract formulation of  living systems and analyze its implications. This is addressed in Section \ref{sec:divCan} below. 

\subsection{Diversity of cancer \label{sec:divCan}}

\begin{figure}[!htb]
\begin{center}
\subfloat[Hypothetical living systems \label{fig:divSystems}]{
\includegraphics[trim=0cm .8cm .8cm .3cm, clip,scale=.3]{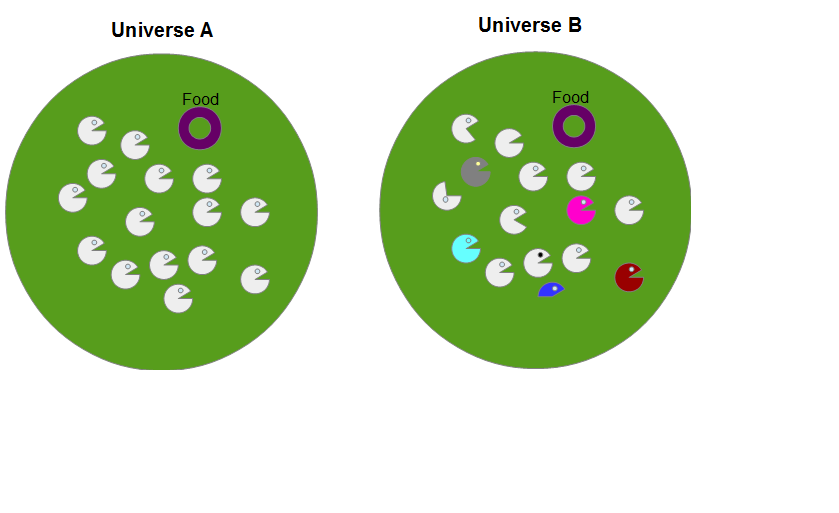} } \quad
\subfloat[Decay profile of the food densities \label{fig:randomRep}]{\includegraphics[trim=.5cm 0cm .5cm .5cm, clip, scale=.3]{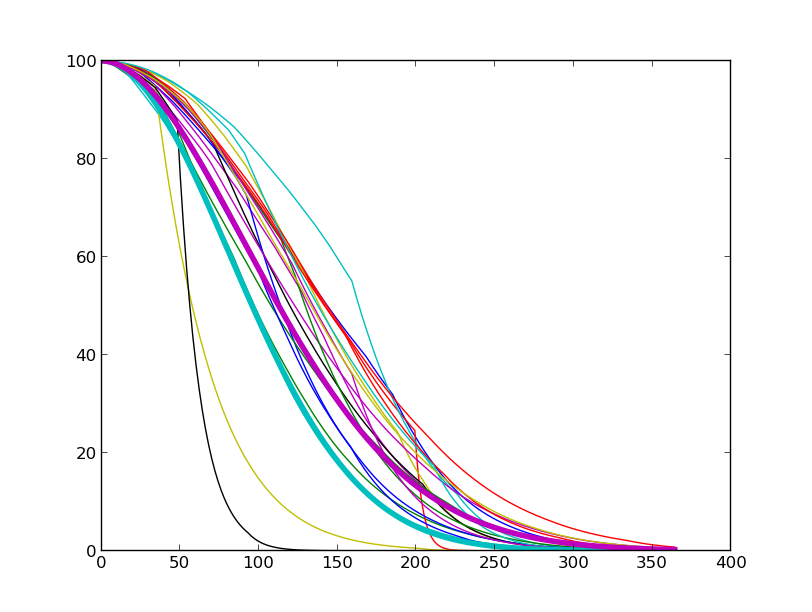}}
\caption{Diversification of living systems. The figure on the left is a cartoon of two hypothetical systems $X$ and $Y$ belonging to Universe A and Universe B, respectively. The difference between system $X$ and $Y$ is that the former replicates identically while the latter replicates with a certain modification. These modifications manifest as different rates of consumption of the food. The figure on the right (Figure \ref{fig:randomRep}) depicts the decay profile of the food density. The thick cyan colored curve represent the decay profile of the food density A, while the thick pink colored curve represents the mean decay of the food density B.}
\end{center}
\end{figure}

For modeling any natural phenomenon it is important to understand the fundamental motive and purpose as to why such phenomenon occurs in the first place, and in the context of biology (especially population biology) it is imperative to understand the purpose and aim of the underlying phenomenon being modeled. There is a fundamental difference between mere chemical systems and  organisms or general living systems (which are more than just chemical systems). The term living system makes it obvious as to what the goal or the aim of the system is, it is simply to stay alive for as long as possible. Unlike lifeless chemical or physical systems where the dynamics of certain phenomena are governed by physical and chemical laws, the dynamics of a living system is mainly driven by the laws of survivability. Briefly, what this means is that a living system (which is inherently assumed to be complex) is capable of bending or evading or circumventing certain phsio-chemical laws, in favor of certain other physio-chemical laws, with the sole purpose of surviving. It is absolutely important to note here that the previous statement does not state that the physical or chemical laws are not applicable. Of course everything that happens is indeed according to some physical laws, but the point is that under different conditions and different circumstances different laws are in action and this changes the results of the outcome. Thus living systems strive to be alive by altering the conditions and/or circumstances and/or environment such that the laws of physics that are in action favor their survivability. Therefore already this basic law of survivability provides a glimpse for the reason for complexity of living systems. In order to elaborate more on this we provide an novel abstract characterization of a generic living systems using the following properties:\\[.5ex]
{\bf Autonomy:} A living system must have necessary mechanisms to produce or harvest energy in order to sustain itself and to survive and reproduce successfully. Such mechanisms need to exist only up until successful reproduction and upbringing of the progeny to independence. This is also the stage where new information is acquired and the system undergoes suitable changes and adaptations in order to better sustain itself.\\[.5ex]
{\bf Replication:} A living system must be able to successfully replicate and create either a clone or a modified version of itself. The functionality of replication also includes the job of raising and nurturing the newborn up to its independence, i.e. up to the stage where the newborn achieves its autonomy.   In this sense the replicant  must also be a living entity, although not immediately but after reaching the state of independence. Here again it is implicitly assumed that the new born inherits the mechanism of reproduction from its parent and is able to reproduce upon attaining autonomy.  We shall omit those systems that are not capable of replication from being called as a living system, since there is no possibility for such systems to continue their tree and consequently are destined for extinction. The reproduction stage is where the transfer of information happens. All or part of  the  information, mainly the mechanism of replication and the tricks and trades of survivals, acquired up until the time of replication is passed on to the child living system during the process of creation or/and  during the time of nurturing until it reaches independence.
The above two properties are general enough to not only cover all known biological life forms till now but also to include non-biological life forms like computer programs and abstract logical systems. The two properties can be in short called  \emph{axioms of life} and they are the fundamental driving force of all things alive. 
Another crucial aspect that can be deduced from the axioms is that: in order to increase the chances of survival and thus increase the chances of reproduction, it is necessary that the living systems diversify their mode of operation or sustenance. What this means is that: replication should not result in an exact copy of the parent but instead it should be a modified copy. More precisely, the progeny should be a modified version of the progenitor. This is not a necessity but just a beneficial strategy for increasing the odds of survival. We shall elaborate on this via an example.\\
 Consider two living systems $X$ and $Y$ belonging to Universe A and Universe B respectively. Universes A and B are independent of each other and have the same type of food resource named as $A$ and $B$, respectively. Each living entity $X$ and $Y$,  for its sustenance, consumes food $A$ and $B$  at  rates $k_x$ and $k_y$, respectively. Letting $a(t)$ and $b(t)$ denote the concentrations of $A$ and $B$ at time $t$, we see that $a(t) = a(0) e^{-k_x t}$ and $b(t) = b(0) e^{-k_y t}$. If $X$ and $Y$ replicates at time $t_x$ and $t_y$, respectively, then for any $t > \max(t_x, t_y)$ we get that $a(t) = a(t_x)e^{-2 k_x (t-t_x)}$ and $b(t) = b(t_y) e^{-2 k_y (t-t_y)}$. Depending on initial concentrations $a(0)$ and $b(0)$, rates $k_x$ and $k_y$, and replication times $t_x$ and $t_y$, either $a(t) \ge b(t)$ or vice versa. The interesting case for us to consider is $a(0) = b(0) = r_0$, $k_x = k_y = k$ and $t_x = t_y = t_r$.  The only difference between both systems is that, when $A$ replicates it produces an exact copy of itself, while when $B$ replicates it produces a modified copy of itself. This modification  manifests as a change in the consumption rate $k_y$ of $B$. Let $k'$ be the consumption rate of the replicant of $Y$. Then for $t> t_r$ we have that $a(t) = a(t_r) e^{-2 k (t-t_r)}$ while $b(t) = b(t_r) e^{-(k+k')(t-t_r)}$. Now depending of whether $k' > k$ or $k' < k$ we get either $b(t) < a(t)$ or $b(t) > a(t)$. Hence, for $k' < k$, $Y$ improves its chances of survival and for $k' > k$, $Y$ worsens its chances of survival. Therefore of all the three cases, i.e. $k' = k$, $k' > k$ and $k' < k$, only the last case improves the longevity of life of $X$ or $Y$. Now after multiple iterations of replication, population $Y$ (provided $k' < k$) will have outlived  $X$. \\
Figure \ref{fig:randomRep} illustrates the decay profile of the food density of $A$ and $B$ belonging to Universes A and B, respectively. The different consumption rates of each replicant of system $Y$ are sampled from a uniform distribution. Based on this we see that in the expectation system $Y$ has better odds of survival since the food is available for a longer time. It must be noted here that in this example the choice of uniform distribution was purely for the sake of simplicity. In fact in reality it may very well belong to other distributions. Thus the possibility that the modifications of a replicant may belong to different distributions establishes the fundamental basis for the diversity and heterogeneity that prevails among living systems.
\noindent
From this simple example it becomes clear that maximization of survivability requires diversification but of course with the risk of extinction. This example also illustrates that if the process of replication is random or uncertain, then the game of survival  also becomes random. This motivate one to consider that randomness and uncertainty are one of the primary, if not the only, source of diversity in living systems. Since all biological living systems replicate randomly, diversification and randomness becomes an inherent property of their evolution. Consequently, we may conclude that the \emph{axioms of life} consist of three properties: \emph{autonomy, replication} and \emph{diversification}.\\
\noindent
Motivated by these considerations, the rest of the work focuses on the topic of multiscale modeling of complex living systems such as cancer. To this end, in Section \ref{Chap6} we present a novel framework for deducing stochastic equations at the macroscopic level, which in some sense captures the diversity of living systems. 

\section[Modeling diversity in cancer]{Modeling diversity in Cancer \label{Chap6}}
Before discussing about the modeling, we shall first introduce some basic relevant concepts.

\subsection{L\'{e}vy processes and their properties \label{sec:procLevyProp}}

\begin{define}[L\'{e}vy process] A stochastic process $X = (X_t)_{t \ge0}$ is said to be a L\'{e}vy process if:
\begin{enumerate}
\item $X_0$ is almost surely equal to 0.
\item $X$ has independent and stationary distribution, i.e. $X_{t_1} - X_{t_2}$ is independent from $X_{t_2} - X_{t_3}$ for $t_3 \le t_2 \le t_1$ and $\Lu(X_{t_1} - X_{t_2})$ $\sim$ $\Lu(X_{t_1 - t_2} - X_{0})$ for $t_1 \ge t_2$, where $\Lu(X_t)$ denotes the law of the r.v. $X_t$.
\item $X$ is stochastically continuous, i.e. \[ \lim_{t \to s} \P(|X_t - X_s| > a) = 0, \: \forall \: a > 0, \: \forall \: s \ge 0 \]
\end{enumerate}
\end{define}

\begin{define}[\LM]
Let $\nu$ be a Borel measure defined on $\R^d - \{0\}$. Then $\nu$ is said to be a \LM if
\[ \int_{\R^d - \{0\}} (|y|^2 \land 1) \nu(dy) < \infty \]
\end{define}


\begin{thm}[{[\citenum{Applebaum2004}]}] If $X = (X_t)_{t \ge0}$ is a L\'{e}vy process, then its characteristic function $\phi_{X_t} := \E[e^{i(\cdot,X_t)}]$ is given as: \[ \phi_{X_t}(\xi) = e^{-t \psi(\xi)}\, \xi \in \R^d\, t \ge 0\] where $\psi$ is called the L\'{e}vy symbol (or simply symbol) of $X_1$.
\end{thm}


\begin{thm}[L\'{e}vy-Khintchine formula\label{rem:LevyKhin}]
In general, if $X$ is a \LP then its symbol $\psi$ takes the following form 
\[ \psi(\xi) = - i( b, \xi) + {1 \over 2}(\xi, Q \xi) + \int_{\R^d - \{0\}} [1 - e^{i (\xi, y)} + i (\xi, y) \mathbbm{1}_{_{B(0,1)}}(y) ] \nu(dy), \] where $(b, \xi)$ denotes the Euclidean inner product and ${B(0,1)}$ denotes the unit open ball in $\R^d$ centered at the origin. The triple $(b, Q, \nu)$ are called the characteristics of $X$. This is the famous L\'{e}vy-Khintchine formula and can be proved using the well known L\'{e}vy-Ito decomposition.
\end{thm}

\begin{thm}[L\'{e}vy-Khintchine {[\citenum{Schilling2012}]}]\label{thm:nonNegLevyKhin} A function $\psi: \R^d \to \C$ is continuous and negative definite (CNDF) if and only if there exists a number $a \ge 0$, a vector $ b \in \R^d$, a symmetric and positive semi-definite matrix $Q \in \R^{d \times d}$ and a L\'{e}vy measure $\nu$ on $\R^d \backslash \{0\}$ such that
\begin{align}
\label{eq:nonNegLevyKhin}
 \psi(\xi) = a + i (b,\xi) + {1 \over 2} (\xi, Q \xi) + \int_{y \ne 0} \Big( 1 - e^{i (\xi, y)} + i {(\xi, y) \over 1 + |y|^2} \Big) \nu(dy). 
 \end{align}
The quadruple $(a, b, Q, \nu)$ is uniquely determined by $\psi$ and vice versa.
\end{thm}

\noindent More specifically, for a real valued continuous negative definite function (CNDF) we have the following characterization:

\begin{cor}\label{cor:realCNDF} Let $\psi: \R^d \to \R$ be a real valued CNDF. Then it has the following representation
\begin{align}
\label{eq:realCndf}
 \psi(\xi) = a + {1 \over 2} (\xi, Q \xi) + \int_{y \ne 0} \Big( 1 - \cos( (\xi, y) ) \Big) \nu(dy)
\end{align}
where $a$, $Q $ and $\nu$ are as in Theorem \ref{thm:nonNegLevyKhin}.
\end{cor}

\begin{thm}[{[\citenum{Berg1993}]}]\label{thm:compCndfBf} Let  $\psi_1: \R^+ \to \R$ and $\psi_2:\R^d \to \R^+$ be CNDF, then their $\psi_1 \circ \psi_2$ is again a CNDF. 
\end{thm}

\begin{prop}\label{prop:estEta}For every fixed $(b, Q, \nu)$, the characteristics of a \LP $X$, we have that $|\psi(\xi)| \le C (1 + |\xi|^2)$.
\end{prop}
\begin{proof}Let $(b, Q, \nu)$ be the characteristics of a \LP $X$, then the symbol $\psi$, given by the L\'{e}vy-Khintchine formula, is
\[ \psi(\xi) = -i ( b, \xi) + {1 \over 2}(\xi, Q \xi) + \int_{\R^d - \{0\}} [1 - e^{i (\xi, y)} + i (\xi, y) \mathbbm{1}_{_{B(0,1)}}(y) ] \nu(dy). \]
So taking the eucledian norm, we have that
\begin{align*}
|\psi(\xi)|^2 &\le (b,\xi)^2 + {1\over 4}(\xi, Q \xi)^2 + (\int_{\R^d - \{0\}} \cdot )^2 \\
&\le |b|^2|\xi|^2 + \|Q\|^2_{\infty} |\xi|^4 + \int_{\R^d - \{0\}} |e^{i(\xi,y)} -1|^2 \nu(dy)  + \int_{_{B(0,1)}} |(y,\xi)|^2 \nu(dy) \\
&\le |b|^2|\xi|^2 + \|Q\|^2_{\infty} |\xi|^4 + \int_{_{{B^c(0,1)}}} |e^{i(\xi,y)} -1|^2 \nu(dy)  + \int_{_{B(0,1)}} |e^{i(\xi,y)} -1 - i(y,\xi)|^2 \nu(dy) \\
&\le |b|^2|\xi|^2 + \|Q\|^2_{\infty} |\xi|^4 + 4 \nu(_{{B^c(0,1)}}) + \int_{_{B(0,1)}} \Big|(y,\xi)^2 {e^{i(\xi,y)} -1 - i(y,\xi) \over (y,\xi)^2 }\Big|^2 \nu(dy) \\
&\le |b|^2|\xi|^2 + \|Q\|^2_{\infty} |\xi|^4 + 4 \nu(_{{B^c(0,1)}}) + \int_{_{B(0,1)}} (y,\xi)^4 \Big|{e^{i(\xi,y)} -1 - i(y,\xi) \over (y,\xi)^2 } \Big|^2 \nu(dy) \\
&\le |b|^2|\xi|^2 + \|Q\|^2_{\infty} |\xi|^4 + 4 \nu(_{{B^c(0,1)}}) + \int_{_{B(0,1)}} |y|^4|\xi|^4 \Big|{e^{i(\xi,y)} -1 - i(y,\xi) \over (y,\xi)^2 } \Big|^2 \nu(dy) \\
&\le |b|^2|\xi|^2 + \|Q\|^2_{\infty} |\xi|^4 + 4 \nu(_{{B^c(0,1)}}) \\
&\hspace{2cm}+ \int_{_{B(0,1)}} |y|^4|\xi|^4 \Big|{(1- \cos((y,\xi))^2 + ((y,\xi) - \sin((y,\xi)))^2 \over (y,\xi)^2 } \Big|^2 \nu(dy) \\
&\le |b|^2|\xi|^2 + \|Q\|^2_{\infty} |\xi|^4 + 4 \nu(_{{B^c(0,1)}}) + 2 \int_{_{B(0,1)}} |y|^4|\xi|^4  \nu(dy) \\
&\le |b|^2|\xi|^2 + \|Q\|^2_{\infty} |\xi|^4 + 4 \nu(_{{B^c(0,1)}}) + 2 |\xi|^4 \int_{_{B(0,1)}} |y|^2  \nu(dy) \\
&\le |b|^2|\xi|^2 + \|Q\|^2_{\infty} |\xi|^4 + 4 \nu(_{{B^c(0,1)}}) + 2 C_{\nu} |\xi|^4 \\
&\le C_b |\xi|^2 + C_A |\xi|^4 +  C  + 2 C_{\nu} |\xi|^4 \\
&\le C^2 (1 + |\xi|^2)^2 \\
\Ra |\psi(\xi)| &\le C(1 + |\xi|^2).
\end{align*}
\end{proof}
\noindent This provides the basis for characterizing function spaces associated with the symbols $\psi$ of the L\'evy processes.

\subsection{Function spaces associated to negative definite functions \label{sec:ndfFunSp}}
In this section we shall discuss about the function spaces associated with pseudo-differential operators and how this can be made use to characterize the domains of the operators generated by L\'evy processes. Furthermore, we also provide some connections to the Bessel potential spaces $H^s_p(\dom)$. 
\noindent
First let us recall the definition of the Fourier transform and its inverse.\\
Fourier transform $\FT$ is a bijective mapping from $\S(\R^d)$ to $\S(\R^d)$, i.e. $\FT: \S(\R^d) \mapsto \S(\R^d)$. For every $\xi \in \R^d$ the Fourier transform of $f \in \S(\R^d)$ is defined as 
\[(\FT f)(\xi) = \hat f(\xi) = {1 \over (2\pi)^{d/2}} \int_{\R^d} e^{-i (\xi, x)} f(x) dx.\]
For every $x \in \R^d$, the inverse Fourier transform of $f \in \S(\R^d)$ is given by
\[ (\IFT f)(x) = \check f(x) = {1 \over (2\pi)^{d/2}} \int_{\R^d} e^{i (\xi, x)} f(x) dx \]

\begin{define}[Fourier multiplier]
\label{def:multFourier}
A Fourier multiplier $m$ is a mapping $\xi \mapsto m(\xi)$ from $\R^d$ to $\C$, i.e. $m: \R^d \to \C$ where $m(\xi)$ is defined as
\[ \FT [T_m f](\xi) = m(\xi) \hat f, \]
for all $f \in \SW(\R^d)$ and some operator $T_m$ on $\SW(\R^d)$.
\end{define}

\begin{define}[Space of Fourier multipliers  {[\citenum{Farkas2001}]}] \label{def:spFoMul}
Let $p,q \in [1,\infty)$, then the space $\mathcal{M}_{p,q}$ of Fourier multipliers is given as 
\begin{align}
\label{eq:spFoMul}
\mathcal{M}_{p,q} = \{ m \in \SW'(\R^d) : \|m\|_{\mathcal{M}_{p,q}} < \infty \}
\end{align}
where
\begin{align}
\label{eq:spFoMulNorm}
\|m\|_{\mathcal{M}_{p,q}} = \sup \Big\{ {\| \IFT[m \FT \phi] \|_{L^p(\R^d)} \over \|\phi \|_{L^p(\R^d)}}: 0 \ne \phi \in \SW(\R^d) \Big\}.
\end{align}
\end{define}

Now we give a simple result about the composition of two pseudodifferential operators, which will be needed in Section \ref{sec:M2MACP}
\begin{lem}\label{lem:psdOpComp}
Let $T_1$ and $T_2$ be two pseudodifferential operators such that $T_1 f = \IFT[ m_1 \FT f]$ and $T_2 f = \IFT[ m_2 \FT f]$. Then we have that 
\begin{align}
T_3 f(x) := T_2  T_1  f(x) = \IFT[ m_2(x,\xi) m_1(x,\xi) \FT f].
\end{align}
\end{lem} 
\begin{proof}
Let $T_1 f(x) = g(x)$, then $\FT[ T_1 f(x)] = \FT[g(x)]$, then we have that
\begin{align}
\label{eq:psdOpComp}
\FT[T_2 T_1  f(x)] &= \FT[T_2 g(x)] \notag\\
&= m_2(\xi) \FT [g(x)] \notag\\
&= m_2(\xi) \FT [ T_1 f(x)] \notag\\
&= m_2(\xi) m_1(\xi) \FT[f(x)] \notag\\
\Ra T_2 T_1 f(x) &= \IFT[m_2(\xi) m_1(\xi) \FT[f(x)].
\end{align}
\end{proof}

\begin{define}[$\psi$-Bessel potential space of order $s\ge 0$ {[\citenum{Farkas2001}]}]\label{def:cndfBsps}
 Let $\psi: \R^d \to \R$ be a real valued CNDF having the representation  \eqref{eq:realCndf}. Then the $\psi$-Bessel potential space of order $s \ge 0$ with respect to $L^p(\R^d)$ for $p \in [1, \infty)$ is the space
\begin{align}
\label{eq:cndfBsp2}
H^{\psi,s}_p(\R^d) = \big\{ u \in L^p(\R^d) : \|u \|_{H^{\psi,s}_p} < \infty \}
\end{align}
where
\begin{align}
\label{eq:cndfBsp2Norm}
	\|u\|_{H^{\psi,s}_p} := \| (1 + \psi(D))^{s \over 2} u \|_{L^p}.
\end{align}
and 
\begin{align*}
	\psi(D) u &= \IFT(\psi \FT u), \quad u \in \S(\R^d)
\end{align*}
\end{define}
\noindent
It is worth mentioning some relations between classical Bessel potential spaces and $\psi$-Bessel potential spaces:

\begin{rem}[\citenum{Farkas2001}]
\begin{enumerate}
\item For $\psi(\xi) = |\xi|^2$, $H^{\psi,s}_p(\R^d)$ is the classical Bessel potential space $H^s_p(\R^d)$. And for $s \in \N$, it is in turn equivalent to the Sobolev space $W^s_p(\R^d)$. 
\item For $p=2$, we have that $H^{\psi,s}_2(\R^d)$ is a Hilbert space for all $s \ge 0$, and in this case it is simply denoted by $H^{\psi,s}(\R^d)$.
\end{enumerate}
\end{rem}
\noindent

We now state some important embedding results that later help us in specifying continuously varying random operators in Section \ref{sec:conRandOp}.

\begin{thm}[\citenum{Farkas2001}]\label{thm:cndfBspEmb}
Let $\psi: \R^d \to \R$ be a real valued CNDF having the representation  \eqref{eq:realCndf}, $p\in (1,\infty)$ and $r \ge 0$. Then for any $s \in \R$, 
\begin{align}
\label{eq:cndfBspEmb}
	H^{\psi,s+r}_p(\R^d) \xhookrightarrow{} H^{\psi,s}_p(\R^d).
\end{align}
If $s \ge0$, then \eqref{eq:cndfBspEmb} also holds for $p = 1$.
\end{thm}

\begin{thm}[\citenum{Farkas2001}]\label{thm:cndfBspEmb2}
Let $\psi_1, \psi_2: \R^d \to \R$ be two CNDFs. Let $p,q\in [1,\infty)$ and $s, r \in \R$. Then for any $s \in \R$, 
\begin{align}
\label{eq:cndfBspEmb2}
	H^{\psi_1,s}_p(\R^d) \xhookrightarrow{} H^{\psi_2,r}_p(\R^d).
\end{align}
if and only if 
\begin{align}
m:= (1+\psi_2)^{r\over2} (1 + \psi_1)^{-s\over2} \in \mathcal{M}_{p,q}
\end{align}
where, $\mathcal{M}_{p,q}$ is the space of Fourier multipliers of type $(p,q)$ (see Definition \ref{def:spFoMul}).
\end{thm}

\begin{cor}[\citenum{Farkas2001}]\label{cor:cndfBspEmb2}
\begin{enumerate}
\item Let $p \in (1,\infty)$ and let $\psi : \R^d \to \R$ be an unbounded CNDF. Then $H^{\psi,s}_p(\R^d) \xhookrightarrow{} H^{\psi,r}_p(\R^d)$ if and only if $s \ge r$.
\item Let $p \in (1,\infty)$, $s >0$ and  $\psi_1, \psi_2: \R^d \to \R$ be two CNDFs. Then the embedding $H^{\psi_1,s}_p(\R^d) \xhookrightarrow{} H^{\psi_2,s}_p(\R^d)$ implies that there exists a constant $c > 0$ such that 
\[ 1 + \psi_2(\xi) \le c (1 + \psi_1(\xi))\, \xi \in \R^d. \]
If, in particular, $H^{\psi_1,s}_p(\R^d) = H^{\psi_2,s}_p(\R^d)$, then we have 
\[  {1 \over c} (1 + \psi_1(\xi)) \le 1 + \psi_2(\xi) \le c (1 + \psi_1(\xi))\, \xi \in \R^d. \]
More interestingly, the converse relation holds for $p = 2$.
\end{enumerate}
\end{cor}

\begin{thm}[Operator domain associated to $H^{\psi,s}_p(\R^d)$ [\citenum{Farkas2001}]] \label{thm:opDomEtaBsp}
For $s \ge 0$, $p \in [1, \infty)$ and $\psi: \R^d \to \R$ a CNDF, we have $H^{\psi,s}_p(\R^d) = D((1 - \psi(D))^{s\over2})$, where
\[ (1 - \psi(D))^{s\over2} u = \IFT[ (1 + \psi(\xi))^{s\over2} \FT u] , \quad \forall u \in \SW(\R^d). \]
\end{thm}

\subsection{From a microscale SDE to a macroscale abstract Cauchy problem \label{sec:M2MACP}}
Consider the following simple SDE driven by a L\'{e}vy process $N_t$ on $\R^d$ which has $A$ as the IG of the corresponding $C_0$-semigroup (in short as $C_0$-SG) $T_t$.
\begin{align*}
	dY_t &= dN_t \\
	Y_0	&= 0.
\end{align*}
Then clearly, $Y_t = N_t$. Now the $C_0$ SG $T_t$ on a Banach space $X$ associated with $N_t$ can be prescribed by $ T_t f(y) = \E(f(Y_t)| Y_0 = y)$ for every $f \in X$. Since $A$ is the IG of the $C_0$ SG $T_t$, for all $f \in D(A)$ we have that 
\begin{align*}
{d \over dt} T_t f &= A T_t f \\
\Rightarrow {d \over dt} \E( f(Y_t) | Y_0 = y) &= A \E(f(Y_t) | Y_0 = y)\\
\Rightarrow {d \over dt} \int_{\R^d} f(x) p_t(y,dx) &= A \int_{\R^d} f(x) p_t(y,dx).
\end{align*}
Now let $p_t(\cdot, dx) \ll dx$, then 
\begin{align*}
{d \over dt} \int_{\R^d} f(x) p_t(y,x) dx &= A \int_{\R^d} f(x) p_t(y,x) dx
\Rightarrow \int_{\R^d} f(x) \partial_t p_t(y,x) dx = \int_{\R^d} f(x) A p_t(y,x) dx. \\
\end{align*}
If this holds for all $f \in C^{\infty}_c(\R^d)$ then we get that
\[ \partial_t p_t(y, x) = A p_t(y, x), \text{ for a.e. } x \in \R^d. \]
This is nothing but the backward Kolmogorov's equation. \\

Similarly, using the identity $ {d \over dt} T_t f = T_t A f$ we get
\begin{align*}
{d \over dt} T_t f &= T_t A f \\
\Rightarrow {d \over dt} \E( f(Y_t) | Y_0 = y) &= \E( A f(Y_t) | Y_0 = y)\\
\Rightarrow {d \over dt} \int_{\R^d} f(x) p_t(y,dx) &= \int_{\R^d} A f(x) p_t(y,dx)
\end{align*}
Now let $p_t(\cdot, dx) \ll dx$, then 
\begin{align*}
{d \over dt} \int_{\R^d} f(x) p_t(y,x) dx &= \int_{\R^d} A f(x) p_t(y,x) dx
\Rightarrow \int_{\R^d} f(x) \partial_t p_t(y,x) dx = \int_{\R^d} f(x) A^{\dagger} p_t(y,x) dx \\
\end{align*}
where $A^{\dagger}$ is the formal adjoint of the operator $A$. If this holds for all $f \in C^{\infty}_c(\R^d)$ then we get that
\[ \partial_t p_t(y, x) = A^{\dagger} p_t(y, x), \text{ for a.e. } y \in \R^d. \] This is called the forward Kolmogorov equation.

\begin{eg}\label{eg:Op1} Let $N_t$ be a Brownian motion with drift on $\R^d$. It is a L\'{e}vy process with characteristics $(b, a, 0)$. The corresponding IG $A$ of the $C_0$ SG $T_t$ associated to $N_t$ is given by \[ A = \sum_k^d b_k \partial_k + {1 \over 2} \sum_{k,j}^d a^j_k \partial_j \partial_k. \]
The formal adjoint $A^{\dagger}$ takes the form
\[ A^{\dagger} = -\sum_k^d  \partial_k  b_k + {1 \over 2} \sum_{k,j}^d \partial_j \partial_k a^j_k  . \]
\end{eg}

\begin{eg}\label{eg:Op2} Let $N_t$ be a Poisson process on $\R^d$ with characteristics $(0, 0 , \lambda \delta_1)$. The corresponding IG $A$ of the $C_0$ SG $T_t$ associated to $N_t$ is given by \[ A f \:(x) = \lambda [f(x+1) - f(x)]. \]
The formal adjoint $A^{\dagger}$ takes the form
\[ A^{\dagger}f \: (x) = \lambda [f(x) - f(x -1)]. \]
\end{eg}

\begin{eg}\label{eg:Op3} Let $N_t$ be a compound Poisson process on $\R^d$ with characteristics $(0, 0 , \lambda \mu)$. The corresponding IG $A$ of the $C_0$ SG $T_t$ associated to $N_t$ is given by \[ A f \:(x) = \lambda  \int_{\R^d}  [f(x+y) - f(x)] \mu(dy). \]
The formal adjoint $A^{\dagger}$ takes the form
\[ A^{\dagger}f \: (x) = \lambda \int_{\R^d} [f(x) - f(x -y)] \mu(dy). \]
\end{eg}

\begin{eg}\label{eg:Op4} Let $N_t$ be a L\'{e}vy process on $\R^d$ with characteristics $(b, Q, \lambda \nu)$. The corresponding IG $A$ of the $C_0$ SG $T_t$ associated to $N_t$ is given by \[ Af \:(x) = \sum_k b_k  \partial_k  f + {1 \over 2} \sum_{k, j} Q_{kj} \partial_k \partial_j f + \lam\int_{\R^d} [ f(x+y) - f(x) \nu(dy). \]
The formal adjoint $A^{\dagger}$ takes the form
\[ A^{\dagger}f \:(x) = -\sum_k \partial_k b_k   f + {1 \over 2} \sum_{k, j} \partial_k \partial_j Q_{kj} f + \lam \int_{\R^d} [ f(x) - f(x-y) \nu(dy). \]
\end{eg}
\begin{eg} \label{eg:alphaStable} Let $N_t$ be a symmetrically invariant stable process of index $\alpha \in (0,2)$ with symbol given by $\psi(\xi) = - |\xi|^{\alpha}$ for all $\xi \in \R^d$. The corresponding IG $A$ of the $C_0$ SG $T_t$ associated to $N_t$ is given by
\[ A = - \Big(\sum_k^d - (\partial_k)^2\Big)^{\alpha/2} = -(-\lap)^{\alpha/2}. \]
The formal adjoint $A^{\dagger}$ is the same as the operator $A$. In fact we actually have that $A$ is a self-adjoint operator.
\end{eg}
\noindent
Above we used a very simple SDE to describe the microscopic dynamics. In practical applications, however, one encounters more generic equations of the following form:
\begin{align}
\label{eq:LevyIntg}
	X_t &= X_0 + \int_0^t b(X_s) ds + \int_0^t q(X_s) dW_s + \int_0^t \int_K g(s,x) \tilde N(ds\, dx) + \int_0 \int_{K^c} h(s,x) N(ds, dx) \notag\\
	&=: X_0 + \hat X + \int_0^t \int_K g(s,x) \tilde N(ds\, dx) + \int_0 \int_{K^c} h(s,x) N(ds, dx) \\
	dX_t &= b(X_t) dt + q(X_t) dW_t + \int_K g(t,x) \tilde N(dt\, dx) + \int_{K^c} h(t,x) N(dt\, dx)\notag \\
	 &= d\hat X + \int_K g(t,x) \tilde N(dt\, dx) + \int_{K^c} h(t,x) N(dt\, dx)\notag 
\end{align}
where, $X_t \in \R^d$, $q : \R^d \in \R^d$, $b \in \R^d \in \R^{d \times n}$, $W_t \in \R^n$, $g \in \R^+ \times \R^d \to \R^d$, $h \in \R^+ \times \R^d \to \R^d$ and $N_t(K) \in \R^d$ for all $K \in \R^d \backslash \{0\}$. \\
Since the above SDE is time homogeneous, the solution $X_t$ is a time homogeneous Markov process. Thus we get a $\CSG$ via the prescription $\E(f(X_t)|X_0 = x)$ (details to be provided later). In order to be able to find the associated IG, we need Ito's formula for L\'{e}vy type stochastic integrals.

\begin{thm}[Ito-type formula for L\'evy integrals [\citenum{Applebaum2004}]]\label{thm:ItoLevy} Let $X_t$ be a L\'{e}vy type stochastic integral taking the form \eqref{eq:LevyIntg}. Then for each $f \in C^2(\R^d)$, $t \ge 0$, with probability 1 we have that
\begin{align}
\label{eq:ItoLevy}
f(X_t) - f(X_0) &= \int_0^t \sum_{k=1}^d \partial_k f(X_s) d \hat X_s + {1 \over 2}\int_0^t \sum_{j,k=1}^d  \partial_k \partial_j f(X_s) d[\hat X^k,\hat X^j]_s \notag\\
& \hspace{-1.5cm}+ \int_0^t \int_K f(X_s^- + g(s,x)) - f(X_s^-) \tilde N(ds,dx) + \int_0^t \int_{K^c} f(X_s^- + h(s,x)) - f(X_s^-) N(ds\,dx)\notag\\
&+ \int_0^t \int_{K} [f(X_s^- + g(s,x)) - f(X_s^-) - g(s,x) \sum_{k=1}^d \partial_k f(X_s) ] \nu(dx) ds,
\end{align}
where $X_s^{-} = X_{s-}$ is the left limit point of $X$ at $s$. 
\end{thm}
\noindent
Now using the above theorem we would like to determine the generator of a generic \LP as given by \eqref{eq:LevyIntg}.

\begin{prop}\label{prop:LevyGen} Let $T_t f = \E(f (X_t) | X_0 = x)$ for $f \in \Lp$, $1 \le p < \infty$ be the $\Lp$-SG associated with the \LP $X$. Then the IG $A$ associated to $(T_t)_{t \ge 0}$ is given by

\begin{align*}
A f (X_s) &= \sum_{k=1}^d \partial_k b_k(X_s) f(X_s) + {1 \over 2} \sum_{k,j=1}^d  Q_{kj}(X_s) \partial_k \partial_j f(X_s)  + \int_{K^c} f(X_s^- + h(s,x)) - f(X_s^-) \nu(dx) \\
&\hspace{2cm}+ \int_{K} [f(X_s^- + g(s,x)) - f(X_s^-) - g(s,x) \sum_k^d \partial_k f(X_s) ] \nu(dx),
\end{align*}
for $f \in C^2(\R^d) \cap L^p(\R^d)$.
\end{prop}
\begin{proof} 
First let $f \in C^2(\R^d) \cap L^p(\R^d)$. Then,
\begin{align*}
T_{t+h}f - T_t f &= \E\Big(f(X_{t+h}) - f(X_t) | X_0 =x \Big) \\
&\overset{\eqref{eq:ItoLevy}}{=}  \E \Bigg(\Big[\int_t^{t+h} \sum_{k=1}^d \partial_k f(X_s) dX_s + {1 \over 2}\int_t^{t+h} \sum_{k,j=1}^d  \partial_k \partial_j f(X_s) d[X^k,X^j]_s \\
&\hspace{.1cm}+ \int_t^{t+h} \int_K f(X_s^- + g(s,x)) - f(X_s^-) \tilde N(ds,dx) \\
&\hspace{.1cm}+ \int_t^{t+h} \int_{K^c} f(X_s^- + h(s,x)) - f(X_s^-) N(ds\,dx)\\
&\hspace{.1cm}+ \int_t^{t+h} \int_{K} [f(X_s^- + g(s,x)) - f(X_s^-) - g(s,x) \sum_k^d \partial_k f(X_s) ] \nu(dx) ds\Big] \Big | X_0 = x \Bigg)
\end{align*}
\begin{align*}
&\overset{\text{Martingale}}{\underset{\text{property}}{=}}  \E \Bigg( \Big[\int_t^{t+h} \sum_{k=1}^d b_k(X_s) \partial_k f(X_s) ds + {1 \over 2}\int_t^{t+h} \sum_{k,j=1}^d  Q_{kj}(X_s) \partial_k \partial_j f(X_s) ds \\
&\hspace{.1cm} + \int_t^{t+h} \int_{K^c} f(X_s^- + h(s,x)) - f(X_s^-) \nu(dx) ds\\
&\hspace{.1cm}+ \int_t^{t+h} \int_{K} [f(X_s^- + g(s,x)) - f(X_s^-) \nu(dx) ds  - \int_t^{t+h} g(s,x) \sum_k^d \partial_k f(X_s) ] \nu(dx) ds\Big] \Big | X_0 = x \Bigg)
\end{align*}
\begin{align*}
&= \int_t^{t+h} \Bigg[ \E \Big(\sum_{k=1}^d b_k(X_s) \partial_k f(X_s) + {1 \over 2} \sum_{k,j=1}^d   Q_{kj}(X_s)\partial_k \partial_j f(X_s)  + \int_{K^c} f(X_s^- + h(s,x)) - f(X_s^-) \nu(dx) \\
&\hspace{1cm} + \int_{K} [f(X_s^- + g(s,x)) - f(X_s^-) - g(s,x) \sum_k^d \partial_k f(X_s) ] \nu(dx) \Big] \Big | X_0 = x \Big)\Bigg] \: ds \\
&= \int_t^{t+h} \E( \hat A f(X_s) | X_0 = x) ds \quad = \int_t^{t+h}T_s \hat A f \: ds 
\end{align*}
where, 
\begin{align*} 
\hat A f (x) &= \sum_{k=1}^d \partial_k b_k(x) f(x) + {1 \over 2} \sum_{k,j=1}^d  Q_{kj}(x) \partial_k \partial_j f(x)  + \int_{K^c} f(x^- + h(s,x)) - f(x^-) \nu(dx) \\
&\hspace{2cm}+ \int_{K} [f(x^- + g(s,x)) - f(x^-) - g(s,x) \sum_k^d \partial_k f(x) ] \nu(dx)
\end{align*}
and $\Big(Q_{kj}(x)\Big)_{k,j \in \{1, \dots, d\}} := a(x) a(x)^T$. 

Therefore, if $A$ is the generator of the $\Lp$-SG $T_t$, then
\begin{align*}
\int_t^{t+h} T_s A f \: ds &= \int_t^{t+h}T_s \hat A f \: ds \\
\Ra A f &= \hat A f ,\quad \text{ for all } f \in C^2(\R^d) \cap \Lp.
\end{align*}
Since $C^{\infty}_c(\R^d) \subset \Lp$ is dense and since $D(A) \subset \Lp$, in particular it holds that $A f = \hat A f$, for $f \in D(A)$ such that $f_n \overset{n \to \infty}{\to} f$ in $\Lp$ with $f_n \subset C^{\infty}_c$. Indeed, since $A$ ( the generator of a $\Csg$) is closed, $A f_n \to A f$ in $\Lp$ as $n \to \infty$. But, since $\hat A f_n = A f_n$ we get that the limits coincide. Here $A$ can be seen as the smallest closed extension of $\hat A$.
\end{proof}

\begin{lem}\label{lem:ndfClosed} Let $\psi$ be a negative definite function. Then the pseudo-differential operator defined as 
\[ A f = \IFT[ \psi \FT f] , \quad \text{ for all $f \in \S(\R^d)$} \] is closable in $L^2(\R^d)$. 
\end{lem}
\begin{proof} Let $(f_n)_{n \in \N}$ be a sequence in $\S(\R^d)$ such that $f_n \to 0$ in $L^2(\R^d)$, then by Plancherel we have that $\hat f_n \to 0$ in $L^2(\R^d)$, in particular $\hat f_{n_k} \to 0$ a.e. along some subsequence $\hat f_{n_k}$. Now, let $ A f_n \to g$ in $L^2(\R^d)$ then by Plancherel we get that $\widehat{A f_n} \to \hat g$ in $L^2(\R^d)$. But the definition of $A $  implies that $ \psi \hat f_n \to \hat g$ in $L^2(\R^d)$. So by dominated convergence we get that $\psi \hat f_n \to 0$ in $L^2(\R^d)$ i.e. $g = 0$ and hence $A$ is closable. 
\end{proof}
\noindent
Thus we see that depending on the type of noise acting at the microscale we get different operators at the macroscale. \\

\subsection[Examples of transport equations]{Examples of transport equations  as a particularization of the above generators}
Many transport equations that appear in the context of biology takes the following form 
\begin{align}
\label{eq:trans}
\begin{split}
\partial_t p_t(x,v) &= - v \cdot \grad p_t(x,v) + \int_{V} [K(v,v') p_t(x,v') - K(v,v') p_t(x,v) ] dv' \\
\int_V K(v,v') dv' &= 1\, K(v,v') = K(v',v) \: \text{ for all }  v,v' \in V \subset \R^3, 
\end{split}
\end{align}
where $p_t(x,v)$ is the density of the required biological quantity being studied (e.g. population density, chemical concentration, etc.) depending on the time $t$, the space variable $x \in \R^3$ and the velocity variable $v \in V \subset \R^3$. Usually, the velocity state space $V$ is taken to be a finite state space having the form $[v_1, v_2]\times \mathbb{S}^2$, where $[v_1, v_2]$ is a finite interval for the magnitude of the velocity vector $v$ and $\mathbb{S}^2$ is the two dimensional sphere representing the orientation of the velocity vector $v$. Here we just consider a generic finite state space $V \subset \R^3$.
To arrive at the transport equation \eqref{eq:trans} we employ the following L\'{e}vy-type SDE:
\begin{align}
\label{eq:motEqn}
\begin{split}
dX_t = V_t dt, \\
dV_t = dL_t
\end{split}
\end{align}
where, $X_t$ and $V_t$ are the position and velocity, respectively, of an individual (e.g. cell) at time $t$, and $X_t$ is a jump type L\'{e}vy process. In order to arrive at the particular form of equation given by \eqref{eq:trans} we assume that $L_t$ is a compound Poisson process with intensity $\lambda $ and symmetric jump probability measure $\mu$ on $\R^d$. Using L\'{e}vy-It\^{o} decomposition we can write $X_t$ as 
\begin{align*}
	L_t &= c t + \int_{|v| < 1} v \tilde N(t,dv) + \int_{|v| \ge 1} v N(t, dv)
\end{align*}
Now using It\^{o}'s formula (see Theorem \ref{thm:ItoLevy}) and Proposition \ref{prop:LevyGen}, for $b(x, v) = (v,0)^T$, $a(x, v) = (0,0)$, $g(x,v) = v$, and $h(x,v) = v$ we get that
\begin{align}
\label{eq:genTransEqn}
{d \over dt} p_t(x,v) = -v \cdot \partial_x p_t(x,v) + \lambda \int_{\R^d} [p_t(x,v) - p_t(x,v - v')] \mu(dv')
\end{align}
Now if $\mu(dv) \ll dv$ i.e. $\mu$  is dense with respect Lebesgue measure such that ${d\mu \over dv}(dv) = f(v)dv$ we get that 
\begin{align*}
{d \over dt} p_t(x,v) &= -v \cdot \partial_x p_t(x,v) + \lambda \int_{\R^d}[ p_t(x,v) - p_t(x,v - v') ] f(v') dv' \\
&= -v \cdot \partial_x p_t(x,v) + \lambda \int_{\R^d}[ p_t(x,v') - p_t(x,v)]  f(v- v') dv'.
\end{align*}
Now if $\mu$ has a finite support $V \subset \R^d$ then the above equation reduces to 
\begin{align*}
{d \over dt} p_t(x,v) &= -v \cdot \partial_x p_t(x,v) + \lambda \int_{V}[ p_t(x,v') - p_t(x,v)]  f(v-v') dv'
\end{align*}
This can now be written in the form of \eqref{eq:trans} by letting $K(v',v) = f(v-v')$. \\
So \eqref{eq:genTransEqn} can be seen as a general form of the transport equation (at least for the biological applications). Based on the applications the jump measure $\mu$ can be particularized in the following different forms:
\begin{enumerate}
\item $\mu(dv) = f(v) dv$. This case is feasible when the velocity of individuals is an observed (say by cell  or particle tracking in images and/or videos) quantity and its probability or density can be estimated.
\item $\mu(dv) = \Big( \int_{\Theta} f(v,\theta) d\theta \Big) dv  = \Big(\int_{\Theta} f(v|\theta) g(\theta) d \theta \Big) dv$. This case is feasible when the velocity is dependent on an additional observable variable $\theta$ (say orientation data, e.g. fiber orientation, track orientation) which can be estimated. 
\end{enumerate}
\noindent Similarly, the microscopic equation 
\begin{align*}
\label{eq:diffmotEqn}
\begin{split}
dX_t &= V_t dt + \sigma(X_t) dW_t, \\
dV_t &= dL_t
\end{split}
\end{align*}
results in the following macroscopic model:
\begin{align*}
{d \over dt} p_t(x,v) &= -v \cdot \partial_x p_t(x,v) + {1 \over 2} \sum_{j,k=1}^3   \partial_{x_k} \partial_{x_j} \sigma_k^j(x) p_t(x,v) + \lambda \int_{V}[ p_t(x,v') - p_t(x,v)]  f(v-v') dv'
\end{align*}
\noindent To obtain reaction terms, one needs to consider functionals of the process $Y_t = (X_t,V_t)$. Let $u(t,y)$ be defined as $$u(t,y) := \E \Big[ \varphi(Y_t) e^{\int_0^{t+h} f(s,Y_s)ds } \Big| Y_0 = y \Big].$$ Now a macroscopic equation associated to $u$ can be obtained in the following way:
\begin{align*}
    u(t+h,y) &= \E \Big[ \varphi(Y_{t+h}) e^{\int_0^{t+h} f(s,Y_s)ds } \Big] \Big| Y_0 = y \Big] \\
    &=\E \Big[ \E \Big[ \varphi(Y_{t+h}) e^{\int_0^{t+h} f(s,Y_s) ds} \Big] \Big| \F_h) | Y_0 = y \Big] \\
    &=\E \Big[ \E \Big[ \varphi(Y_{t+h}) e^{\int_h^{t+h} f(s,Y_s) ds} \Big] \Big| \F_h \Big] e^{\int_0^{h} f(s,Y_s) ds} \Big| Y_0 = y \Big] \\
    &=\E \Big[ u(t,Y_h) e^{\int_0^{h} f(s,Y_s) ds} \Big| Y_0 = y \Big] \\
\Ra u(t+h,y) - u(t,y) &=\E \Big[ u(t,Y_h) e^{\int_0^{h} f(s,Y_s) ds} - u(t,y)\Big| Y_0 = y \Big] \\
\Ra \lim_{h \to 0} {u(t+h,y) - u(t,y) \over h} &= \lim_{h \to 0} \E \Big[ {1 \over h}\Big(u(t,Y_h) e^{\int_0^{h} f(s,Y_s) ds} - u(t,y) \Big) \Big| Y_0 = y \Big] \\
&= \lim_{h \to 0} \E \Big[ e^{\int_0^{h} f(s,Y_s) ds} \partial_y u(t,Y_h) + u(t,Y_h) \partial_t e^{\int_0^{h} f(s,Y_s) ds} \\
&\hspace*{2cm} + \mathcal{O}(h^2)\Big) \Big| Y_0 = y \Big] \\
\Ra \partial_t u(t,y) &= A u(t,y) + u(t,y) f(s,y), \quad u(0,y) = \varphi(y).
\end{align*}
\noindent Analogously, letting $u(t,y) := \E \Big[ \varphi(Y_t) e^{\int_0^{t+h} f(s,Y_s)ds } + \int_0^t g(s,Y_s) ds \Big| Y_0 = y \Big]$, yields the following macroscopic equation:
\begin{align*}
\partial_t u(t,y) &= A u(t,y) + u(t,y) f(t,y) + g(t,y), \quad  u(0,y) = \varphi(y).
\end{align*}
Going a step further, taking $f$ and $g$ to be functions of $u$ would yield nonlinear macroscopic equation-
\begin{align*}
\partial_t u(t,y) &= A u(t,y) + u(t,y) f(u) + g(u), \quad  u(0,y) = \varphi(y).
\end{align*}
Thus we see that, macroscopic equations arise a result of averaging out the dynamics of a functional of the microscopic dynamics. If the macroscopic functional itself is random then this gives rise to random macroscoipc equations. We shall discuss this in more detail in the next section. 

\subsection{Construction of random operators \label{sec:conRandOp}}
In this section we provide a framework for the construction of random time dependent operators. The most simplest idea is to use an operator valued process $(A_t)_{t\in (0,T)}$ which, for each fixed $\om \in \Om$, is constant in a given time interval $(t_k,t_{k+1}) \subset (0,T)$ for any $k \in \N$. More precisely, we have the following: \\
Let $(\tau_k)_{k\in\N}$ be a  sequence of exponentially distributed independent random variables. Let $N_t$ be the standard Poisson process with intensity $\lam$, and let $t_n = \sum_{k=0}^{N_t=n} \tau_k$ with $\tau_0 = 0$. Here $\tau_k$ represents the waiting time and $t_n$ represents the jumping time, for each $k,n \in \N$, respectively. Then we can define a piecewise constant (w.r.t time) operator $A_t$ as
\begin{align}
A_t(\om) = A_n \quad \text{for $t \in (t_n(\om), t_{n+1}(\om)) \subset (0,T]$}
\end{align}
where $(A_n)_{n \le N}$, $N \in \N$ is a sequence of operators, each generated by a different \LP and its associated $L^p$ semigroup. Based on this, for every fixed $\om \in \Om$ and a given realization of jumping times $(t_n)_{n\in\N}$ we get the following random initial value problem:
\begin{align}
{d \over dt} T^n_t f &= A_n T^n_t f \quad \text{for $t \in (t_n, t_{n+1})$} \\
T^n_0 f &= \lim_{t \uparrow t_{n-1}} T^{n-1}_t f \\
T^0_0 f &= f.
\end{align}
Now the idea is to generalize this concept for continuously varying operators. To this end we use the relation between characteristic functionals, negative definite functions, and the associated generators to construct random operators. Here we assume that we are given the L\'evy symbol $\psi$ (which in general may be random and time dependent as well) and from this we define the generator via Fourier transform which will be both random and time dependent.\\

\newcommand{\Lv}{\mathfrak{L}}
\newcommand{\Ag}{\mathfrak{A}}

\begin{define}\label{def:nonAtndfOp} Let $\Theta = (\Theta_t)_{t \in [0,T]}$ be a  stochastic process taking values in $CN(\Xi)$, the convex cone of continuous negative definite functions on $\Xi$, where $\Xi = \R^d$. Then for such $\Theta$ we define an operator valued process $\Ag = (A_t)_{t \in [0,T]}$ as  
\begin{align}
\label{eq:nonAtndfOp}
 A_t f = \IFT[ \Theta_t \FT f], \quad \text{ for each $f \in \S(\R^d)$ and $t \ge 0$} \end{align}
\end{define}
\noindent
Here we can apply Theorem \ref{thm:nonNegLevyKhin} and represent $\Theta_t$ as  $(a_t, b_t, Q_t, \nu_t)$ for each $t \in [0,T]$, respectively. The elements are such that $ a_t \in \R^+$, $b_t \times \R$,  $Q_t \in \R^{d\times d}$ and $\nu_t \in \Lv(\Xi)$, where $\Lv(\Xi)$ is the space of L\'evy measures on $\Xi$. In the light of Lemma \ref{lem:ndfClosed}, the operator family $\Ag$ defined via \eqref{eq:nonAtndfOp} consists of a set of closable operators $A_t$ on $L^2(\R)$ for each $t \in [0,T]$. In order to be able to deduce a differential equation associated to the time dependent operators defined above, we need to construct, for each $\om \in \Om$, an evolution operator $\evolOp$ associated with the process $\Ag$. To this end we make the following ansatz:\\[-5ex]
\begin{ansatz}
\label{asm:nonAtndfOp}
\begin{enumerate}
\item $\Ag = (A_t)_{t\in(0,T]}$ is a family of sectorial operators for each $\om \in \Omega$. 
\item The family $\Ag$ of sectorial operators satisfies  the following Assumption \ref{asm:nonAtAcp} for each $\om \in \Omega$. \\[-5ex]
\end{enumerate}
\end{ansatz}
\begin{assmp}
\label{asm:nonAtAcp}
\begin{enumerate}
\item Let $A_t \in \secOpSp(X)$ have an uniform sector angle $\secAng{A}$ for all $t \in [0,T]$, i.e.
\begin{align}
\label{eq:asmOpPos}
\begin{split}
\sigma(A_t) \subset \Sigma_{\secAng{A}} &= \Big\{\lam \in \C \backslash \{0\}: |\arg \lam| < \secAng{A} \Big\}, \quad t \in [0,T],\\
\|(\lam - A_t)^{-1} \|_{_{L(X)}} &\le { M \over |\lam|}, \quad \lam \notin \Sigma_{\secAng{A}}, t \in [0,T].
\end{split}
\end{align}
\item There exists an exponent $\nu >0$ such that 
\begin{align}
\label{eq:asmFracDom}
D(A_s) \subset D(A_t^{\nu}), \quad \forall t,s \in [0,T]
\end{align}
\item In addition $A_t^{-1}$ is H\"{o}lder continuous with respect to $t$, i.e. 
\begin{align}
\label{eq:asmFracResLip}
\|A_t^{\nu} [ A_t^{-1} - A_s^{-1}]\|_{_{L(X)}} \le \const_A |t-s|^{\mu}, \quad s,t \in [0,T]
\end{align}
for some constant $\const_A > 0$ and some exponent $\mu \in (0,1)$ such that $1 < \mu + \nu$.
\end{enumerate}   
\end{assmp}
\begin{rem} Equation \eqref{eq:asmOpPos} of Assumption \ref{asm:nonAtAcp} implicitly assumes that $0 \in \rho(A_t)$ for all $t \in [0,T]$, i.e. $(A_t)_{t \in [0,T]} \in P(X)$ i.e. $(A_t)_{t \in [0,T]}$ is a family of positive operators on $X$. 
\end{rem}

\begin{thm}\label{thm:evolOpSecOp} Let $(A_t)_{t \in [0,T]}$ be a family of sectorial operators satisfying the Assumption \ref{asm:nonAtAcp}, then there exists a unique two parameter evolution operator $\evolOp$ having the properties specified by Properties \ref{prop:evolOp} below. 
\end{thm}
\begin{proof} Refer to [\citenum{PAZY83, YAGI09}]
\end{proof}

\begin{ppts*}[Evolution operator] \label{prop:evolOp} Let $\delInt = \{(t,s) \in (\R^+)^2 : 0 \le s < t \le T \}$, then $U(t,s)$ satisfies the following properties:
\begin{enumerate}
\item Semigroup property: 
\begin{align}
\label{eq:propSgEvolOp}
\begin{split}
U(t,s) &= U(t,r)U(r,s), \quad \text{ for all $0 \le s \le r \le t \le T$}, \\
U(t,t) &= 1, \quad \text{ for all $t \in [0,T]$}.
\end{split}
\end{align}
\item The mapping $(t,s) \mapsto U(t,s)$ belongs to $C(\overline{\delInt};L(X))$ and satisfies the estimate
\begin{align}
\label{eq:evolOpUniBd}
\|U(t,s)\|_{_{L(X)}} \le \arbconst , \quad \text{for some arbitrary constant $\arbconst >0$}.
\end{align}
\item The mapping $(t,s) \mapsto A_tU(t,s)$ belongs to $C(\delInt; L(X))$  and satisfies the estimate
\begin{align}
\label{eq:secOpevolOpSingBd}
\|A_t U(t,s)\|_{_{L(X)}} \le \arbconst (t-s)^{-1}, \quad (t,s) \in \delInt.
\end{align}
\item The mapping $(t,s) \mapsto U(t,s)$ satisfies the following initial value problem:
\begin{align}
\label{eq:evolOpEqn1}
\begin{split}
    {\partial \over \partial t} U(t,s) &= - A(t) \: U(t,s) , \quad s \le t \le T,\\
    U(s,s) &= 1,
\end{split}
\end{align}
and
\begin{align}
\label{eq:evolOpEqn2}
\begin{split}
    {\partial \over \partial t} U(t,s) &= U(t,s) \:  A(t),  \quad s \le t \le T,\\
    U(s,s) &= 1.
\end{split}
\end{align}
\end{enumerate}
\end{ppts*}

\noindent Based on this we can construct a two parameter evolution operator $ \Ug := \evolOp$ on $L^2(\R^d)$ associated to $\Ag$, for every $\om \in \Omega$ fixed. We shall now illustrate this process.\\
\noindent
Let $A$ be a sectorial operator of angle $\secAng{A}$ and let 
$\psi_{A}$ be a real valued CNDF associated with the operator $A$. Let $\psi_S : \R^+ \to \R$ be another CNDF, then  due to Theorem \ref{thm:compCndfBf}, we can compose $\psi_S$ with $\psi_A$ to get another real valued CNDF $\psi_Z$ which is given as $\psi_Z := \psi_S \circ \psi_A$. Now applying the sector mapping theorem ([\citenum{Berg1993}]) we get that $\psi_{Z}$ is also a sectorial operator of angle $\secAng{\psi_{_Z}} \le \secAng{A}$. Thus one may compose different CNDF to generate new sectorial operators. Now in order to get time dependent random operators, we let $\psi_t = \beta_t \psi_Z$, where $\beta_t$ is a real valued bounded continuous stochastic process. Letting $\beta_t \subset [\beta_1, \beta_2]$ for $0 < \beta_1 < \beta_2$, we see that $(1 + \beta_1 \psi_Z) \le (1 + \beta_2 \psi_Z)$. Thus applying the embedding result in Corollary \ref{cor:cndfBspEmb2}, we get that $H^{\beta_2 \psi_Z, r}(\R^d) \xhookrightarrow{} H^{\beta_1\psi_Z,r}(\R^d)$ for all $r \in \R$. Now based on Theorem \ref{thm:opDomEtaBsp}, for each $t \in [0,T]$ and $\om \in \Om$ we can associate an operator $A_t = (1 - \psi_t(D))^{s\over2}$, such that $D(A_t) = H^{\psi_t,s}(\R^d)$. Letting $\Theta_t = (1 + \psi_t )^{s \over 2}$ we get the desired continuous negative definite function. Note that since $\psi_t(A)$ is a sectorial operator of angle $\secAng{\psi_t} \le \secAng{A}$ we also get that $A_t$ is a sectorial operator of angle $\secAng{A_t} \le {s \over 2} \secAng{\psi_t}$. 
These considerations can be summarized in the following lemma.

\begin{lem}\label{lem:randOpLem} Let $(\Theta_t)_{t \in [0,T]}$ be defined as $(1+\psi_t)^{s\over 2}$ with $\psi_t = \beta_t \psi$ where $\beta_t \in [\beta_1, \beta_2]$, $0 < \beta_1 < \beta_2$ is a bounded stochastic process and $\psi$ is a real-valued CNDF generated by a L\'evy process. Correspondingly, let $(A_t)_{t\in[0,T]}$ be the family of differential operators given by $A_t = (1 - \psi_t(D))^{s \over 2}$. Then we have that:
\begin{enumerate}
\item $\Theta_t = (1 + \psi_t )^{s \over 2}$ is a stochastic process such that  for every fixed $t \in [0,T]$ and $\om \in \Omega$, $\Theta_t(\om)$ is a real valued CNDF.
\item The operator $A_t$ generated by $\Theta_t$ is a sectorial operator of angle $\secAng{A_t} \le {s \over 2} \secAng{\psi_t} \le \piT$.
\item $D(A_{t_1}) = H^{\psi_{t_1}, s}(\R^d) \xhookrightarrow{} H^{\beta_1 \psi_Z,s}(\R^d) \xhookrightarrow{} H^{\beta_1 \psi_Z,r}(\R^d) = D((A_{t_2})^{\nu})$ for  $s > r > 0$, $\nu = {r \over s} < 1$, and all $t_1, t_2 \in [0,T]$. 
\end{enumerate}
\end{lem}
\noindent
The only thing that needs to be checked is that $\|A_t^{\nu}[A_t^{-1} - A_s^{-1}]\| \le C |t-s|^{\mu}$ with $1 < \mu + \nu$, which we shall establish now.\\
\noindent
Since $A_t u = \IFT[ \Theta_t \FT u ]$ for all $u \in \SW(\R^d)$, we see that
$ \IFT [\Theta_t^{-1} \FT A_t u ] = u$. By letting $A_t u = v$ and $u = A_t^{-1} v \in \SW(\R^d)$, (since $A_t$ is a sectorial operator and $0 \in \sigma(A_t)$ for every $t \in [0,T]$), we observe that $ \IFT [\Theta_t^{-1} \FT v ] = A_t^{-1} v$. 
Thus (in general, by replacing $v$ by $u$) we have that $\FT[A_t^{-1} - A_s^{-1} u] = [\Theta_t^{-1} - \Theta_s^{-1}]\FT u$ for all $u \in \SW(\R^d)$. Now for $t_1, t_2 \in [0,T]$, $\beta'_t \in [\beta_3\, \beta_4]$, $\beta_3 \le \beta_4$, we can estimate $\Theta_t^{-1} - \Theta_s^{-1}$ in the following way:
\begin{align}
\label{eq:conLipCndf}
|\Theta_{t_1}^{-1} - \Theta_{t_2}^{-1}| &= |(1 + \psi_{t_1} )^{-s \over 2} - (1 + \psi_{t_2} )^{-s \over 2}| ={ |(1 + \psi_{t_2} )^{s \over 2} - (1 + \psi_{t_1} )^{s \over 2}| \over (1 + \psi_{t_2} )^{s \over 2}(1 + \psi_{t_1} )^{s \over 2}} \notag\\
&= {|(1 + \beta_{t_1} \psi_Z)^{s\over 2} - (1 + \beta_{t_2} \psi_Z)^{s \over 2} |\over  (1  + \beta_{t_1}\psi_Z)^{s\over 2}(1 +  \beta_{t_2}\psi_Z)^{s\over 2}} \notag\\
&\le {|(1 + \beta_{t_1} \psi_Z)^{s\over 2} - (1 + \beta_{t_2} \psi_Z)^{s \over 2} |\over (1  + \beta_1 \psi_Z)^{s\over 2}(1 + \beta_1 \psi_Z)^{s\over 2}} \notag\\
&\le {1 \over (1 + \beta_1 \psi_Z)^s} (1 + \beta_t \psi_Z)^{s - 2 \over 2} |\beta'_t| \: | |(\beta_{t_1} - \beta_{t_2})| \psi_Z , \quad \text{for $t \in (t_1, t_2)$} \notag\\
%
&\le {C |\beta_4| \psi_Z \over (1 + \beta_1 \psi_Z)^{s + 2 \over 2}} |\beta_{t_1} - \beta_{t_2}| , \quad \text{ for $\beta'_t \in [\beta_3, \beta_4]$} \notag\\
&\le {C |\beta_4|({1 \over \beta_1} + \psi_Z) \over (1 + \beta_1 \psi_Z)^{s + 2 \over 2}} |\beta_{t_1} - \beta_{t_2}|, \quad \text{since $\beta_1 > 0$} \notag\\
&\le {C {|\beta_4| \over \beta_1}(1 + \beta_1 \psi_Z) \over (1 + \beta_1 \psi_Z)^{s + 2 \over 2}} |\beta_{t_1} - \beta_{t_2}| = {{|\beta_4| \over \beta_1}\over (1 + \beta_1 \psi_Z)^{s\over 2}} |\beta_{t_1} - \beta_{t_2}| \notag\\
%
&= {C |\beta_4| \over \beta_1} (1 + \beta_1 \psi_Z)^{-s \over 2} | \beta_{t_1} - \beta_{t_2}| \\ 
&\le {C |\beta_4| \over \beta_1} (1 + \beta_1 \psi_Z)^{-r \over 2} | \beta_{tvvvv_1} - \beta_{t_2}|, \quad \text{for $1 < r \le s $}
\end{align}
Recalling that $(A_t)^{\nu} = (1 - \psi_t(D))^{r \over 2}$, $\FT[(A_t)^{\nu} u] = (1 + \psi_t)^{r \over 2} \FT u$ and applying Lemma \ref{lem:psdOpComp} we see that 
\begin{align}
\label{eq:conLipCndfOp}
\|A_{t_1}^{\nu}(A_{t_1}^{-1} - A_{t_2}^{-1}) u \|_2 &= \| \FT[A_{t_1}^{\nu}(A_{t_1}^{-1} - A_{t_2}^{-1}) u ] \|_2 = \| (1 + \psi_t)^{r \over 2} \FT[(A_{t_1}^{-1} - A_{t_2}^{-1}) u ] \|_2 \notag\\
&\le \| (1 + \psi_t)^{r \over 2} {C |\beta_4| \over \beta_1} (1 + \beta_1 \psi_Z)^{-r \over 2} | \beta_{t_1} - \beta_{t_2}| \FT u ] \|_2 \notag\\
&= {C |\beta_4| \over \beta_1}  | \beta_{t_1} - \beta_{t_2}| \| (1 + \psi_t)^{r \over 2} (1 + \beta_1 \psi_Z)^{-r \over 2} \FT u\|_2 \notag
\end{align}
\begin{align}
&\le {C |\beta_4| \over \beta_1}  | \beta_{t_1} - \beta_{t_2}| \| (1 + \beta_2
 \psi_Z)^{r \over 2} (1 + \beta_1 \psi_Z)^{-r \over 2} \FT u\|_2 \notag\\
 &\le {C |\beta_4| \over \beta_1}  | \beta_{t_1} - \beta_{t_2}| \| (1 + \beta_2 \psi_Z)^{r \over 2} (1 + \beta_1 \psi_Z)^{-r \over 2} \FT u\|_2 \notag\\
  &\le {C |\beta_4| \over \beta_1}  | \beta_{t_1} - \beta_{t_2}| \| (\beta_2)^{r \over 2}({1 \over \beta_2}+ \psi_Z)^{r \over 2} (1 + \beta_1 \psi_Z)^{-r \over 2} \FT u\|_2 \notag
\end{align}
\begin{align}
  &\le {C |\beta_4| \over \beta_1}  | \beta_{t_1} - \beta_{t_2}| \| (\beta_2)^{r \over 2}({1 \over \beta_1}+ \psi_Z)^{r \over 2} (1 + \beta_1 \psi_Z)^{-r \over 2} \FT u\|_2 \notag\\
 &\le {C |\beta_4| \over \beta_1}  | \beta_{t_1} - \beta_{t_2}| \| (\beta_2/\beta_1)^{r \over 2}( 1 + \beta_1 \psi_Z)^{r \over 2} (1 + \beta_1 \psi_Z)^{-r \over 2} \FT u\|_2 \notag\\
 &\le {C |\beta_4| \over \beta_1}  (\beta_2/\beta_1)^{r \over 2} \| \FT u \|_2 \: | \beta_{t_1} - \beta_{t_2}| \notag\\
 &\le {C |\beta_4| \over \beta_1}  (\beta_2/\beta_1)^{r \over 2} \| \FT u \|_2 \: |\beta'| \: | t_1 - t_2| \notag\\
 &\le {C |\beta_4|^2 \over \beta_1}  (\beta_2/\beta_1)^{r \over 2} \| u \|_2 \:  | t_1 - t_2| 
\end{align}
\noindent
Based on \eqref{eq:conLipCndfOp} and \eqref{eq:conLipCndf} we get the following result:
\begin{lem}\label{lem:randOpLem2} Let $(\beta_t)_{t \in [0,T]}$ be a bounded stochastic process that is continuous and has bounded derivatives, i.e. $\beta_t \in [\beta_1, \beta_2]$ and $\beta'_t \in [\beta_3, \beta_4]$ for all $t \in [0,T]$, with $\beta_1 > 0$. Moreover, let $s \in \R+$, $0 < r < s$. Then for $\nu = {r \over s}$ we get that 
\begin{align}
\|A_{t_1}^{\nu}(A_{t_1}^{-1} - A_{t_2}^{-1}) u \|_2 \le (\beta_2/\beta_1)^{r \over 2}{|\beta_4|^2 \over \beta_1}  | t_1 - t_2 | 
\end{align}
\end{lem}

\begin{figure}[h!]
\centering
	\includegraphics[scale=.5]{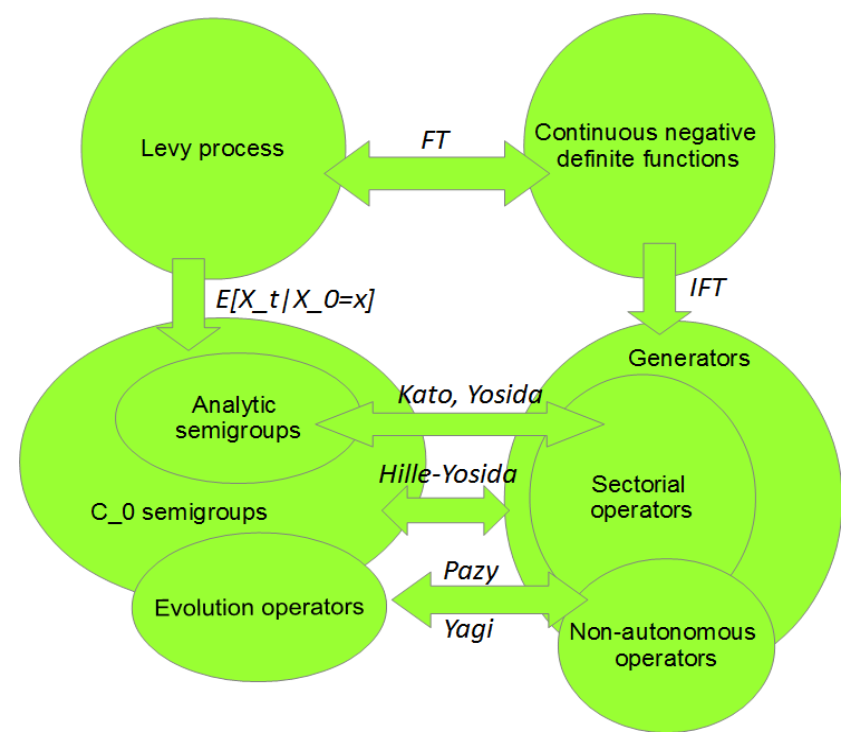} 
\caption{The framework for obtaining random abstract Cauchy problems via random operators. \label{fig:randOpFrm}}
\end{figure}

\begin{thm}\label{thm:randACP} Let $(\Theta_t)_{t \in [0,T]}$, $\Ag := (A_t)_{t \in [0,T]}$ ,$(\beta_t)_{t\in[0,T]}$ be as in Lemma \ref{lem:randOpLem} and Lemma \ref{lem:randOpLem2}. Then, for every $\om \in \Omega$ fixed, there exists a unique  two parameter evolution operator $ \Ug := \evolOp$ on $L^2(\R^d)$ associated to $\Ag$,  such that for every $\om \in \Om$ fixed, and $f \in D((A_t)^{\nu}) \subset L^2(\R^d)$ we get the following random abstract Cauchy problem
\begin{align}
\label{eq:randACP}
\begin{split}
{d \over dt}U_{t,0}(\om) f &= -A_t(\om) U_{t,s}(\om) f , \quad \text{ in $L^2(\R^d)$}\\
U_{0,0}f &= f.
\end{split}
\end{align}
\end{thm}
\begin{proof} Follows from Lemma \ref{lem:ndfClosed} Lemma \ref{lem:randOpLem}, Lemma \ref{lem:randOpLem2}, and Theorem \ref{thm:evolOpSecOp}.
\end{proof}
\noindent
The preceding discussion gives an abstract method for constructing random operators that satisfy the Ansatz \ref{asm:nonAtndfOp} and thus generate a two parameter evolution semigroup which establishes the abstract Cauchy problem which is stochastic by construction. The source of randomness stems from the switching noises operating at the microscopic level. Figure \ref{fig:randOpFrm} illustrates the mathematical framework used for this construction. 

\noindent
Now we provide some concrete examples.

\begin{eg}\label{eg:randOp1}Let $(\beta_t)_{t\in[0,T]}$ be a bounded continuous process defined in the following way:
\begin{align*}
\beta_t := \beta_1 + {\beta_2 I_t \over 1 + I_t}, \text{ with }  I_t := \int_0^t \sin(B_s^{\beta_3,\beta_4})^2 ds, 
\end{align*} 
where
\[ B_t^{\beta_3, \beta_4} = {1 \over T} \Big( (T-t) \beta_3 + T W_t + t (\beta_4  - W_T) \Big), \text{ for $t\in[0,T]$} \]
is the Brownian bridge process starting at $\beta_3$ and ending at $\beta_4$ and $(W_t)_{t \in [0,T]}$ is the standard Wiener process. Let $\psi_A(\xi) = |\xi|^2$, $\psi_S(\lam) = \lam^{\alp}$ for $\alp \in (0,1)$. Then $\psi_t = \beta_t (\psi_S \circ \psi_A) = \beta_t |\xi|^{2\alp}$. Then $\Theta_t$ is given as 
\[ \Theta_t(\xi) = (1 + \beta_t |\xi|^{2\alp})^{s \over 2} \]
where $s \in (0,2)$, $\beta_1 > 1$, $\beta_2 > \beta_1$ and $0 < \beta_3 < \beta_4$.
\end{eg}

\begin{eg}\label{eg:randOp2} Let $\beta_t$ and $\psi_A$ be defined as in Example \ref{eg:randOp1}, letting $\psi_S =1 $ we get 
\[ \Theta_t(\xi) = (1 + \beta_t |\xi|^{2})^{s \over 2}. \] This operator is nothing but the fractional Laplacian.
\end{eg}

\begin{eg}\label{eg:randOp3} In this example we use a more direct approach and define $\Theta_t$ as
\[ \Theta_t(\xi) = |\xi|^{2\alp_t} \]
where $\alp_t = {\alp_1 + \alp_2 I^2_t \over 1 + I_t^2}$, such that $\alp_t \in [\alp_1, \alp_2] \subset (0,1)$. The fact that the operator associated to $\psi_t = \alp_t$ satisfies the Ansatz \ref{asm:nonAtndfOp} will be established by Lemma \ref{eq:fracOpLipCon} in Section \ref{sec:fracSPDE}.
\end{eg}
\noindent
Next we provide a generalization of the above procedure and consider a general Banach space valued random equation and deduce an analogous functional equation. \\
\noindent
Let $\evolOp$ be the two parameter semigroup on $L^2(\R^d)$ constructed in Theorem \ref{thm:randACP} and $\Ag=(A_t)_{t>0}$ the corresponding family of sectorial operators. Then we have that, $U_{t,s}$ fulfills the following (operator) equations

\begin{equation*}
\begin{split}
{\partial \over \partial t} U_{t,s} &= -A_t U_{t,s} \\
U_{s,s} &= I 
\end{split}
\quad \quad \text{ and } \quad \quad
\begin{split}
{d \over ds}U_{t,s} &= U_{t,s} A_s\\
U_{t,t} &= I
\end{split}
\end{equation*}

In particular, for any $f \in D((A_t)^{\nu}) \subset L^2(\R^d)$ we have that 
\begin{align}\label{eq:nonAutoACP}
{d \over dt} U_{t,0}f &= -A_t U_{t,0} f \quad \text{ in $L^2(\R^d)$}\\
U_{0,0}f &= f \notag
\end{align}

\begin{prop}
\label{prop:uMP}
Let us assume that there exists a unique solution $u_t = U_{t,0} f$ to the non-autonomous abstract Cauchy problem \eqref{eq:nonAutoACP}. Then we get that $u = (u_t)_{t \ge 0}$ is a $L^2(\R^d)$ valued Markov process. 
\end{prop}
\begin{proof} Without loss of generality, let $t > s > r$. Since $u_t = U_{t,0} f = U_{t,s}U_{s,r}U_{r,0}f$, for some given value of $u_s$ say  $g$ (i.e. $u_s = U_{s,0} f = g$) we have that $u_t = U_{t,s}g$ which implies that 
\[(u_t)_{t >s} \underset{u_s = g}{\coprod} (u_r)_{r< s},\] i.e.  given $g$, $u_t$ for $t > s$ is independent of the values of $u_r$ for all $r < s$. Since $u_t \in L^2(\R^d)$, altogether we get that $u_t$ is an $L^2(\R^d)$ valued Markov process.
\end{proof}
\noindent
With this we want to define a semigroup on $C_b(L^2(\R^d))$ via the prescription $T_{s,t} f = \E(f(u_t) | u_s)$ for all $f \in C_b(L^2(\R^d))$ and $t \ge s \ge 0$. 
\begin{prop} Let $u_t = U_{t,0} f$ be a unique solution to the non-autonomous abstract Cauchy problem \eqref{eq:nonAutoACP}, then $T_{s,t} f (u_s) = \E(f(u_t) |u_s)$, for $0 \le s \le t < \infty$ is a bounded linear operator on $C_b(L^2(\R^d))$ and is a $\CSG$.
\end{prop}
\begin{proof} Let $f \in C_b(L^2(\R^d))$ and let $\| \cdot \|_0$ denote the norm on $C_b(L^2(\R^d))$ i.e. \[  \| f \|_0 = \sup_{x \in L^2(\R^d)} |f(x)|, \quad \: \forall f \in C_b(L^2(\R^d)),\] then we have the following assertions:
\begin{enumerate}
\item $T_{s,t} \in \Lu(C_b(L^2(\R^d)))$: Let $t \ge s \ge 0$, then 
\begin{align*}
\|T_{s,t} f\|_0 = \|\E(f(u_t) | \cdot)\|_0 \overset{Jensen's}{\underset{inequality}{\le}} \E(\|f(u_{\cdot})\|_0 | u_s) \le \|f\|_0.
\end{align*}
Thus $T_{s,t}$ is a contraction operator.
\item $T_{r,t} = T_{r,s}T_{s,t}$: Let $0 \le r \le s \le t < \infty$, then 
\begin{align*}
T_{r,t}f(u_r) = \E(f(u_t)|u_r) &\overset{Tower}{\underset{property}{=}} \E( \E(f(u_t) | \F_s) | u_r) \overset{Markov}{\underset{property}{=}}\E( \E(f(u_t) |u_s) | u_r) \\
&= \E((T_{s,t} f)(u_s) | u_r) = T_{r,s}T_{s,t} f(u_r)
\end{align*}
\item $\underset{ t \to s}{\lim} \|T_{s,t}f - f\|_0 = f$: For this we make use of the fact that for a.e. $\om \in \Om$ the evolution operator $U_{t,s}$ is strongly continuous w.r.t $(t,s)$. Thus $u_t = U_{t,0} f$ is continuous w.r.t $t$ in $L^2(\R^d)$, i.e. $\underset{t \to s}{\lim} \|u_t - u_s\|_{L^2(\R^d)} = 0$ for $0 \le s \le t$.
\begin{align*}
\|T_{s,t}f - f\|_0 &= \|\E(f(u_t) | \cdot ) - f(u_s)\|_0 = \| \E(f(u_t) - f(u_s)| \cdot) \|_0 \\
&\le \E(\|f(u_{\cdot}) - f(u_{\cdot})\|_0 | u_s) \le \|f(u_{\cdot}) - f(u_{\cdot})\|_0
\end{align*}
Since $ u_t \to u_s$ as $t \to s$ for a.a. $\om \in \Om$ and $f \in C_b(L^2(\R^d))$ we have that $f(u_t) \to f(u_s)$ as $t \to s$. As a result we get that $T_{s,t} \to I$ strongly as $t \to s$. In fact, since $U_{t,s}$ is strongly continuous in $(s,t) \in \R^+\times \R^+$, i.e. $U_{t,s} \to I$ as $(t,s) \to (r,r)$, it holds that $(s,t) \mapsto T_{s,t}$ is strongly continuous as well.
\end{enumerate}
\end{proof}
\noindent
In order to establish a Cauchy problem associated to the (time inhomogeneous) semigroup $(T_{s,t})_{0\le s<t}$, we find it convenient to make $T$ time homogeneous. This is possible by augmenting the state space of the Markov process $(u_t)_{t \ge 0}$ (defined in Proposition \ref{prop:uMP}) with the time state space, i.e. by converting the nonautonomous ODE to an autonomous ODE by taking the deterministic time variable to be a part of the state space.
\begin{align}
\label{eq:augStateSp}
{d \over d \tau} u_{\tau} &= A_{\tau} u_{\tau} & {d \over dt} u_{_{\tau_t}} &= A_{_{\tau_t}} u_{_{\tau_t}} \notag\\
u_0 &= f  \hspace*{3.5cm} {\Longrightarrow} &{d \over dt} \tau_t &= 1 \\
& & u_{_{\tau_0}} &= f\, \tau_0 = 0. \notag
\end{align}
\noindent
Following this transformation we should also change the probability transition function to incorporate the augmented deterministic time variable. It turns out that the new probability transition function is time homogeneous in some sense, as shown in the proposition below.

\begin{prop} 
\label{prop:HomMarkovProc}
Let $(E, \|\cdot\|_E)$ be a Banach space and $\B$ the Borel sigma algebra on $E$ induced by its norm $\|\cdot\|_E$. Let $(u_{\tau})_{\tau \ge 0}$ be an $(E, \B)$ valued Markov process (given by the LHS of \eqref{eq:augStateSp}) with the probability transition function $p_{\tau_1,\tau_2}(x, B )$ for $0 \le \tau_1 \le \tau_2$ and $x \in E$ and $B \in \B$.  Then the process $(\hat u_t)_{t \ge 0} = (\tau_t, u_{t})$ (given by the RHS of \eqref{eq:augStateSp}) with the probability transition function $\hat p_{s,t}$ defined as 
\[ \hat p_{s,t}(\hat x, \hat B) = \P(u_{t} \in B | u_{s} = x) \otimes \delta_{t}( I | \tau_s = s) \] with $\hat x = (s, x)$, $\hat B = B \times I$, and $I$ a Borel subset of $\R^+$ is a time homogeneous Markov process.
\end{prop}
\begin{proof}

For $(\hat u_t)_{t \ge 0} = (\tau_t, u_{_{\tau_t}})$, and $\hat x = (s, x)$ we have that 
\begin{align*}
\hat{p}_{s,t}(\hat x, \hat B) = \P(\hat u_{t} \in \hat B | \hat u_{s} = \hat x) &= \P\Big((\tau_t,u_{_{\tau_t}}) \in I \times B \Big| (\tau_0,u_{_{\tau_0}}) = (s,x)\Big)
\end{align*}
and
\begin{align*}
\hat p_{0, t-s}(\hat x , \hat B) &=\P\Big((\tau_{t-s},u_{_{\tau_{t-s}}}) \in I \times B \Big| (\tau_0,u_{_{\tau_0}}) = (0,x)\Big)
\end{align*}
Since the initial condition $\hat x = (s, x)$ is fixed and $\tau_t = t - \tau_0$, 
we see that  $\tau_t = t-s$ if $\tau_0 = s$ and $\tau_{t-s} = t-s$ if $\tau_{0}= 0$. Thus the evolution of $A_{\tau_t}$ and hence of $u_{\tau_t}$ in the time interval $(s,t]$ and $(0,t-s]$ remains the same. Consequently,  we get that $\hat p_{s,t}(\hat x, \hat B) = \hat p_{0,t-s}(\hat x, \hat B)$.
\end{proof}
\noindent Now if $(u_t)_{t \ge 0}$ satisfies the Cauchy problem \ref{eq:nonAutoACP}, we get that $(\hat u_t)_{t \ge 0}$ is a $(\R^+ \times L^2(\R^d))$ valued homogeneous Markov process. From now on we drop the hat and refer $(u_t)_{t \ge 0}$ to be our newly constructed homogeneous Markov process $(\hat u_t)_{t \ge 0}$.\\
\noindent
Now letting $D_u f$ denote the Fr\'{e}chet derivative of $f$ at point $u$ , applying Taylor's Theorem (which yields us the remainder term $R_{u,e}$),we can deduce the result. \\
\noindent
\begin{prop}
Let $(T_t f)(u0) := \E(f(u_t) | u_0 =  u0)$ be a $\CSG$ on $C_b(\hLt)$, where $\hLt := \R^+ \times L^2(\R^d)$.  Let $e = u_{t+h} - u_t$, then for sufficiently smooth $f$, we have that 
\[ {d \over dt} T_t f = T_t (G f), \text{ with } G f = (D_{u_t} f)(A_t u_t).\]
\end{prop}
\begin{proof}
\begin{align*}
{T_{t+h} f - T_t f \over h} &= {1 \over h} \E\Big[ f(u_{t+h}) - f(u_t) \Big| u_0 \Big]\\
&= {1 \over h} \E\Big[f(u_t + e) - f(u_t) \Big| u_0\Big] \\
&= {1 \over h}\E\Big[ (D_{u_t} f)(e) + \big(R_{u_t, e}\big)(e) \Big| u_0 \Big] \\
&=  \E \Big[ (D_{u_t} f)({e \over h}) + \big(R(u_t,e)\big)({e \over h}) \Big| u_0 \Big] \\
&=  \E \Big[ (D_{u_t} f)({u_{t+h} - u_t \over h}) + \big(R_{u_t,e}\big)({u_{t+h} - u_t \over h}) \Big| u_0 \Big]
\end{align*}
Now taking the $h \dto 0$ on both sides, we get 
\begin{align*}
{d \over dt} T_t f &= \lim_{h \dto 0} \E \Big[ \big( f(u_t +e) - f(u_t) \big) \Big| u_0 \Big] \\
&= \E \Big[ \lim_{h \dto 0} \big( f(u_t +e) - f(u_t) \big) \Big| u_0 \Big] \\
&= \E \Big[ \lim_{h \dto 0} \big( (D_{u_t} f)({u_{t+h} - u_t \over h}) + \big(R_{u_t,e}\big)({u_{t+h} - u_t \over h}) \big)\Big| u_0 \Big] \\
&= \E \Big[ (D_{u_t} f)(A_t u_t)+ \lim_{h \dto 0} \big(R_{u_t,e}\big)({u_{t+h} - u_t \over h}) \Big| u_0 \Big] \\
&= \E \Big[ (D_{u_t} f)(A_t u_t) \Big| u_0 \Big] \\
&= T_t (G f)
\end{align*}
For $f \in C^1(\hLt) \cap C_b(\hLt)$ we have that the mapping $ \hLt \ni u \mapsto D_u f \in \Lu(\hLt; \R)$ is continuous in $u$. Moreover, the remainder term $R_{u_t, e}(e^2) \in o(e^2)$ and $R_{u_t, e} \to 0$ in $\Lu(\hLt; \R)$ as $e \to 0$. Since, $u \in C([0,T]; L^2(\R^d))$ we have that $e = u_{t+h} - u_t$ converges to $0$ in $\hLt$ as $h \to 0$.
\end{proof}
\noindent
This completes the general multiscale framework for deducing stochastic equations at the macroscopic level.  Next we formally apply this method and propose a macroscopic stochastic equation for the dynamics of cancer invasion under the influence of tissue acidity.


\newcommand{\fracL}{\text{fractional Laplacian }}
\newcommand{\fracLap}{(-\lap)^{\alp}}
\newcommand{\fLap}{(-\lap)^{{\alpha \over 2}}}
\newcommand{\sign}{\text{sign}}

\section{Numerical examples \label{sec:fracSPDE}}
In this section we shall provide some numerical examples that illustrate the methodology for modeling diversity of a complex biological systems. Here we choose cancer as a prototypical example of such a system, and we specifically focus on modeling the diverse invasion capability of cancer cells under the influence of tissue and cellular acidity. To this end we provide two approaches: namely microscopic and macroscopic approaches. In the first approach we consider the dynamics at the cell-level interactions while in the second approach we consider the dynamics at the tissue-level interactions. 

\subsection{Microscopic modeling:}
\begin{figure}
	\centering
    \subfloat[$X_0$ and $\nabla N_0$]{\includegraphics[trim = 0cm 0cm 0cm 0cm, clip,scale=.225]{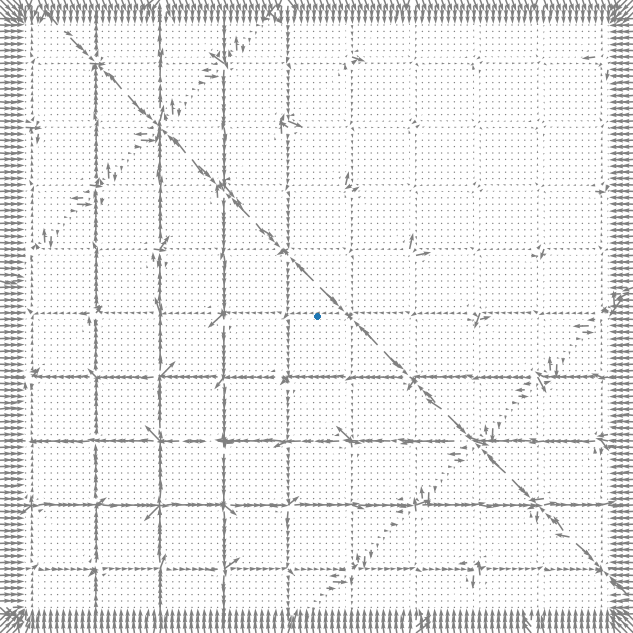}}\qquad \qquad
	\subfloat[$H_i(0,X_0)$]{\includegraphics[trim = 0cm 0cm 0cm 0cm, clip,scale=.275]{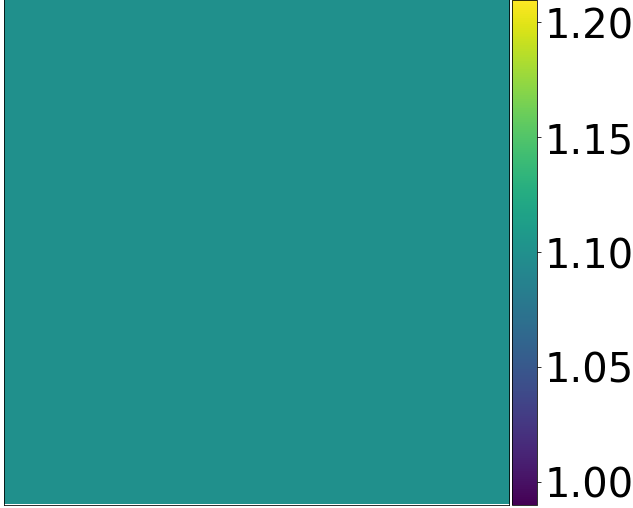}}\\
	\subfloat[$H_e(0,x)$]{\includegraphics[trim = 0cm 0cm 0cm 0cm, clip,scale=.275]{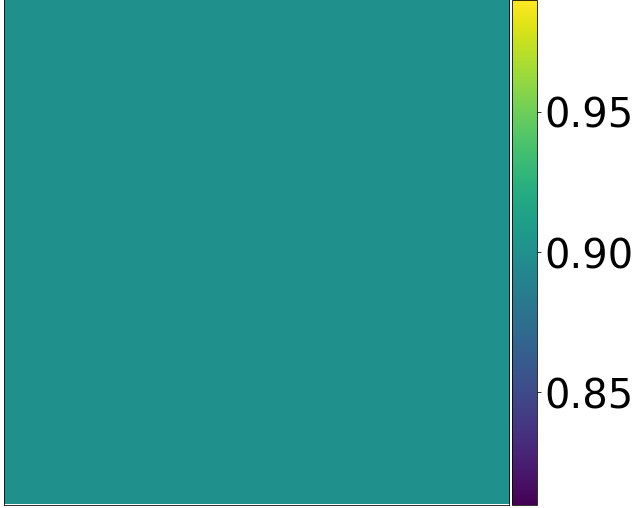}}\qquad
	\subfloat[$N_i(0,x)$]{\includegraphics[trim = 0cm 0cm 0cm 0cm, clip,scale=.275]{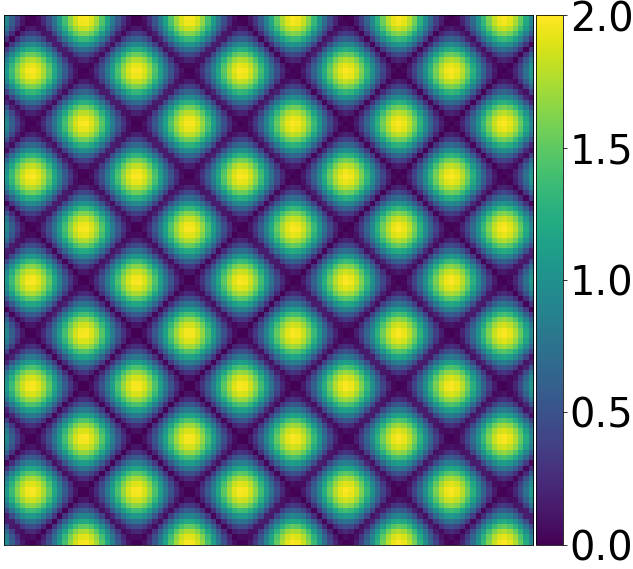}}
\caption{Initial condition for the microscopic model \eqref{eq:micAcidCan}. Starting from the left, the first plot shows the directional vector field generated by underlying tissue structure in gray background. The blue dots show different cell particles and their velocity vectors are indicated by red arrow. Second plot from the left shows the initial concentration of intra-cellular protons $H_i$ of all $M=2500$ particles arranged in as 2D grid. Third plot from the left shows the concentration of extra-cellular protons $H_e$. The fourth plot from the right shows the density of normal cells i.e the structure of the underlying tissue.  \label{fig:micAicdCanInit}}
\end{figure}

In this section we consider the modeling of invasive dynamics of cancer cells under the influence of acidity while keeping in mind to take into account diversity aspect. Based on this the main quantities under considerations are $H_i$, $H_e$ which capture the intra- and extra-cellular concentration of protons. The motion of cancer cells itself is modeled via first order Newtonian dynamics, namely in terms of velocity variable $V$ and position variable $X$. Finally, the structure or the density of the underlying tissue structure is modeled by $N$. Based on this, the microscopic dynamics is given as: 
\begin{subequations}
\label{eq:micAcidCan}
\begin{align}
    {d V_t} &=  \grad N_i dt + dL_t \label{eq:micAcidV}\\
    {d X_t} &=  V_t dt \label{eq:micAcidX}\\
    {d H_i\over dt}(t,X_t)  &= -T(H_i,H_e) - S_1(H_i) + Q(H_i) \label{eq:micAcidHi}\\
    {d H_e\over dt}(t,X_t)  &= T(H_i,H_e) - S_2 \label{eq:micAcidHe}\\
    {d N_i \over dt}(t,X_t) &= -\gamma N_i H_e \label{eq:micAcidN}
\end{align}
\end{subequations}
The intra-cellular proton dynamics is affected by: (i) production of protons due to glycolysis $Q(H_i)$, (ii) buffering of protons due to intracelluar and organellar buffer $S_1(H_i)$ and (iii) efflux of intracellular protons due to membrane proton transporters $T_1, T_2$ and $T_3$. These terms are summed up and represented by a single term $T$. Similarly, the extra-cellular proton dynamics is affected by: (i) influx of protons by membrane transporters $T$ and (ii) sequestering of excess protons by the tissue vasculature $S_2$. The effect of excess protons in the extra-cellular region accelerate the activity protolytic activity of matrix-degrading enzymes because of which there is a loss of extra-cellular tissue. Thus $N$ is just modeled by a decay term. Lastly, the movement of cancer cells is specified (in general) by a jump velocity model where in the velocity process $V$ is given by the SDE \eqref{eq:micAcidX} and the position process $X$ is simply the integral of the velocity term \eqref{eq:micAcidV}. Additionally, the mean velocity of the cancer cell is assumed to be specified by the fiber orientation $\nabla N$  of the underlying tissue structure. Altogether the model is as given in \eqref{eq:micAcidCan}. Since, $N$ and $H_e$ are typically macroscopic variables, their dynamics is affected only at the position of the cancer cells, thus in the equation we see that $N$ and $H_e$ are written as a function $X_t$. Their macroscopic version $N(t,x)$ and $H_e(t,x)$, respectively, are obtained simply by convolving $N(t,X_t)$ and $H_e(t,X_t)$ by the density function of $X_t$. 
\begin{align*}
    H_e(t,x) &=  \int_{\dom} f_X(x-y)H_e(t,x) dy,  \quad    N(t,x) =  \int_{\dom} f_X(x-y)N(t,x) dy
\end{align*}
We now perform some numerical experiments for the above model. We consider a 2D spatial domain $\dom = [0,1]\times[0,1]$ and discretize it uniformly. Similarly, the temporal component is discretized into $N$ steps with time step $\tau = .05$. Since we are interested in microscopic simulation, we consider $M$ different cancer cells which are idealized as particles. Each particle undergo motion as per \eqref{eq:micAcidV} and \eqref{eq:micAcidX} while exhibiting intra-cellular dynamics (for  $H_i$) as per \eqref{eq:micAcidHi}. Due to the coupling between $H_i(\cdot, X_t)$ and $H_e(\cdot,X_t)$, the macroscopic extracellular proton concentration $H(t,x)$ is affected along the path of the particles. Analogously, the macroscopic normal tissue density is $N(t,x)$ is also affected. Based on this, the initial condition for the above SDE as shown in Figure \ref{fig:micAicdCanInit}. We now run three different simulations with $M=2500$ and $N=25$, each with a different type of noise processes $dL_t$ for the velocity equation \eqref{eq:micAcidV}. The different types of processes used for the simulation are as defined in \eqref{eq:micModCanLP}. In order to understand the 
effect of these different types of noise distribution, we compute a macroscopic quantity at the end of the simulation and compare it for all different cases. The macroscopic variable here is the mean survival percentage of cancer cells. To enable this we incorporate a condition based on which we kill the particle process $X_t$ . To be more precise, all the particles $X^j_t$, with either $H_i(t,X^j_t) < h_1$ or $H_i(t,X^j_t) > h_2$ or $H_e(t, X^j_t) > h_3$ are killed, i.e. are removed from the system. The remaining percentage of particles (cells) determine the survival percentage of the population. 
\begin{align*} 
\begin{array}{l l }
S_{H_1} := \Big \{j \in \N, j \le M: H_i(t,X^j_t) < h_1 \Big \},& S_{H_2} := \Big \{j \in \N, j \le M: H_i(t,X^j_t) > h_2 \Big \}, \\
S_{He} :=  \Big \{j \in \N, j \le M: H_e(t, X^j_t) > h_3 \Big \},& S_0 := \Big\{ j \in \N: j \le M \Big\},  S =  { \# \tilde S \over M}, \\
S_1 := S_{H_1} \bigcup S_{H_2} \bigcup S_{He},& \tilde S := S_0 - S_1
\end{array}
\end{align*}
\vspace*{-.5cm}
\begin{subequations}
\label{eq:micModCanLP}
\begin{align}
    dL_t &\sim \mathcal{N}(0,1) \label{eq:fxdNoise}\\
    dL_t &\sim \begin{cases}
                \mathcal{N}(0,1) & \text{ if } \qquad U \in [0, .3), U \in \mathcal{U}(0,1) \\
                \text{Laplace}(0,1) & \text{ else if} \:\: U \in [.3, .5), \\
                \text{Triangular}(-4,0,8) & \text{ else }
          \end{cases} \label{eq:swNoise}\\
    dL_t &\sim 10 \sin(\sigma_t) \mathcal{N}(0,1), \quad \sigma_t \sim \text{Cauchy}(0,1) \label{eq:swVarNoise}
\end{align}
\end{subequations}
\begin{figure}
    \centering
    \subfloat[$X_T$ and $\nabla N_T$]{\includegraphics[trim = 0cm 0cm 0cm 0cm, clip,scale=.225]{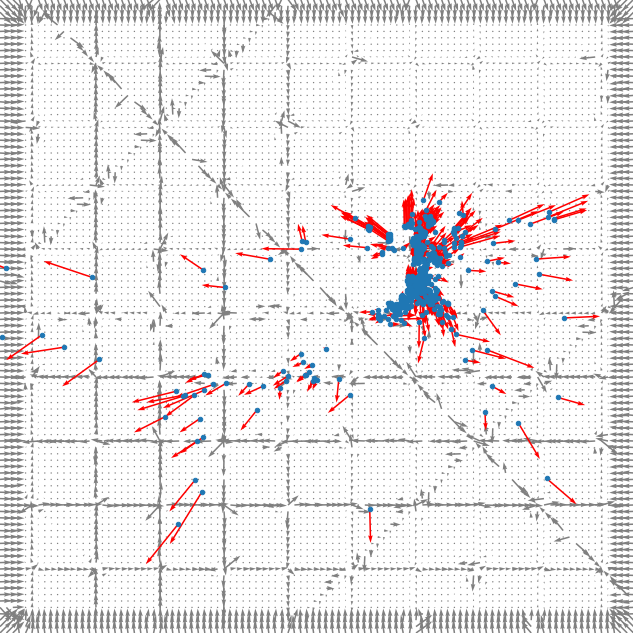}}\qquad  \qquad
	\subfloat[$H_i(T,X_T)$]{\includegraphics[trim = 0cm 0cm 0cm 0cm, clip,scale=.275]{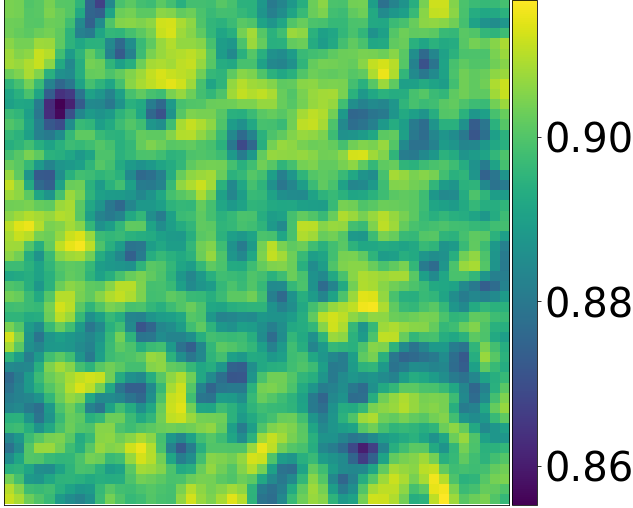}}\\
	\subfloat[$H_e(T,x)$]{\includegraphics[trim = 0cm 0cm 0cm 0cm, clip,scale=.275]{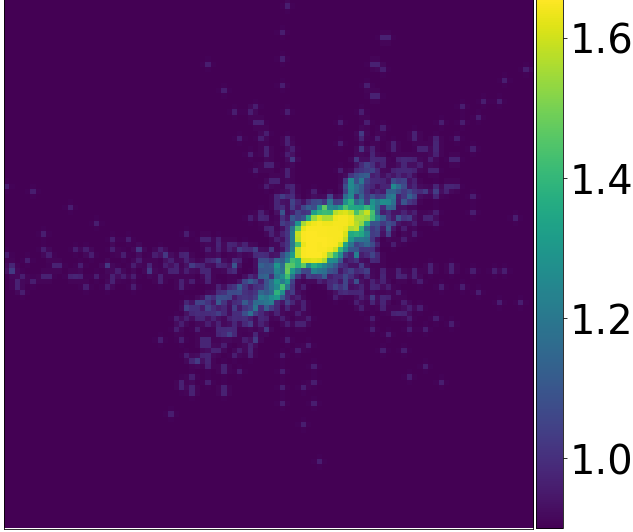}}\qquad
	\subfloat[$N_i(T,x)$]{\includegraphics[trim = 0cm 0cm 0cm 0cm, clip,scale=.275]{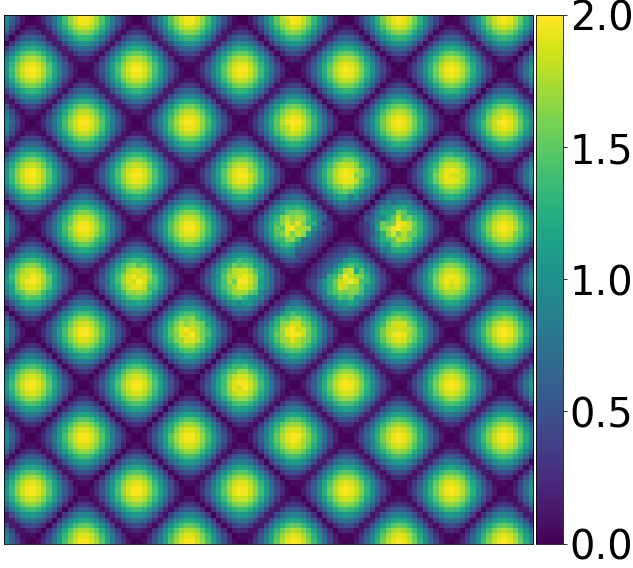}}
\caption{ Simulation results for the microscopic model \eqref{eq:micAcidCan} at time $T=2.5$s when taking $dL_t$ as per \eqref{eq:fxdNoise}. The resulting mean survival percentage $S \approx 35\%$ \label{fig:micAicdCanGauss}}
\end{figure}
\begin{figure}
    \centering
    \subfloat[$X_T$ and $\nabla N_T$]{\includegraphics[trim = 0cm 0cm 0cm 0cm, clip,scale=.225]{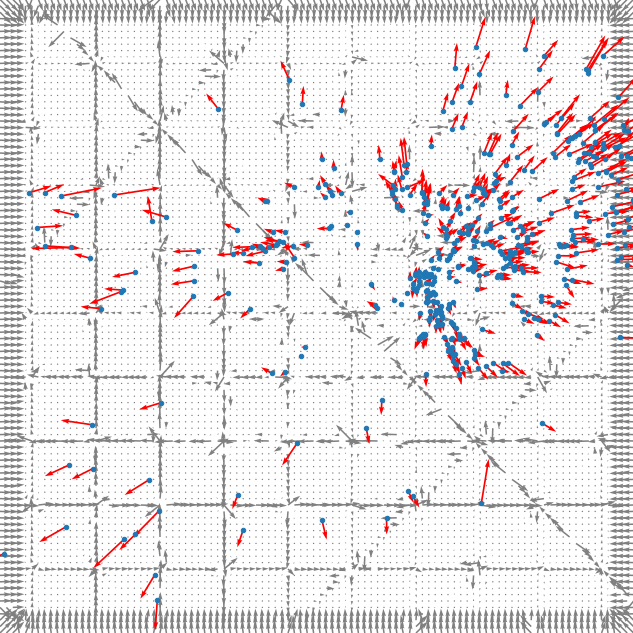}}\qquad  \qquad
	\subfloat[$H_i(T,X_T)$]{\includegraphics[trim = 0cm 0cm 0cm 0cm, clip,scale=.275]{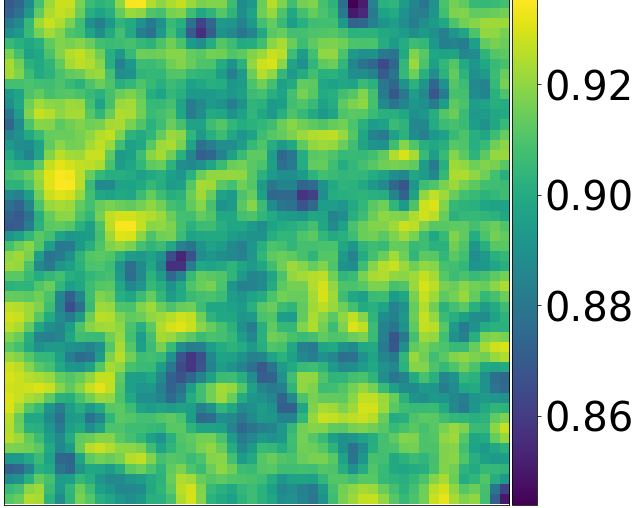}}\\
	\subfloat[$H_e(T,x)$]{\includegraphics[trim = 0cm 0cm 0cm 0cm, clip,scale=.275]{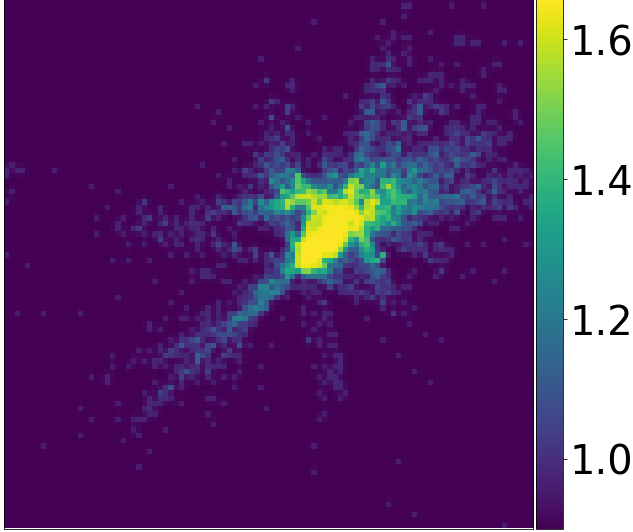}}\qquad
	\subfloat[$N_i(T,x)$]{\includegraphics[trim = 0cm 0cm 0cm 0cm, clip,scale=.275]{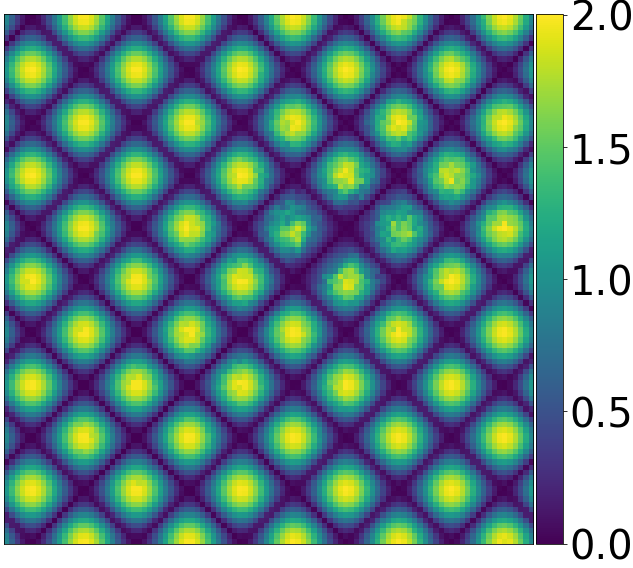}}
\caption{ Simulation results for the microscopic model \eqref{eq:micAcidCan} at time $t=2.5$s when taking $dL_t$ as per \eqref{eq:swNoise}. The resulting mean survival percentage $S \approx 40\%$ \label{fig:micAicdCanJumpLevy}}
\end{figure}
\begin{figure}
    \centering
    \subfloat[$X_T$ and $\nabla N_T$]{\includegraphics[trim = 0cm 0cm 0cm 0cm, clip,scale=.225]{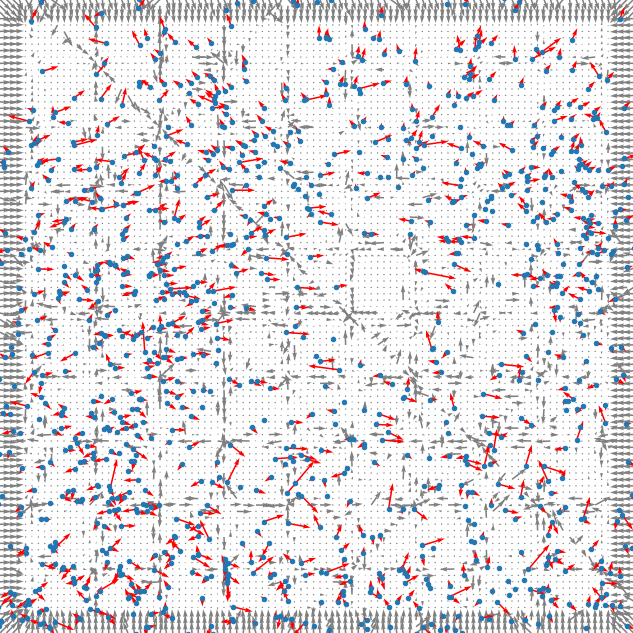}}\qquad  \qquad
	\subfloat[$H_i(T,X_T)$]{\includegraphics[trim = 0cm 0cm 0cm 0cm, clip,scale=.275]{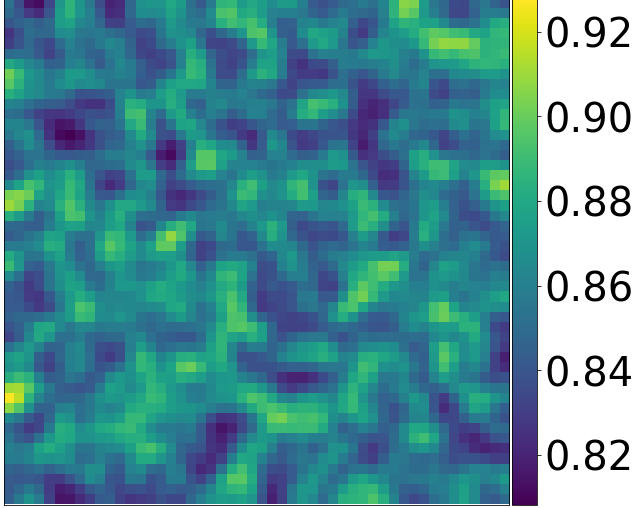}}\\
	\subfloat[$H_e(T,x)$]{\includegraphics[trim = 0cm 0cm 0cm 0cm, clip,scale=.275]{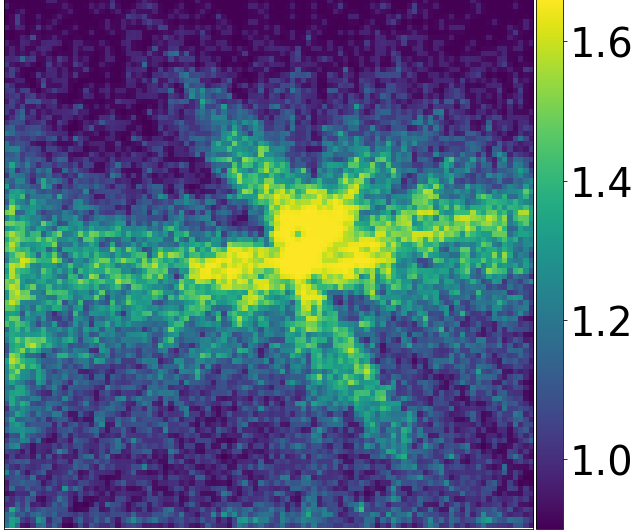}}\qquad
	\subfloat[$N_i(T,x)$]{\includegraphics[trim = 0cm 0cm 0cm 0cm, clip,scale=.275]{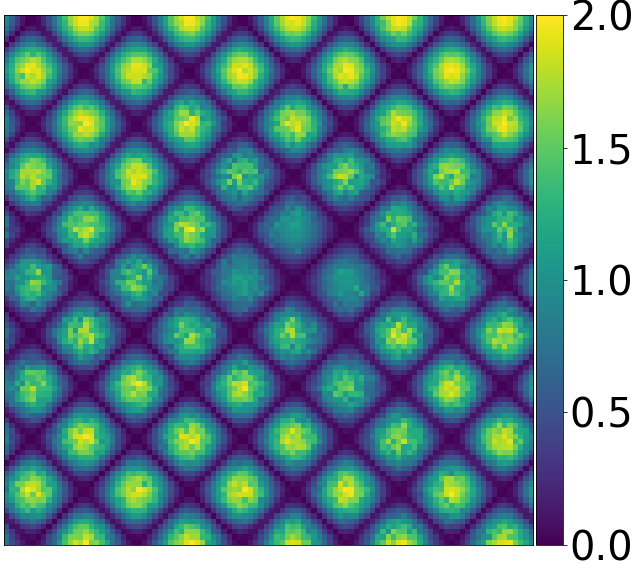}}
\caption{Simulation results for the microscopic model \eqref{eq:micAcidCan} at time $t=2.5$s when taking $dL_t$ as per \eqref{eq:swVarNoise}. The resulting mean survival percentage $S \approx 50\%$ \label{fig:micAicdCanLevy}}
\end{figure}
Based on the mean survival percentage $S$ and the final distribution of cells obtained for the three different cases, we can infer that the migration pattern generated by type-3 noise process \eqref{eq:swVarNoise} is the most invasive while the pattern generated by type-1 noise process \eqref{eq:fxdNoise} is least invasive. Even though one could also generate a highly invasive migration pattern by simply choosing a high variance parameter for the Type-1 noise, it fails to take into account the  intracellular and extracellular disturbances leading to a change in the distribution types such as from uni-modal to multi-modal, thin-tailed to heavy-tailed and so on. Such switching in the distribution types results in switch in the migration regime consequently exhibiting a diverse migration regime. Based on the discussion from Section \ref{sec:M2MACP} we can infer that Type-2 and Type-3 noise processes induce random operators at the macroscopic levels while Type-1 noise does not, thus make it unsuited for modeling diverse behaviors.
Motivated by this and due to the tight coupling, offered by the modeling framework (cf. Fig. \ref{fig:micMacScaleDiag}), between the micro and macro level,  we shall now consider a macroscoipic model for acid mediated cancer invasion involving random operators. 

\subsection{Macroscopic modeling:}
In this section we provide an example of continuous negative definite functions that satisfy Ansatz \ref{asm:nonAtndfOp}, and based on this example we propose a macroscopic model for the acid-mediated tumor invasion. \\
To be precise, in this section we explore the properties of the operator whose symbol is given by $\eta(\xi) = |\xi|^{2\alp}$, for $ \xi \in \R^d$ and $\alp \in (0,1)$. 
The function $|\xi|^{2 \alp}$ appears as a the symbol of the \LP $(Z_t)_{t \ge 0}$ obtained by the subordination of the Wiener process $W = (W_t)_{t\ge0}$ with respect to the $\alpha$-subordinator $S = (S_t)_{t\ge0}$, whose characteristics is given as $(0,\lam)$, with $\lam(dx) = {\alpha \over \Gamma(1-\alpha)}{ dx \over x^{1+\alpha}}$. Consequently, the subordinated process $Z = W_{S} = (W_{S_t})_{t \ge 0}$. Correspondingly, as a consequence of the Phillips theorem we find that the generator $A^Z$ of the subordinated process $Z$ is given by the fractional Laplacian $(-\lap)^{\alp}$ for $\alp \in (0,1)$ (as shown in example \ref{eg:alphaStable}) . In order to specify a time dependent random operator, we just let $\alp$ be a bounded stochastic process with values in $(0,1)$. First we shall recall the definition and properties of the fractional Laplacian $\fracLap$.


\begin{define}[Fractional Laplacian]\label{def:fracLap} Let $f \in L^p(\R^d)$, $p \in \{1, 2\}$. Then the operator $A$ defined on $L^p(\R^d)$ as 
\[ (\FT A f)(\xi) = -|\xi|^{2\alp} (\FT f)(\xi) \quad \alp \in (0,1) \]
is called the $\fracL$ on  $L^p(\R^d)$ with $A = -\fracLap$. 
\end{define}
\noindent
Though this is an easy and convenient way of defining the fractional Laplacian, one can find in the literature various different definitions, namely- Balakrishnan's definition, singular integral definition, Dirichlet form definition, Dynkin's definition and so on.  We now state a very nice theorem, due to [\citenum{Kwasnicki2017}], that establishes an equivalence relation between these various definitions:

\newcommand{\BSP}{\mathscr{B}}

\begin{thm} Let $\BSP$ be any of the spaces $L^p(\R^d)$, $C_0$ or $C_b$, $p \in [1,\infty)$. Let $f \in \BSP$ and $\alp \in (0,2)$. Then the following definitions of the $\fracL$ $A = -\fracLap$ are equivalent:
\begin{enumerate}[(a)]
\item Fourier definition:
\begin{align}
\label{eq:defFracLapF}
(\FT A f)(\xi) = -|\xi|^{\alp} (\FT f)(\xi) 
\end{align}
if $\BSP = L^p(\R^d), p \in [1,2]$.
\item Distributional definition:
\begin{align}
\label{eq:defFracLapD}
\int_{\R^d} A f(y) \phi(y) dy = \int_{\R^d} f(x) A \phi(x) dx
\end{align}
for all $\phi \in \SW(\R^d)$ and with $A \phi$ defined as in \eqref{eq:defFracLapF}.
\item Bochner's definition:
\begin{align}
Af = {1 \over |\Gamma(-{\alp \over 2})|} \int_0^{\infty} (e^{t \lap}f - f) t^{-1 - {\alp \over 2}} dt,
\end{align}
where the integral is Bochner's integral for $\BSP$-valued functions.
\item Balakrishnan's definition:
\begin{align}
Af = {\sin {\alp \pi \over 2} \over \pi} \int_0^{\infty} \lap (s - \lap)^{-1}s^{\alp \over 2 -1} ds,
\end{align}
where the integral is Bochner's integral for $\BSP$-valued functions.
\item Singular integral definition:
\begin{align}
\label{eq:defFracLapS}
Af(x) = \lim_{r \to 0^+} {2^{\alp} \Gamma({d+\alp \over 2}) \over \pi^{d /2} |\Gamma(-{\alp \over 2})|} \int_{\R^d \backslash B(x,r)} {f(x + z) - f(x) \over |z|^{d + \alp}} dz,
\end{align}
with the limit in $\BSP$.
\item Dynkin's definition:
\begin{align}
Af(x) = \lim_{r \to 0^+} {2^{\alp} \Gamma({d+\alp \over 2}) \over \pi^{d /2} |\Gamma(-{\alp \over 2})|} \int_{\R^d \backslash \overline{B}(x,r)} {f(x + z) - f(x) \over |z|^d (|z|^2 - r^2)^{\alp/2}} dz,
\end{align}
with the limit in $\BSP$.
\item Quadratic form or Dirichlet form definition: $\la Af, \phi \ra = \mathcal{E}(f, \phi)$ for all $\phi$ in the Sobolev space $H^{\alp \over 2}$, where
\begin{align}
\mathcal{E}(f,g) = {2^{\alp} \Gamma({d+\alp \over 2}) \over \pi^{d /2} |\Gamma(-{\alp \over 2})|} \int_{\R^d} \int_{\R^d} {\Big(f(x) - f(y)\Big)\Big(\overline {g(x)} - \overline{g(y)}\Big) \over |x - y|^{d + \alp}} dx dy,
\end{align}
when $\BSP = L^2(\R^d)$.
\item Semigroup definition:
\begin{align} 
	A f = \lim_{t \to 0^+} {T_t f - f \over t} ,
\end{align}
with $T_t f = f * p_t$ and $(\FT p_t)(\xi) = e^{-t |\xi|^{\alp}}$.

\item Inverse Riesz potential definition:
\begin{align}
- f(x) = {2^{-\alp} \Gamma({d-\alp \over 2}) \over \pi^{d /2} |\Gamma(-{\alp \over 2})|} \int_{\R^d} { A f(x + z) \over |z|^{d - \alp} } dz
\end{align}
if $\alp < d$ and $\BSP = L^p(\R^d)$, $p \in [1, {d \over \alp})$. 
\item Harmonic extensions' definition:
\begin{align}
\begin{split}
\lap_x u(x,y) + \alp^2 c_{\alp}^{2 \over \alp} y^{2 - 2/\alp} \partial^2_{y} u(x,y) &= 0 \quad for \: y > 0,\\
u(x,0) &= f(x), \\
\partial_y u(x,0) &= A f(x),
\end{split}
\end{align}
where $c_{\alp} = {2^{-\alp} |\Gamma({-\alp \over 2}) | \over \Gamma({\alp \over 2})}$ and $u(\cdot, y)$ is a function of class $\BSP$ which depends continuously on $y \in [0,\infty)$ and $\|u(\cdot, y)\|_{\BSP}$ is bounded in $y \in [0, \infty)$.
\end{enumerate}
In addition, in (c), (e), (f), (h) and (j), the convergence in the uniform norm can be relaxed to pointwise convergence to a function in $\BSP$ when $\BSP = C_0$ or $\BSP = C_b$. Finally, for $\BSP = L^p(\R^d)$, with $p \in [1, \infty)$, the norm convergence in $(e)$, $(f)$, $(h)$ or $(j)$ implies pointwise convergence for almost all $x$.  
\end{thm}
\begin{proof} Refer to [\citenum{Kwasnicki2017}]
\end{proof}
\noindent
We now collect some properties of the operator $(-\lap)^{\alp}$. 

\begin{thm}
\label{thm:domFracLap}
If the IG of the $L^2$-Markov semigroup induced by the \LP $X$ is the fractional Laplacian $(-\lap)^{\alpha}$, for $ \alpha \in (0,1)$, then its domain $D(A)$ is given by $H^{2\alpha}(\R^d) \equiv W^{2\alpha,2}(\R^d)$.
\end{thm}
\begin{proof} Refer to [\citenum{Samko2001,MaRockner2012}]. \end{proof}

\begin{thm} Let $A$ be a sectorial operator of angle $\secAng{A}$ then for $\alp \in (0,1)$ the operator $A^{\alp}$ is a sectorial operator of angle $\alp \secAng{A}$.
\end{thm}
\begin{proof} Refer to [\citenum{YAGI09,Berg1993}]. \end{proof}

\begin{eg} 
\label{eg:secOpFracLap}
Let $A = -\lap$, then we know that $-\lap$ is sectorial operator of angle $0$, thus $(-\lap)^{\alp}$, for $\alp \in (0,1)$,  is also a sectorial operator of angle $0$. 
\end{eg} 
\noindent
Now by letting $\alp$ to be a bounded stochastic process, such that $\alp(\om): [0,T] \to (0,1)$ for every $\om \in \Omega$ fixed, we can specify a stochastic, time dependent Fourier multiplier $\xi^{\alp_t}$. Then due to \eqref{eq:defFracLapD} and Theorem \ref{thm:domFracLap} we  readily observe that, for each fixed $\om$ and $t$, we get  $(A_t(\om), D(A_t(\om))) = ((-\lap)^{\alp_t(\om)}, H^{\alp_t(\om)}(\R^d))$. Additionally, Example \ref{eg:secOpFracLap} yields that for every $\om \in \Om$ and $t \in [0,T]$, $A_t(\om)$ is a sectorial operator of angle $0$. Thus for every $\om$ and $t$ fixed we have that $A_t(\om)$ generates a bounded analytic semigroup of angle $\piT$. 
Thus we are now in a good position to check the existence of the two parameter semigroup for time dependent operator $A_t(\om)$ with $\om \in \Om$ being fixed. 
	To this end we need to check that $(A_t(\om))_{(t \in [0,T])}$ satisfies  Assumption \ref{asm:nonAtAcp}. The only tricky part is to check the condition \eqref{eq:asmFracResLip}, which is established by Lemma \ref{lem:fracOpLipCon} in  Section \ref{sec:fracSpdeMod}.\\
\newcommand{\gam}{\gamma}
\subsubsection{A model for acid mediated cancer invasion \label{sec:fracSpdeMod}}
Based on the above discussion we now propose a macroscopic model for cancer invasion which is mediated by tissue acidity. In this model the main variables of interest are: (i) cancer cell density $C$, (ii) normal cell density $N$, and (iii) the shifted ratio between extracellular and intracellular proton concentrations $H_e$ and $H_i$, respectively, called the proton ratio index or simply proton index $H = H_e/ H_i - 1$. Next we describe the dynamics for each of the aforementioned quantities.\\[1ex]
\noindent
\textbf{Proton dynamics:}
Since we intend to propose just a macroscopic models, the proton dynamics is modeled via the changes in the ratio of intracellular and extracellular protons. Cancer cells exhibit a reversed pH gradient x, where $H_e$ is greater than $H_i$. This relationship is captured by the shift, by a constant of value 1, in the definition of $H$. The shift by a constant $1$ is just for modeling convenience, since positivity of $H$ implies that $H_e$ is greater than $H_i$. Based on this the dynamics of $H$ is governed by spatial spread due to random motion and directed flow towards higher cancer cell density. The  former is modeled by the diffusion operator $\lap$, while the later is modeled via the advection term. Apart from this, the dynamics is also influenced by proton production and loss, which is abstractly modeled by a logistic term $R_1$. This accounts for the consistent accumulation of acid in the interstitial region of the cells. Finally, the dynamics is also influenced by random fluctuations which is captured by a spatial Brownian motion. \\[1ex]
\noindent
\textbf{Cancer cell dynamics:}
According to  [\citenum{VISAF96,VISRAPLU08,MDPMO14,REY15}] many animals exhibit L\'{e}vy flight type movements in order to gather food or to avoid predators.   
As discussed at the end of the previous subsection, L\'{e}vy flight processes can be constructed using a CTRW process and it generates a fractional Laplacian operator at the macroscopic level. Moreover, because we are interested in modeling diverse kinds of movements at the macroscopic level, in order to describe the movement of cancer cells through fibrous tissue and crowded cells, we use a time dependent random fractional Laplacian operator $-(-\lap)^{\alp_t}$. The exponent $\alp_t$ is a bounded stochastic process which is a function of the proton index $H$. Apart from this we also model movement of cells via the pH-taxis and haptotaxis mechanisms. The former models the movement of cells in the direction of a high proton index $H$ while the latter models the movement of cells away from high normal cell (or tissue) density. Apart from spatial movement, the dynamics of cancer cells is also governed by proliferation, which is modeled by a logistic term. \\[1ex]
\noindent
\textbf{Normal cell dynamics:}
The dynamics of normal cell density $N$ is simply given by an ODE which describes the decay of healthy cells due to their interaction with cancer cells and protons. 

\begin{subequations}
\label{eq:fracSpdeMod}
\begin{align}
dH_t &= \sigma_H A H_t dt + R_1(H_t) dt +  f(H_t,C_t) \grad H_t \cdot \grad C_t dt + H_t dW_t \label{eq:fracSpdeModH}\\[2ex]
dC_t  &= \sigma_C A_t C_t dt + R_2(C_t) dt + \grad \cdot (g(H_t,C_t) \grad N_t) dt - \grad \cdot(h(H_t,C_t) \grad H_t) dt \label{eq:fracSpdeModC}\\[2ex]
dN_t &= R_3(C_t, H_t) dt. \label{eq:fracSpdeModN}
\end{align}
\end{subequations}


\vspace{-.5cm}
\begin{align*}
A := \lap \quad \text{ and } \quad A_t = -(-A)^{\alpha_t}, \quad \alpha_t \in [a_1, a_2] \subset (0,1),
\end{align*}
where  $\alpha_t$ is an $\R$ valued, bounded and continuous (w.r.t time $t$)  process such that:
\begin{align*}
\sup_{t,x,\om} \alpha_t = a_2 \in (0,1), \quad
\inf_{t,x,\om} \alpha_t = a_1 \in (0,1),
\end{align*}
and ${1 \over 2} < a_1 < a_2 < 1$. 
Let $\gamma \in (0,1)$ be such that \[ \gamma a_2 = \eta \in (0,1), \quad \gamma a_1 = \delta \in (0,1)\] and $\delta < \eta$. Based on this, we let $\alp_t$ take the following form:
\begin{align}
\label{eq:alpProc}
\alpha_t := a_1 + (a_2-a_1) {a H_t \over 1 + a H_t},  \quad a, a_1, a_2 \in \R \text{ and } a_1 < a_2 . 
\end{align}
Since $a_1 > {1 \over 2}$ there exists $\nu > 0$ such that $(1 - \nu) a_1 > {1 \over 2}$. 
\begin{lem} The operator $(-A_t , D(-A_t)) = ((-A)^{\alpha_t}, H^{\alpha_t}(\R^d))$ is a closed and densely defined operator in $L^2(\R^n)$.
\end{lem}
\begin{lem} 
\label{lem:fracOpLipCon1}
For $ u \in H^{\eta}(\R^d)$, $\|(A_t)^{\gamma} u\|_{ L^2(\R^d)} \le \ubd_{_{A^{\eta}}} $ uniformly in $t$, where $\ubd_{_{A^{\eta}}} < \infty$ is a fixed constant. 
\end{lem}
\begin{proof} Firstly, we note that $(-A_t)^{\gamma} = (-A)^{\gamma \alpha_t}$, so $D(A_t^{\gamma}) = H^{\gamma \alpha_t}$. Since $H^{s_1}(\R^d) \subseteq H^{s_2}(\R^d)$ for $s_1 \ge s_2$ we have that $H^{\eta}(\R^d) \subseteq D((-A_t)^{\gamma})$ i.e. $H^{\eta}(\R^d) \subseteq H^{\gamma \alpha_t}(\R^d)$ (since $\gamma \alpha_t \in [\delta, \eta] \subset (0,1)$ and continuous w.r.t,  $t$). 
\begin{align*}
\sup_{u \in H^{\eta}(\R^d), \|u\| =1} \|(-A_t)^{\gam} u\|  \le \sup_{u \in D((-A_t)^{\gam}), \|u\|=1} \|(-A_t)^{\gam} u\|  = \|(-A_t)^{\gam}\|.
\end{align*}
Thus
\begin{align}
\label{eq:fracOpLipCon1}
{\underset{\om \in \Om}{{{\sup_{t \in \R,}}}}} {\underset{ \|u\| =1}{\sup_{u \in H^{\eta}(\R^d),}}} \|(-A_t)^{\gam} u\| \le {\underset{\om \in \Om}{{{\sup_{t \in \R,}}}}} \|(-A_t)^{\gam}\| < \ubd_{_{A^{\eta}}},
\end{align}
with $\ubd_{_{A^{\eta}}} < \infty$ being a fixed constant.
\end{proof}
\begin{lem} 
\label{lem:fracOpLipCon2}
For a bounded process $\alpha_t$ defined as in \eqref{eq:alpProc}, we have that 
\begin{align}
\label{eq:fracOpLipCon2}
\|[(-A_t)^{-1} - (-A_t)^{-1}] u\|_{H^{\eta}(\R^d)} \le \const(\eta,\delta) \|\alpha_t - \alpha_s\|_{H^2(\R^d)} \|u\|_{L^2(\R^d)}, \quad \forall \: u \in L^2(\R^d).
\end{align}
\end{lem}
\begin{proof}
For $-\alpha_t$, due to Equation 2.108 of [\citenum{YAGI09}, we have that:
\begin{align*}
(-A_t)^{-1} = {-\sin(\pi \alpha_t) \over \pi} \int_0^{\infty} \rho^{-\alpha_t} (\rho + A)^{-1} d\rho .
\end{align*}
\begin{align*}
\|[ (-A_t)^{-1} - (-A_s)^{-1}] u \|_{H^{\eta}} &=  \Big\| {\sin(\pi \alpha_t) \over \pi} \int_0^{\infty} \rho^{-\alpha_t} (\rho + A)^{-1} u d\rho \\
&\hspace*{2cm} - {\sin(\pi \alpha_s) \over \pi} \int_0^{\infty} \rho^{-\alpha_s} (\rho + A)^{-1} u d\rho \Big\|_{H^{\eta}} \\
&\hspace*{-2cm}= {1 \over  \pi} \Big \| \sin(\pi \alpha_t) \int_0^{\infty} \rho^{-\alpha_t} (\rho + A)^{-1} u d\rho  - \sin(\pi \alpha_s) \int_0^{\infty} \rho^{-\alpha_s} (\rho + A)^{-1} u d\rho \Big \|_{H^{\eta}} \\
&\hspace*{-2cm}\le {1 \over \pi} \Bigg( \int_0^{\infty} \Big\| \Big(\sin(\pi \alpha_t)  - \sin(\pi \alpha_s) \Big) \rho^{-\alpha_t} (\rho + A)^{-1} u \Big\|_{H^{\eta}} d\rho \\
%
%
&\hspace*{1cm} +  | \sin(\pi \alpha_s)| \int_0^{\infty}\Big\| |\rho^{-\alpha_t} - \rho^{-\alpha_s}|   (\rho + A)^{-1} u \Big\|_{H^{\eta}} d\rho \Bigg) \\
&\hspace*{-2cm}\le \pi  \Bigg(\int_0^{\infty} \rho^{-\alpha_t} \Big\| |\alpha_t - \alpha_s| (\rho + A)^{-1} u \Big\|_{H^{\eta}} d\rho \\
&\hspace{2cm}+ \int_0^{\infty} |\ln(\rho) \rho^{-\alpha_{\cdot}}| \Big\| | \alpha_t - \alpha_s|  (\rho + A)^{-1} u \Big\|_{H^{\eta}} d\rho \Bigg)\\
&\hspace*{-2cm}\le \pi \Bigg(\int_0^{\infty} \rho^{-\alpha_t} \|\alpha_t - \alpha_s \|_{L^{\infty}} \Big\|(\rho + A)^{-1} u \Big\|_{H^{\eta}} d\rho \\
&\hspace*{2cm}\le \int_0^{\infty} |\ln(\rho) \rho^{-\alpha_{\cdot}}| \|\alpha_t - \alpha_s\|_{L^{\infty}} \Big\| (\rho + A)^{-1} u \Big\|_{H^{\eta}} d\rho \Bigg)\\
\end{align*}
This implies that
\begin{align*}
\|[ (-A_t)^{-1} - (-A_s)^{-1}] u \|_{H^{\eta}} &\le { 2 \pi } \|\alpha_t - \alpha_s\|_{H^2}  \int_0^{\infty} (\rho^{-\eta}  + |\ln(\rho)|\rho^{-\eta}) \|(\rho + A)^{-1} u\|_{H^{\eta}} d\rho \\
&\hspace*{-2cm}\le { 2 \pi}  \Bigg ( \int_0^{\infty} {2\rho^{-\eta}  \over (1 + \rho) } d\rho + \int_0^{1} { \ln(\rho)\rho^{-\eta} \over (1 + \rho)}  d\rho + \int_1^{\infty} { \ln(\rho)\rho^{-\eta} \over (1 + \rho)}  d\rho \Bigg) \|\alpha_t - \alpha_s\|_{H^2} \|u\|_{L^2} \\
&\hspace*{-2cm}\le { 2 \pi}  \Bigg ( \int_0^1 2 \rho^{(1-\eta)-1} (1-\rho)^{1-1} d\rho  + \int_1^{\infty} 2 \rho^{-\eta-1} d\rho + \int_{-\infty}^{0} {z e^{(1-\eta)z} \over 1 + e^z} dz \\
&\hspace{2cm} + \int_{0}^{\infty} {z e^{(1-\eta)z} \over 1 + e^z} dz \Bigg)  \|\alpha_t - \alpha_s\|_{H^2} \|u\|_{L^2}\\
\end{align*}
This further implies that
\begin{align*}
\|[ (-A_t)^{-1} - (-A_s)^{-1}] u \|_{H^{\eta}} &\le { 2 \pi}  \Bigg ( \beta(1-\eta,1) + {2 \over \eta} +  \int_{0}^{\infty} {-z e^{-(1-\eta)z} \over 1 + e^{-z}} dz + \int_{0}^{\infty} {z e^{(1-\eta)z} \over 1 + e^z} dz \Bigg) \\
&\hspace*{2cm} \|\alpha_t - \alpha_s\|_{H^2} \|u\|_{L^2} \\
&\hspace*{-2.5cm}\le { 2 \pi} \Bigg ( \beta(1-\eta,1) + {2 \over \delta} +  |\int_{0}^{\infty} {z e^{-(1-\eta)z} \over 1 + e^{-z}} dz | + \int_{0}^{\infty} {z e^{(-\eta)z} \over 1 + e^{-z}} dz \Bigg)  \|\alpha_t - \alpha_s\|_{H^2} \|u\|_{L^2} \\
&\hspace*{-2.5cm}\le { 2 \pi} \Bigg ( \beta(1-\eta,1) + {2 \over \delta} + \int_{0}^{\infty} {z e^{-(1-\eta)z} } dz |+ \int_{0}^{\infty} {\eta z e^{(-\eta)z} } dz \Bigg)  \|\alpha_t - \alpha_s\|_{H^2} \|u\|_{L^2} \\
&\hspace*{-2.5cm}\le { 2 \pi} \Big( \beta(1-\eta,1) + {2 \over \delta} + {1 \over (1 - \eta)^2} + {1 \over \eta^2} \Big)  \|\alpha_t - \alpha_s\|_{H^2} \|u\|_{L^2} \\
&\hspace*{-2.5cm}\le \const(\eta,\delta)  \|\alpha_t - \alpha_s\|_{H^2} \|u\|_{L^2}.
\end{align*}
\end{proof}

\begin{lem} 
\label{lem:fracOpLipCon}
Let $H \in C^{\mu}([0,T];H^2(\R^d))$ and $\alp_t$ as defined in \eqref{eq:alpProc}. Then we get that 
\begin{align}
\label{eq:fracOpLipCon}
\|(-A_t)^{-\gam}[ (-A_t)^{-1} - (-A_s)^{-1}] u\|_{L^2(\R^d)} \le C |t - s|^{\mu} \|u\|_{L^2(\R^d)}
\end{align}
\end{lem}
\begin{proof} Because $A_t^{-1}$ maps $D(A^{\eta})$ into $L^2(\R^d)$ and $D(A_t^{-\gam}) \subset D(A^{\eta})$ for all $t \in [0,T]$, we are only left to show that $\|\alp_t - \alp_s\|_{H^2(\R^d)}$ is H\:older continuous, uniformly with respect to $t$ and $s$. This follows easily due to \eqref{eq:alpProc} and $H \in C^{\mu}([0,T];H^2(\R^d))$, as shown below
\begin{align*}
\|\alp_t - \alp_s \|_{H^{\eta}} &\le \const \Bigg\| {H_t   - H_s \over (1 + a H_t)(1 + a H_s) } \Bigg\| \\
&\le \const \|H_t - H_s\| _{H^{\eta}}\\
&\le \const \|H_t - H_s\| _{H^2} \\
&\le \const |t - s|^{\mu}, \quad \text{ since $H \in C^{\mu}([0,T];H^2(\R^d))$.}
\end{align*}
Thus due to \eqref{eq:fracOpLipCon1} and \eqref{eq:fracOpLipCon2} we get the required result.
\end{proof}
\noindent
Thus we have verified that $A_t$ and $A$ satisfy Assumption  \ref{asm:nonAtAcp} and consequently generate the corresponding semigroup. Thus for sufficiently regular data, following the lines of [\citenum{YAGI09}], one can prove the local existence of a mild solution to \eqref{eq:fracSpdeMod}.\\

\textbf{Properties of the coefficient functions $g, f, h$}
\begin{enumerate}
\item $g \in C^2(\R \times \R; \R)$. There exist constants $\lbd_g, \ubd_g > 0$ such that $\lbd_g \le g + g' + g''\le \ubd_g$.
\item $f \in C^2(\R \times \R; \R)$. There exist  constants $\lbd_f, \ubd_f > 0$ such that $\lbd_f \le f + f' + f''  \le \ubd_f$.
\item $h \in C^2(\R \times \R; \R)$. There exist constants $\lbd_h, \ubd_h > 0$ such that $\lbd_h \le h + h' + h'' \le \ubd_h$.
\item $g, f, h$ such that $C, H$ and $N$ are in $ L^{12}(\dom)$ and $\grad C, \grad N$ and $\grad H$ are in $L^4(\dom)$.
\item $N \in H^2_p$, $p \ge {12 \over 5}$. 
\item $R_1$, $R_2$ and $R_3$ are of sub-polynomial growth. 
\end{enumerate}

\subsection{Numerical simulations}\label{sec:simulations}
In this section we perform numerical simulations for the proposed model \eqref{eq:fracSpdeMod}.
We assume periodic boundary conditions and restrict ourselves to the finite 2D domain $\dom$. Let $L^2_p(\dom)$ denote the subspace of $L^2$ functions on $\dom$ that are not constants and are even. 
We let $\dom := [a_1, b_1)\times[a_2,b_2)$;  $a_1 < b_1\, a_2 < b_2$; $a_1\,a_2\,b_1\,b_2 \in \mathbb{R}$. Let $M_x \, M_y$ be the number of partitions of the $x-$axis and the $y-$axis of $\dom$. Thus we have $\delta_x := \frac{b_1-a_1}{M_x} \, \delta_y := \frac{b_2- a_2}{M_y}$. 
\noindent
The spatial grid points at which the solution to our problem will be computed are represented as $(x_{k} \, y_{j} )_{k,j} \, k \in \{0\,\dots\,M_x\} \, j \in \{0\,\dots\,M_y\}\, x_{k} := k \delta_x\, y_{j} := j \delta_y$. Also, let the time interval $I := [0 , T] \, T \in \mathbb{R}^+$ be divided into a number of $N_{\tau}$ subintervals with $\tau := \frac{|I|}{N_\tau}$ and the temporal grid points $(t_n)_{n \in \{0,\dots,N_{\tau}\}}$ with $t_n := n \tau$.
\subsubsection{Simulating the $Q$-Wiener process:}
Let $Q$ be a bounded linear operator on $L^2(\dom)$ such that
\[ Q e_n = \lam_n e_n,\] where $(e_n)_{n \in \N}$ is the orthonormal basis of $L^2(\dom)$. The eigenvalues $(\lam_n)_{n \in \N}$  of $Q$ are  such that 
\[ \sum_{k = 1}^{\infty} \lam_k^2 < \infty. \]
Let $L_x = b_1 - a_1$, and $L_y = b_2 - a_2$. Let 
\[ \Bigg(\sqrt{{2 \over L_x}} \cos\Big({2\pi n x \over L_x}\Big)\Bigg)_{n\in\N} \quad \text{and} \quad \Bigg(\sqrt{{2 \over L_y}} \cos \Big({2\pi m x \over L_y}\Big)\Bigg)_{m \in \N} \] be the ONB of $L^2_p([a_1,b_1))$ and $L^2_p([a_2,b_2))$, respectively.  Now if we let $\lam_k = {1 \over (1 + k)}$ for $k \in \{m, n\}$ with $m$, $n \in \N$. Then by Proposition 2.1.6 and Proposition 2.1.10 of [\citenum{Prevot2007}], the $Q$-Wiener process on $L^2_p([a_1,b_1))$ and $L^2_p([a_2,b_2))$ can be represented as
\[ w_t(x) = \sum_{n=1}^{\infty} \xi_n(t) \lam_n e_n(x) \] and
\[ w_t(y) = \sum_{m=1}^{\infty} \xi_m(t) \lam_m e_m(y), \] where
$(\xi_k)_{k \in \N}$ for $k = m$ or $k = n$ is an sequence of independent real valued Brownian motions. Since $L^2_p([a_1,b_1)\times[a_2,b_2))$ is isomorphic to $L^2_p([a_1,b_1)) \times L^2_p([a_2,b_2))$, we get that 
\[ W_t(x,y) = \sum_{m,n=1}^{\infty} \xi_{m,n}(t) \lam_{m,n} e_n(x) \otimes e_m(y) \]  which, for $\lam_{m,n} = \lam_m \lam_n$, is equivalent (in the distributional sense) to $w_t(x) \times w_t(y)$, i.e. 
\[ W_t(x,y) \overset{\text{in Law}}{=} w_t(x) \times w_t(y) = \sum_{n=1}^{\infty} \xi_n(t) \lam_n e_n(x)  \times  \sum_{m=1}^{\infty} \xi_m(t) \lam_m e_m(y). \]
Thus the increments of the $Q$-Wiener process are given as 
\[ dW_t(x,y) = \sum_{m,n=1}^{\infty} d\xi_{m,n}(t) \lam_{m,n} e_n(x) \otimes e_m(y). \] Letting $W^n_{k,j} = W_{t_n}(x_k,y_j)$ we get the following discretization
\begin{align}
\label{eq:discQBM}
dW^n_{k,j} = dW_{k,j} = \sum_{m,n=1}^{\infty} \sqrt{\tau} z_{m,n} \lam_{m,n} e_n(x_k) \otimes e_m(y_j),
\end{align}
where $(z_{m,n})_{m,n \in \N}$ is a sequence of independent random variables having standard Gaussian distribution. 

\subsubsection{Simulating the fractional Laplacian \label{sec:simFracLap}}
For numerical simulation, we use the singular representation of $(-\lap)^{\alp\over 2}$ for $\alp \in (0,2)$, which, based on \eqref{eq:defFracLapS}, is given as 

\begin{align}
\label{eq:fracLapInt}
I f (x) =  (-\lap)^{\alp \over 2}f (x) = \lim_{r \to 0^+} {-2^{\alp} \Gamma({d+\alp \over 2}) \over \pi^{d /2} |\Gamma(-{\alp \over 2})|} \int_{\R^d \backslash B(x,r)} {f(x + z) - f(x) \over |z|^{d + \alp}} dz 
\end{align}
Letting 
\[ c_{d,\alp} := {-2^{\alp} \Gamma({d+\alp \over 2}) \over \pi^{d /2} |\Gamma(-{\alp \over 2})|} = {\alp 2^{\alp-1} \Gamma({d+\alp \over 2}) \over \pi^{d /2} \Gamma(-{2 - \alp \over 2})} \] and $h < 1$,  following [\citenum{Huang2014}], the integral $I$ in \eqref{eq:fracLapInt} can be split into the singular part and the tail part as follows:

\begin{align}
I_s f(x) &= c_{d,\alp} \int_{|y|\le h} {f(x) - f(x-y) \over |x-y|^{d + \alp}} dy  \label{eq:singInt}\\
I_t f(x) &= c_{d,\alp} \int_{|y| \ge h} {f(x) - f(x-y) \over |x-y|^{d + \alp}} dy. \label{eq:tailInt}
\end{align}
The singular integral \eqref{eq:singInt} is implemented using the finite difference scheme described in [\citenum{Huang2014}], which for $d = 1$ is given as

\begin{align}
\label{eq:discSingInt}
I^h_s f(x_k) = {-c_{1,\alp} \over (2-\alp)}\Big[ {f(x_{k+1}) - 2 f(x_k) + f(x_{k-1}) \over \delta_x^{\alp}}\Big] + \mathcal{O}(\delta_x^{4-\alp}).
\end{align}
The tail integral \eqref{eq:tailInt} is implemented using a simple quadrature rule, which for $d=1$ is as follows:
\begin{align}
\label{eq:discTailInt}
I^h_t f(x_k) = c_{1,\alp} \sum_{i=1}^{N} {f(x_k) - f(x_{k +i}) \over (i h)^{\alp}}, \quad \text{for $N<\infty$ large enough.}
\end{align}
\noindent
Based on this we have the following discretization scheme for the system \eqref{eq:fracSpdeMod}. 

\subsubsection{Discretization scheme}
For the normal cell density equation \eqref{eq:fracSpdeModN}, we use the following straightforward discretization:

\begin{align*}
N^{n+1}_{k,j}= N^n_{k,j} + \tau R_3(C^n_{k,j},H^n_{k,j}),
\end{align*}
where, $N^{n+1}_{k,j} := N_{t_{n+1}}(x_k,y_j)$ , $C^{n+1}_{k,j} := C_{t_{n+1}}(x_k,y_j)$, $H^{n+1}_{k,j} := H_{t_{n+1}}(x_k,y_j)$ and $R_3(C^n_{k,j},H^n_{k,j}) = \gam_3 (C^n_{k,j}+H^n_{k,j}) N^n_{k,j}$. \\
\noindent
For discretizing the proton equation \eqref{eq:fracSpdeModH} we use an implicit finite difference scheme for the advection and diffusion terms, while an explicit scheme is employed for the reaction term. The noise term is discretized using \eqref{eq:discQBM}. After letting
\begin{align*}
\begin{array}{c c}
R^{n}_{1;k,j} := \gamma_1 H^n_{k,j} (1 - H^n_{k,j}), & f_{k,j} :=  {C^n_{k,j} \over 1 + C^n_{k,j}}
\end{array}
\end{align*}
we get
\begin{align*}
H^{n+1}_{k,j} &= H^n_{k,j} + \tau R^n_{1;k,j} +  H^n_{k,j} dW_{k,j} \\
&\hspace{1cm} + {\tau \sigma_H  \over \delta^2_x} ( H^{n+1}_{k-1,j} + H^{n+1}_{k+1,j} - 2 H^{n+1}_{k,j}) + {\tau \sigma_H  \over \delta^2_y} ( H^{n+1}_{k,j-1} + H^{n+1}_{k,j+1} - 2 H^{n+1}_{k,j}) \\
&\hspace{1cm} + {\gamma_f \tau \over 4 \delta^2_x}[ f^{n}_{k,j}(H^{n+1}_{k,j+1} - H^{n+1}_{k,j-1})(C^{n}_{k,j+1} - C^{n}_{k,j-1})]\\
&\hspace{1cm} + {\gamma_f \tau \over 4 \delta^2_y}[ f^{n}_{k,j}(H^{n+1}_{k,j+1} - H^{n+1}_{k,j-1})(C^{n}_{k,j+1} - C^{n}_{k,j-1})]
\end{align*}
\noindent
For the cancer cell population dynamics equation \eqref{eq:fracSpdeModC} we use \eqref{eq:discSingInt} and \eqref{eq:discTailInt} to discretize the fractional Laplacian implicitly, while the taxis terms and the reaction term are discretized explicitly. Altogether, after letting 

\begin{align*}
\begin{array}{c c c}
R^{n+1}_{2;k,j} := \gamma_2 C^n_{k,j}(1- C^n_{k,j}), &
g^{n+1}_{1;k,j} := \gamma_g {N^{n+1}_{k,j} C^n_{k,j} \over 1 + (C^n_{k,j} + N^n_{k,j})^2}, &
h^{n+1}_{2,k,j} := \gamma_h {H^{n+1}_{k,j} C^n_{k,j} \over 1 + (C^n_{k,j} + H^n_{k,j})^2}.
\end{array}
\end{align*}
we get
\begin{align*}
C^{n+1}_{k,j} &= C^n_{e,k,j} + \tau R^{n+1}_{2;k,j} + I^h_t(x_k,y_j) + \\
&\hspace{1cm} + {\tau \sigma_C c_{2,\alp} \over (2-\alp) \delta^{\alpha}_x} ( C^{n+1}_{k-1,j} +  C^{n+1}_{k+1,j} - 2 C^{n+1}_{k,j})  + {\tau  \sigma_C c_{2,\alp} \over (2-\alp) \delta^{\alpha}_y} (  C^{n+1}_{k,j-1} + C^{n+1}_{k,j+1} - 2  C^{n+1}_{k,j}) \\
&\hspace{1cm} + {\gamma_g \tau \over 2 \delta^2_x}[( g^{n+1}_{k,j-1} + g^{n+1}_{k,j})(H^n_{k,j-1} - H^n_{k,j}) + ( g^{n+1}_{k,j+1} + g^{n+1}_{k,j})(H^n_{k,j+1} - H^n_{k,j}) ] \\
&\hspace{1cm} + {\gamma_g \tau \over 2 \delta^2_y}[( g^{n+1}_{k,j-1} + g^{n+1}_{k,j})(H^n_{k,j-1} - H^n_{k,j}) + ( g^{n+1}_{k,j+1} + g^{n+1}_{k,j})(H^n_{k,j+1} - H^n_{k,j}) ] \\
&\hspace{1cm} + {\gamma_h \tau \over 2 \delta^2_x}[( h^{n}_{k,j-1} + h^{n}_{k,j})(N^n_{k,j-1} - N^n_{k,j}) + ( h^{n}_{k,j+1} + h^{n}_{k,j})(N^n_{k,j+1} - N^n_{k,j}) ] \\
&\hspace{1cm} + {\gamma_h \tau \over 2 \delta^2_y}[( h^{n}_{k,j-1} + h^{n}_{k,j})(N^n_{k,j-1} - N^n_{k,j}) + ( h^{n}_{k,j+1} + h^{n}_{k,j})(N^n_{k,j+1} - N^n_{k,j}) ]
\end{align*}

\subsubsection{Numerical simulation of the system \eqref{eq:fracSpdeMod}}
In this section we provide the results obtained from the above numerical discretization scheme. In order to highlight the effectiveness of the model we illustrate the simulations just for the 2D case. 

\begin{table}[h!]
\begin{center}
        \addtolength{\tabcolsep}{-4pt}
        \caption{ \small{ Simulation parameters } \label{tab:simparFracSPDE}}
        \begin{center}
        \begin{tabular}{|l|l|l|}
        \hline
        \multicolumn{3}{|l|}{\qquad \qquad \small{ Numerical parameters} }\\ 
        \hline
        & \quad \small{ Phenomenological relevance } &   \\\hline
        N & {\small{\# time steps}}& 150 \\\hline
        M & {\small{\#  Monte Carlo simulations}}& 500 \\\hline
        $\tau$ & {\small{temporal step size}} & .1 \\\hline
        $h_{x_1}$ & {\small{spatial step size along $x_1$}}& .1 \\\hline 
        $h_{x_2}$ & {\small{spatial step size along $x_2$}}& .1\\\hline 
        $N_{x_1}$ & {\small{grid resolution along $x_1$}}&  21\\\hline
        $N_{x_2}$ & {\small{grid resolution along $x_2$}}&  21\\\hline
        \end{tabular}
        \end{center}
    \end{center}
\end{table}
\begin{table}[h!]
\begin{center}
        \addtolength{\tabcolsep}{-4pt}
        \caption{\small{ Model parameters }\label{tab:modParFracSPDE}}
        \begin{tabular}{|l|l|l|}
        \hline
        \multicolumn{3}{|l|}{\qquad \qquad \small{ Growth and decay parameters $\Xi_R$}}\\\hline
        & \qquad \qquad Phenomenological relevance &   \\\hline
        $\gamma_{1}$& growth rate for cancer cells& .005\\\hline
        $\gamma_{2}$& growth rate for proton index value &  .05\\\hline
        $\gamma_{3}$& decay rate of normal cells& .015\\\hline
        $ \sigma_{W}$ & intensity of noise for the proton index & .131 \\\hline
        \hline
        \multicolumn{3}{|l|}{\qquad \qquad \small{ Migration parameters $\Xi_M$}}\\\hline
        &  \qquad \qquad \small{ Phenomenological relevance } &  \\\hline
        $\sigma_{_H}$& diffusion coefficient for the proton index& .0008\\\hline
        $\gamma_{C}$& diffusion coefficient for cancer cells &  .00035\\\hline
        $\gamma_{g}$& speed of cell movement due to hapto-taxis &  .0007\\\hline
        $\gamma_{h}$& speed of cell movement due to chemo-taxis & .0037\\\hline
        $\gamma_{f}$& speed of advection for the proton index &  .0082\\\hline
        \end{tabular}
\end{center}
\end{table}


\begin{figure}[h!]
	\centering
	\includegraphics[trim = 0cm 0cm 0cm 0cm,clip,scale=.175]{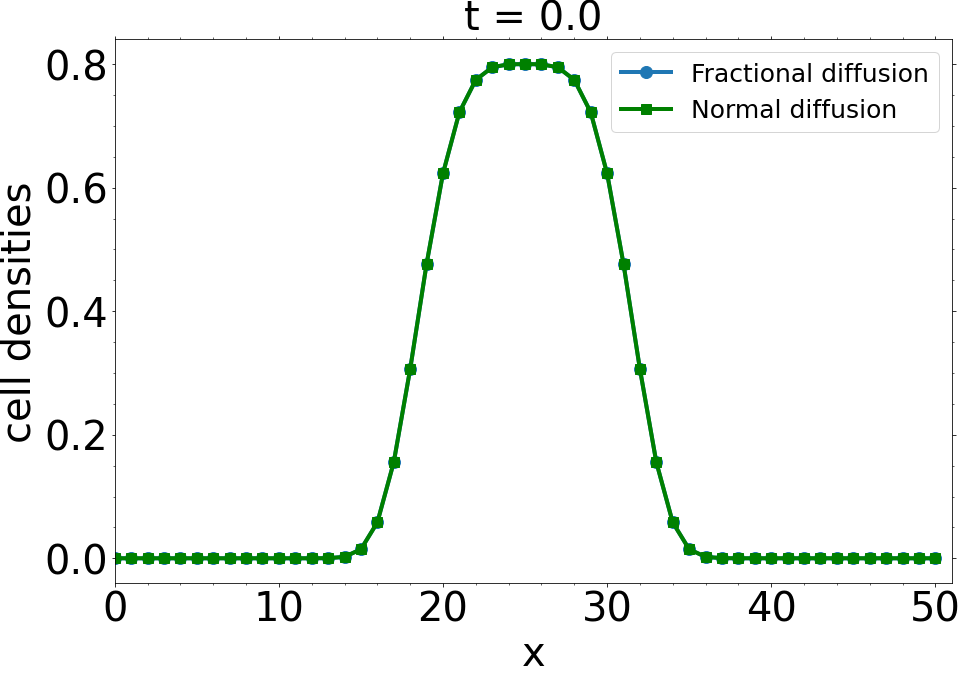}\qquad
	\includegraphics[trim = 0cm 0cm 0cm 0cm,clip,scale=.175]{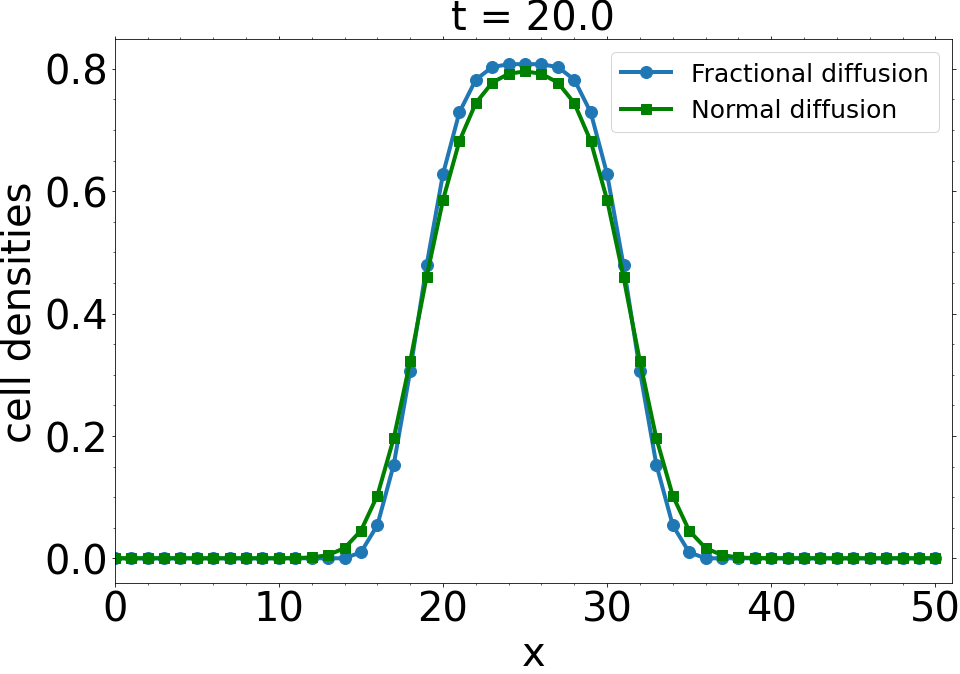}\\
	\includegraphics[trim = 0cm 0cm 0cm 0cm,clip,scale=.175]{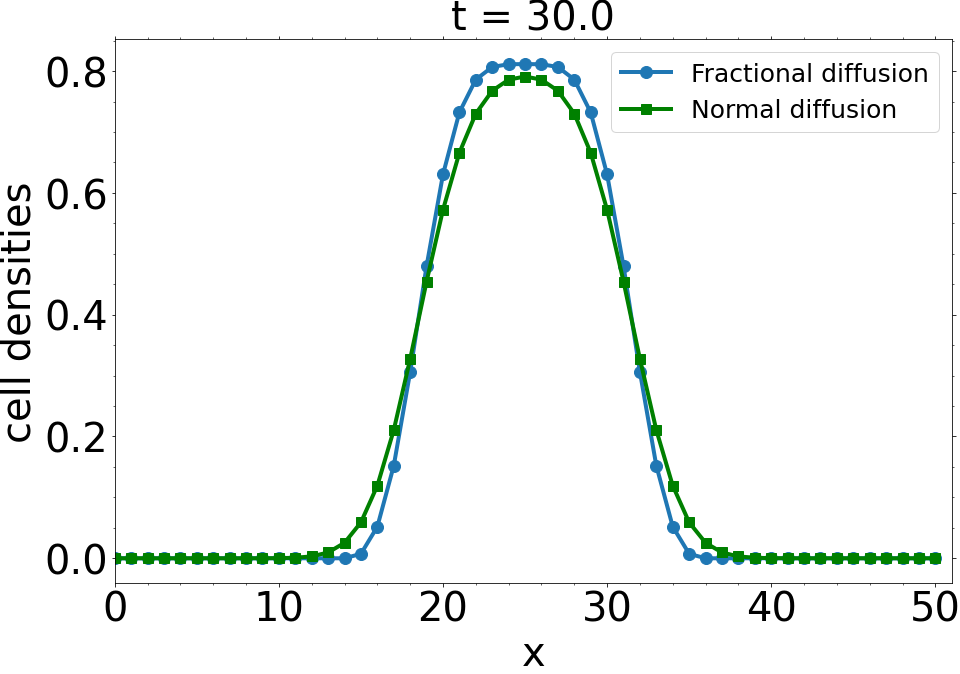}\qquad
	\includegraphics[trim = 0cm 0cm 0cm 0cm,clip,scale=.175]{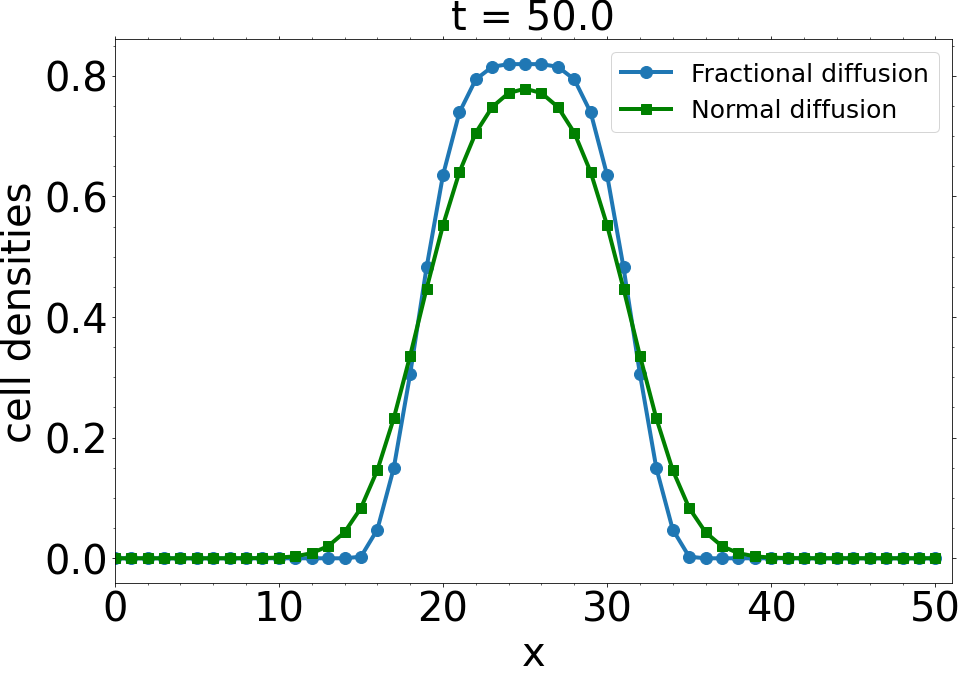}
\caption{Time snapshots of solutions to a fractional diffusion and a standard diffusion equation. The blue curve with dots represents the solution to the fractional diffusion equation, while the green curve with bars represents the solution to the standard diffusion equation. \label{fig:simFracVsStdDiff}}
\end{figure}

\begin{figure}
	\centering
	\includegraphics[trim = 0cm 0cm 0cm 0cm,clip,scale=.225]{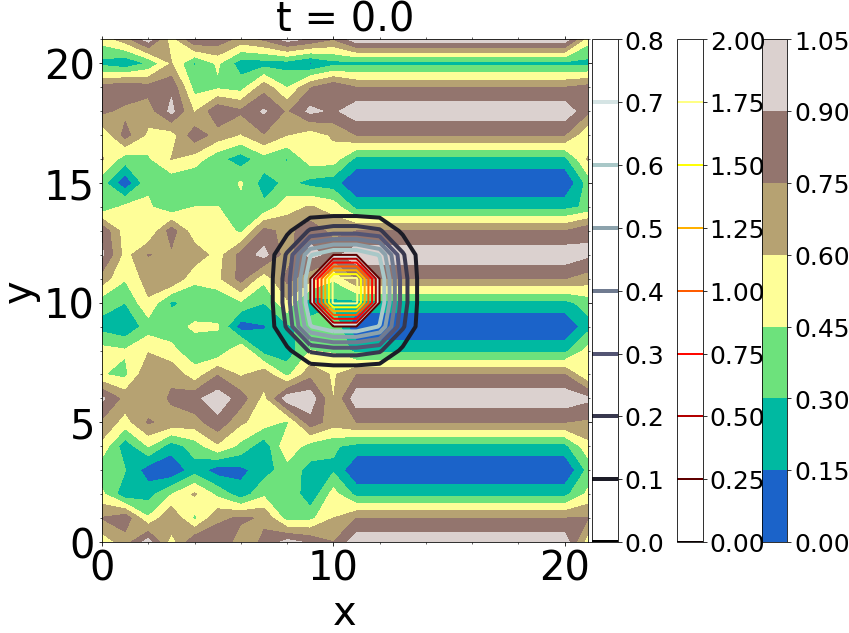}
	\includegraphics[trim = 0cm 0cm 0cm 0cm,clip,scale=.225]{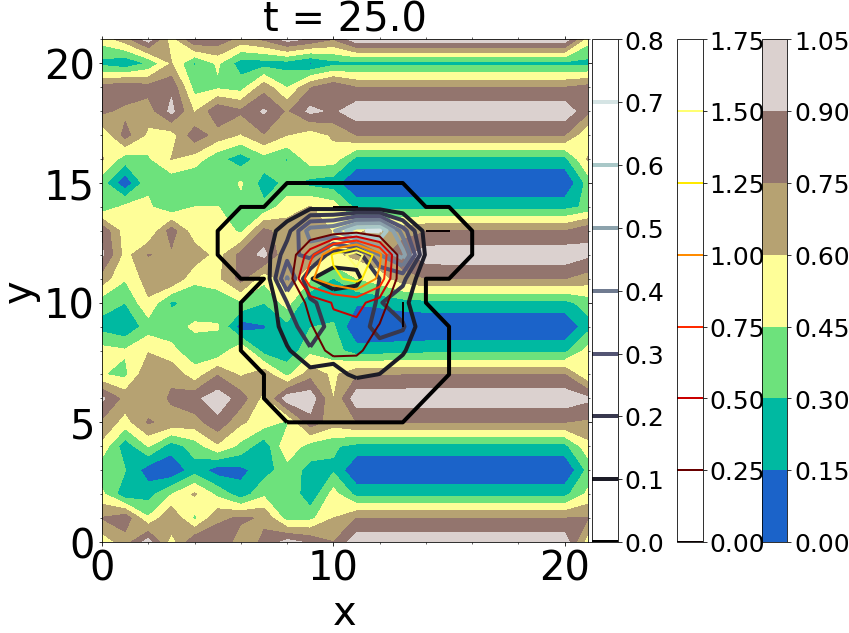}
	\includegraphics[trim = 0cm 0cm 0cm 0cm,clip,scale=.225]{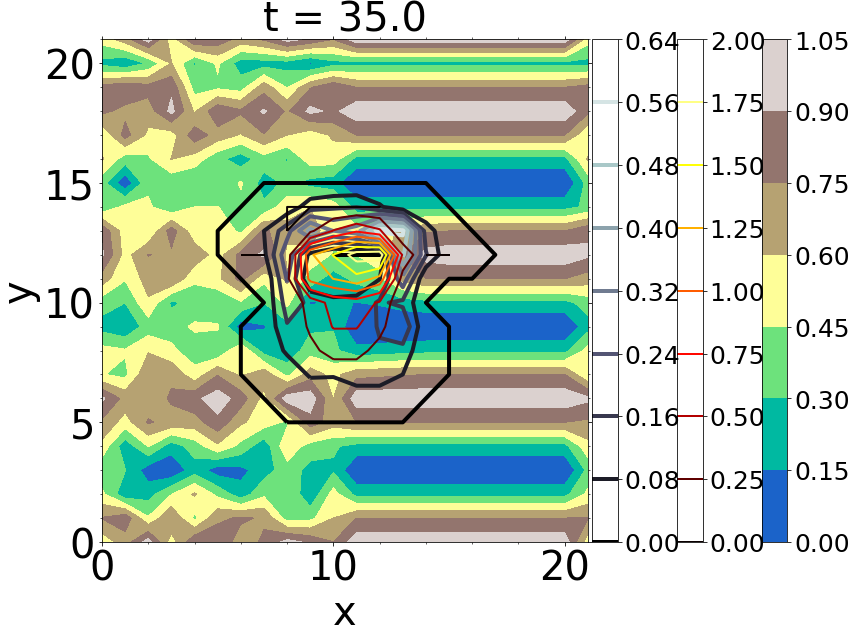}
	\includegraphics[trim = 0cm 0cm 0cm 0cm,clip,scale=.225]{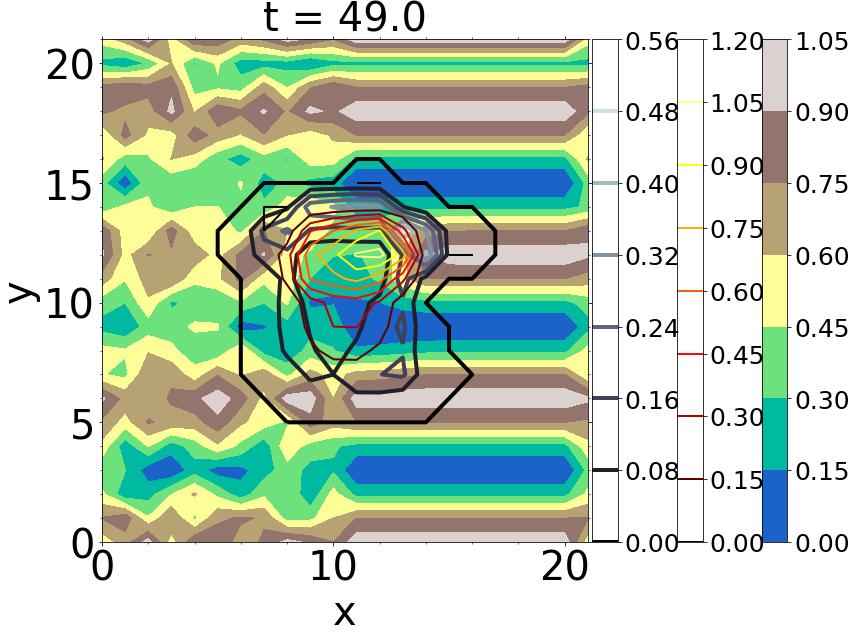}
\caption{Time snapshots of the sample solution number 10.  \label{fig:sim2DMC1}}
\end{figure}
\begin{figure}
	\centering
	\includegraphics[trim = 0cm 0cm 0cm 0cm,clip,scale=.225]{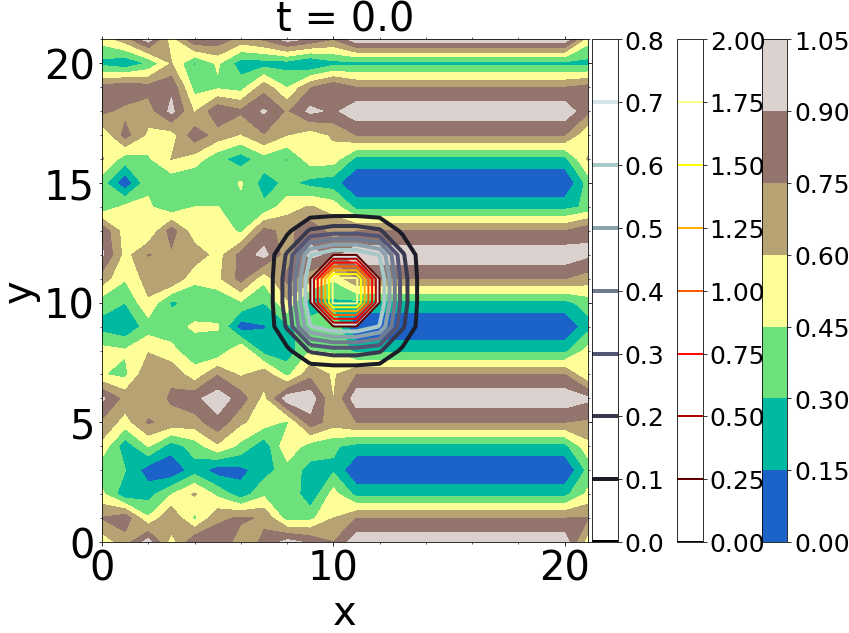}
	\includegraphics[trim = 0cm 0cm 0cm 0cm,clip,scale=.225]{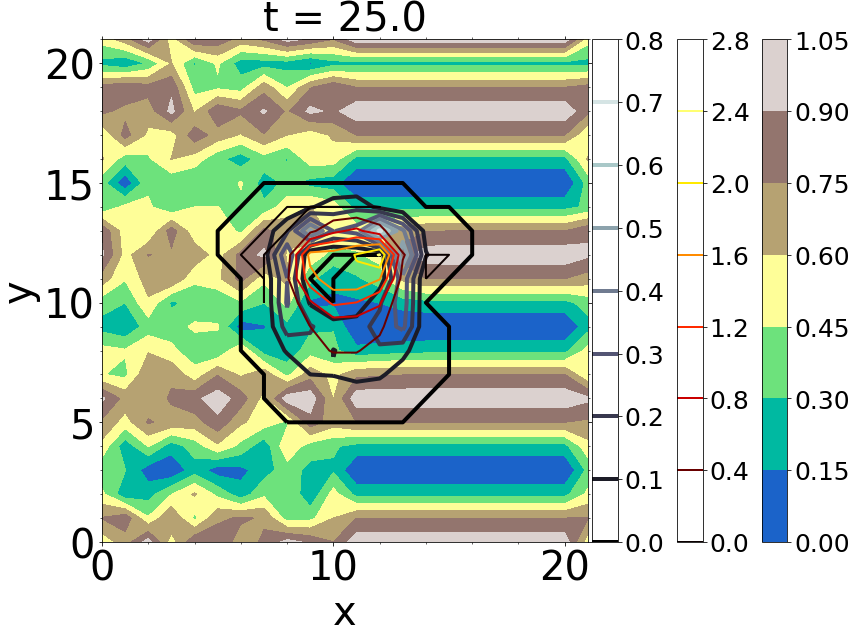}
	\includegraphics[trim = 0cm 0cm 0cm 0cm,clip,scale=.225]{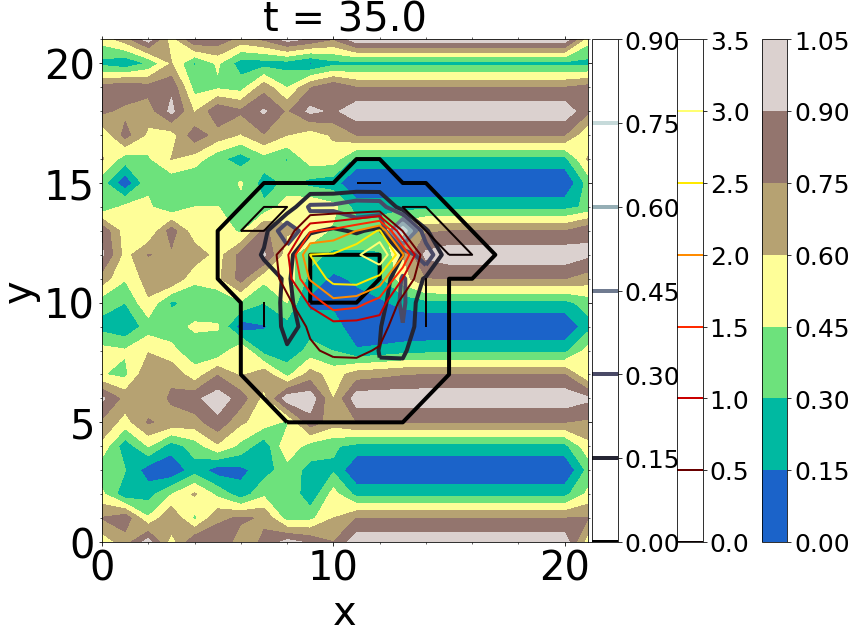}
	\includegraphics[trim = 0cm 0cm 0cm 0cm,clip,scale=.225]{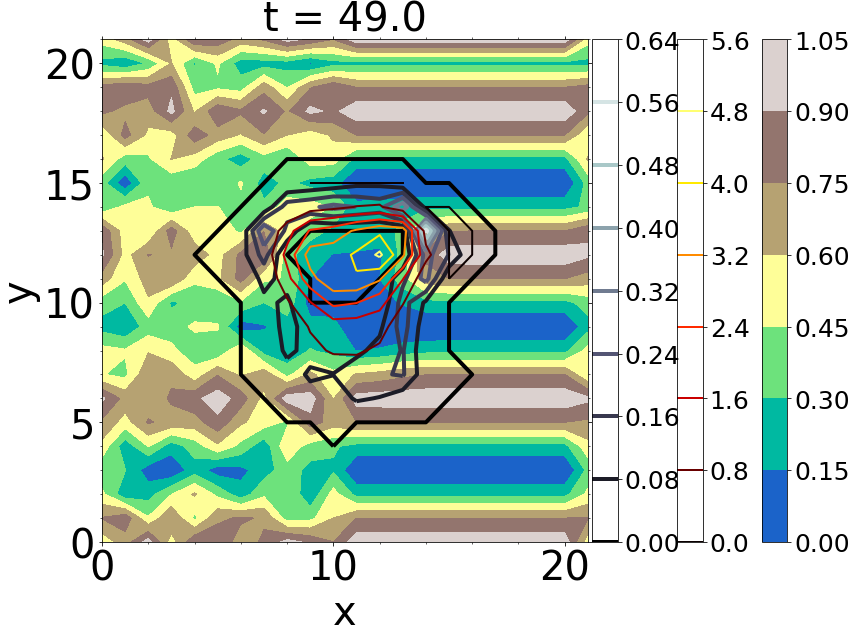}
\caption{Time snapshots of the sample solution number 77.  \label{fig:sim2DMC2}}
\end{figure}
\begin{figure}
	\centering
	\includegraphics[trim = 0cm 0cm 0cm 0cm,clip,scale=.225]{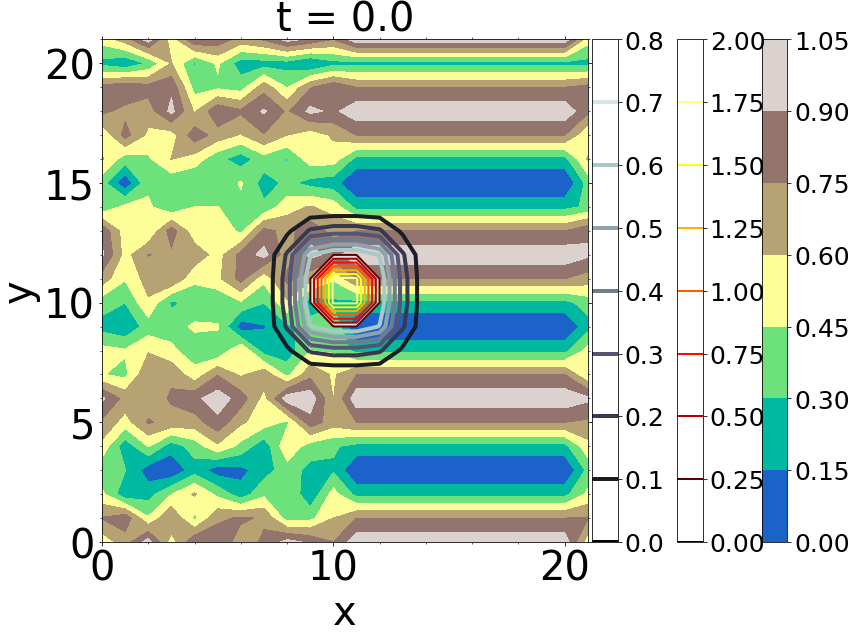}
	\includegraphics[trim = 0cm 0cm 0cm 0cm,clip,scale=.225]{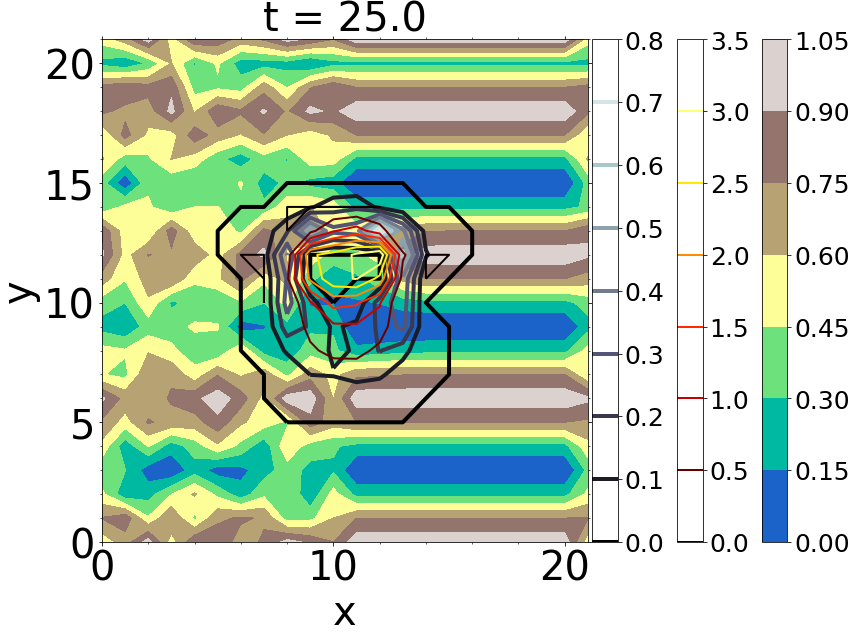}
	\includegraphics[trim = 0cm 0cm 0cm 0cm,clip,scale=.225]{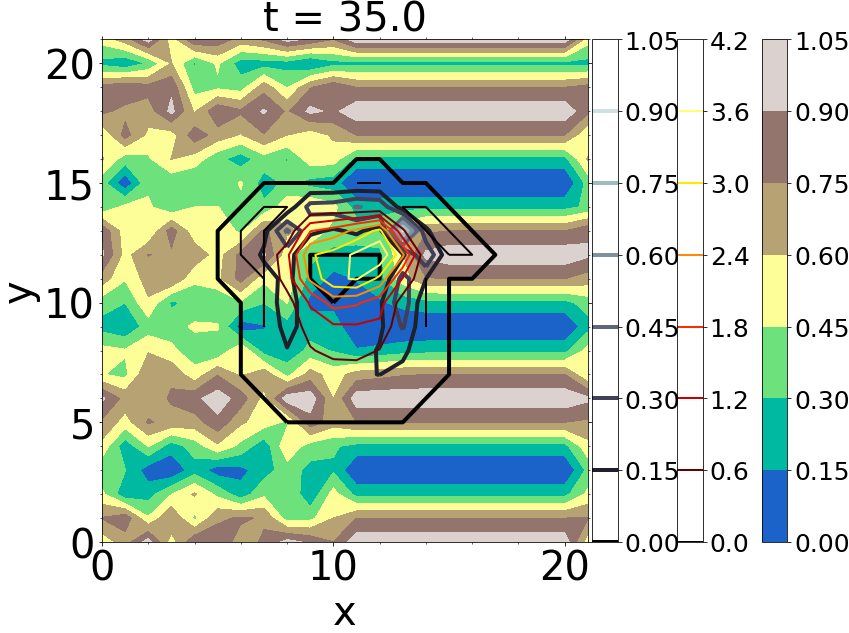}
	\includegraphics[trim = 0cm 0cm 0cm 0cm,clip,scale=.225]{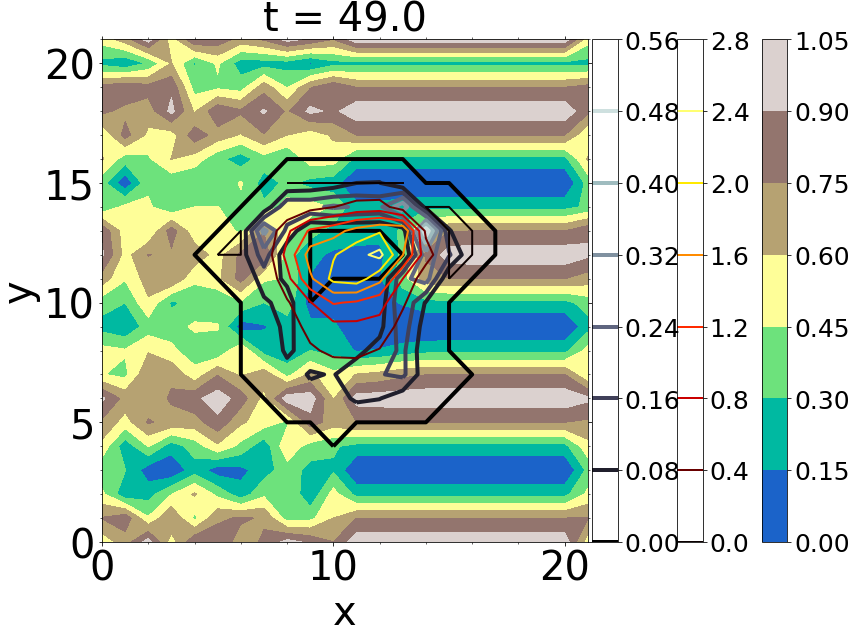}
\caption{Time snapshots of the sample solution number 328.  \label{fig:sim2DMC3}}
\end{figure}
\begin{figure}
	\centering
	\includegraphics[trim = 0cm 0cm 0cm 0cm,clip,scale=.225]{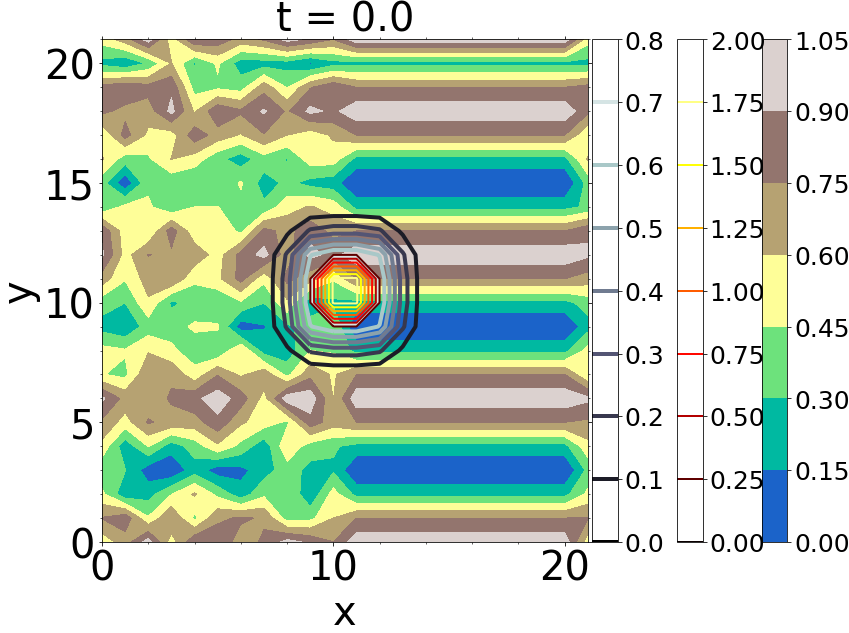}
	\includegraphics[trim = 0cm 0cm 0cm 0cm,clip,scale=.225]{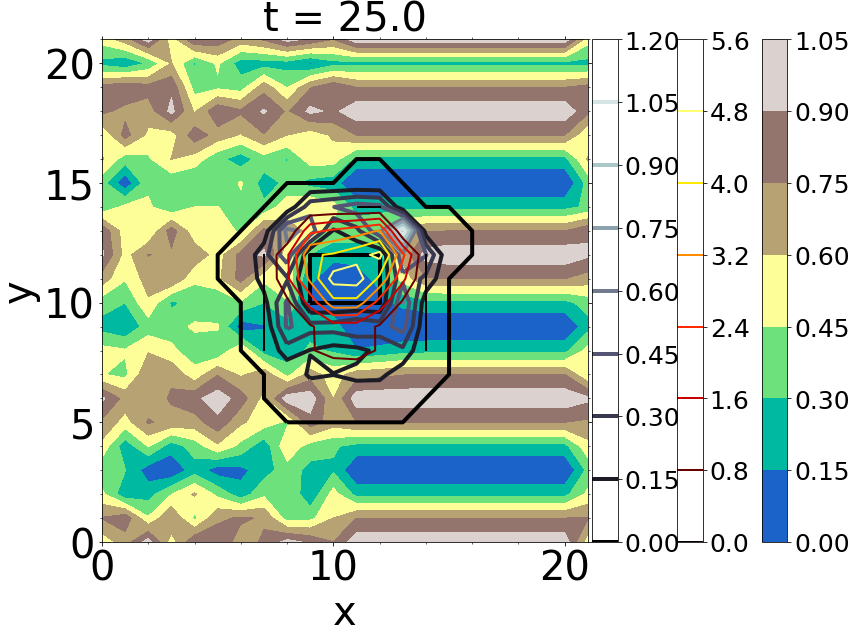}
	\includegraphics[trim = 0cm 0cm 0cm 0cm,clip,scale=.225]{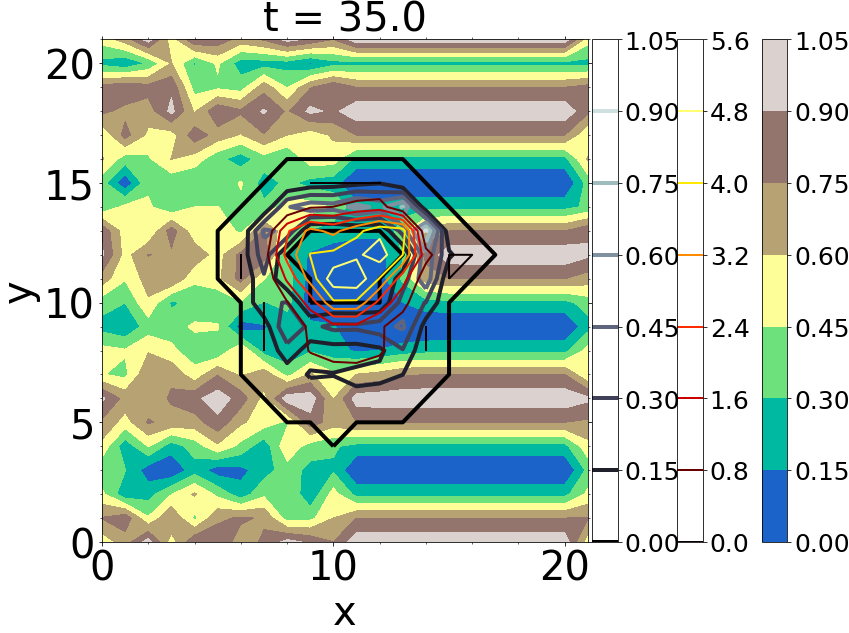}
	\includegraphics[trim = 0cm 0cm 0cm 0cm,clip,scale=.225]{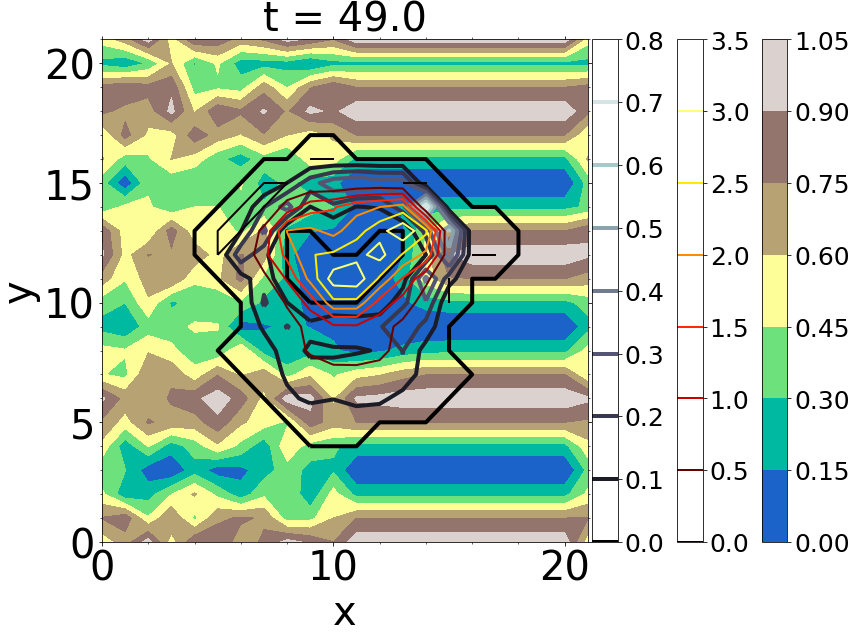}
\caption{Time snapshots of the sample solution number 427.  \label{fig:sim2DMC4}}
\end{figure}
\begin{figure}
	\centering
	\includegraphics[trim = 0cm 0cm 0cm 0cm,clip,scale=.225]{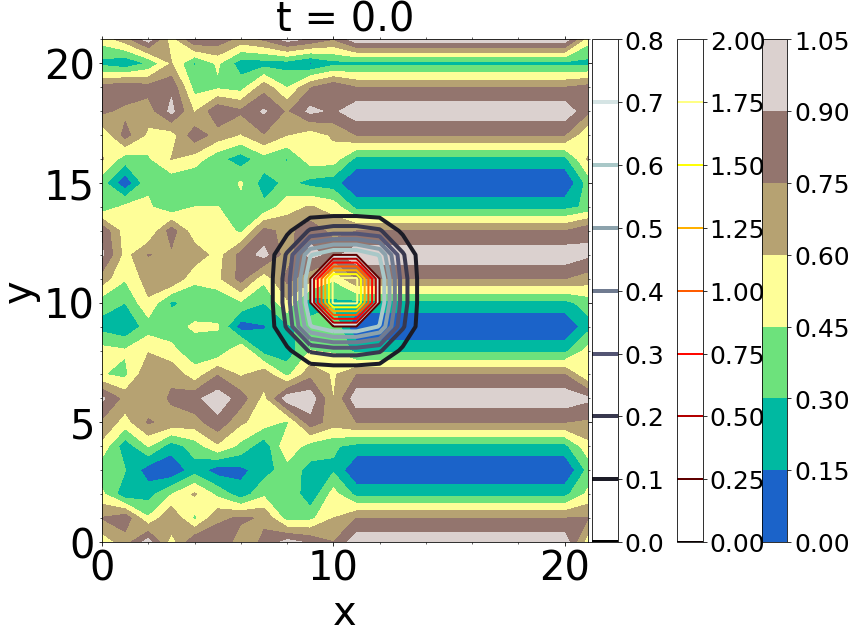}
	\includegraphics[trim = 0cm 0cm 0cm 0cm,clip,scale=.225]{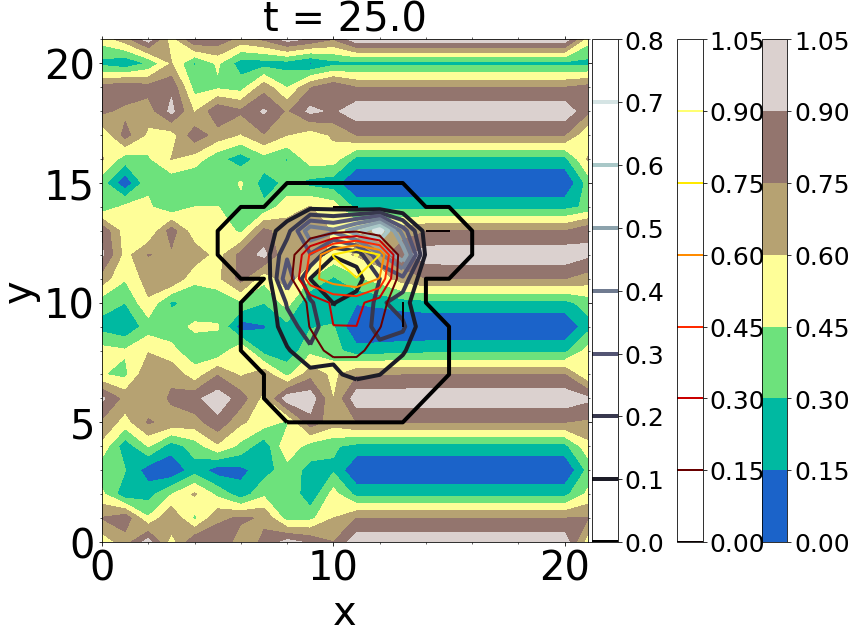}
	\includegraphics[trim = 0cm 0cm 0cm 0cm,clip,scale=.225]{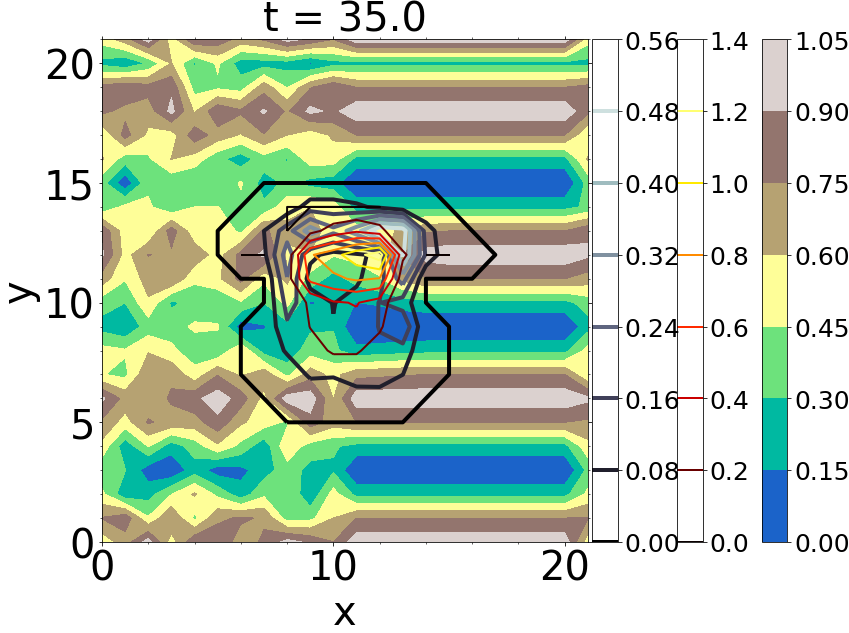}
	\includegraphics[trim = 0cm 0cm 0cm 0cm,clip,scale=.225]{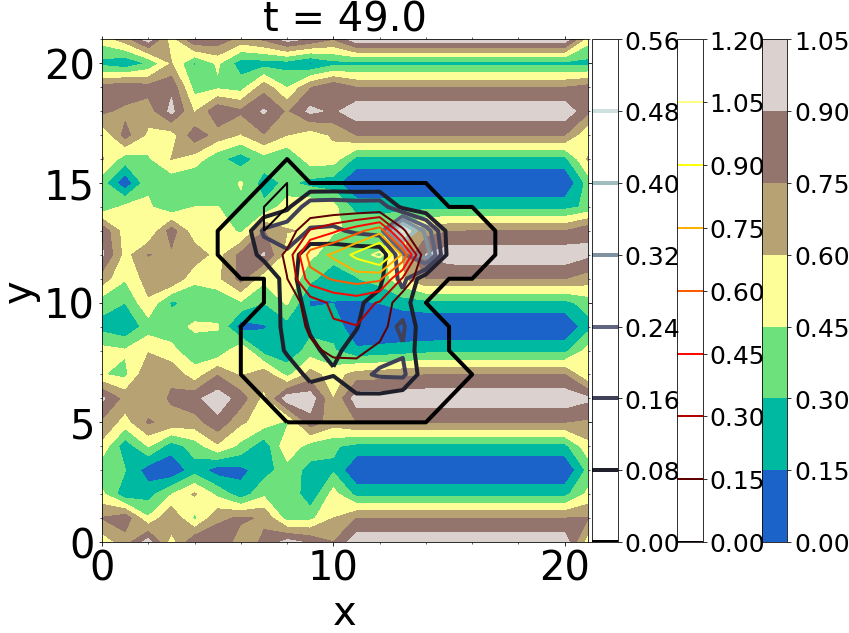}
\caption{Time snapshots of the sample solution number 456.  \label{fig:sim2DMC5}}
\end{figure}
\begin{figure}
	\centering
	\includegraphics[trim = 0cm 0cm 0cm 0cm,clip,scale=.225]{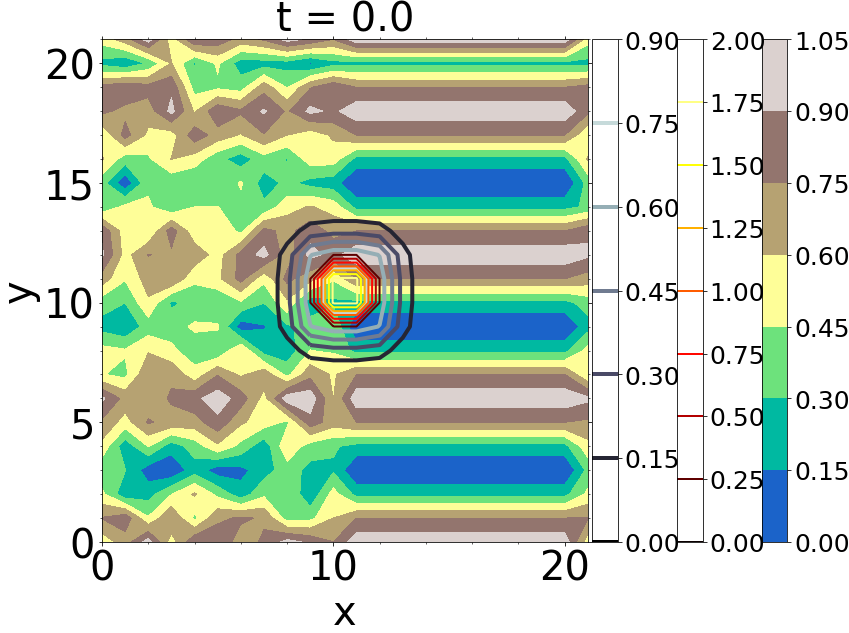}
	\includegraphics[trim = 0cm 0cm 0cm 0cm,clip,scale=.225]{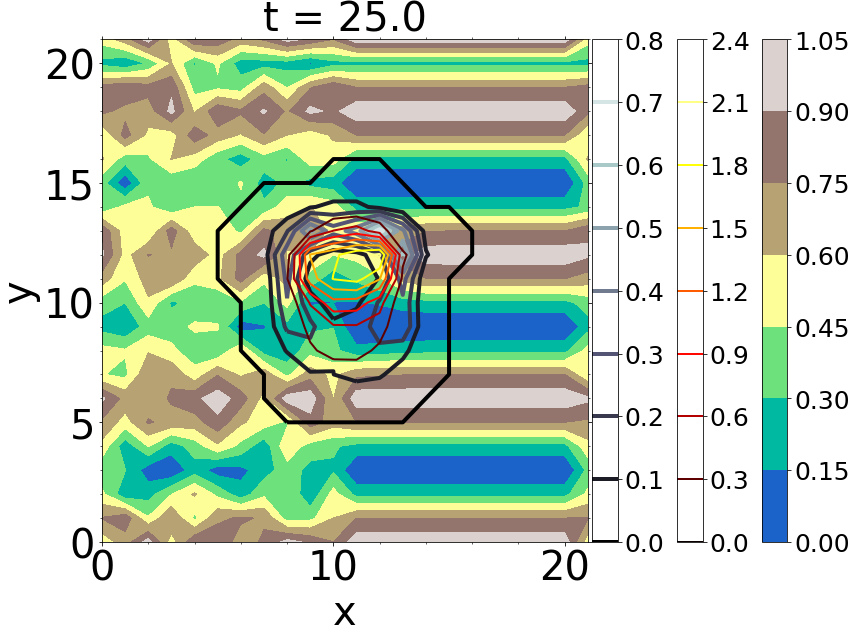}
	\includegraphics[trim = 0cm 0cm 0cm 0cm,clip,scale=.225]{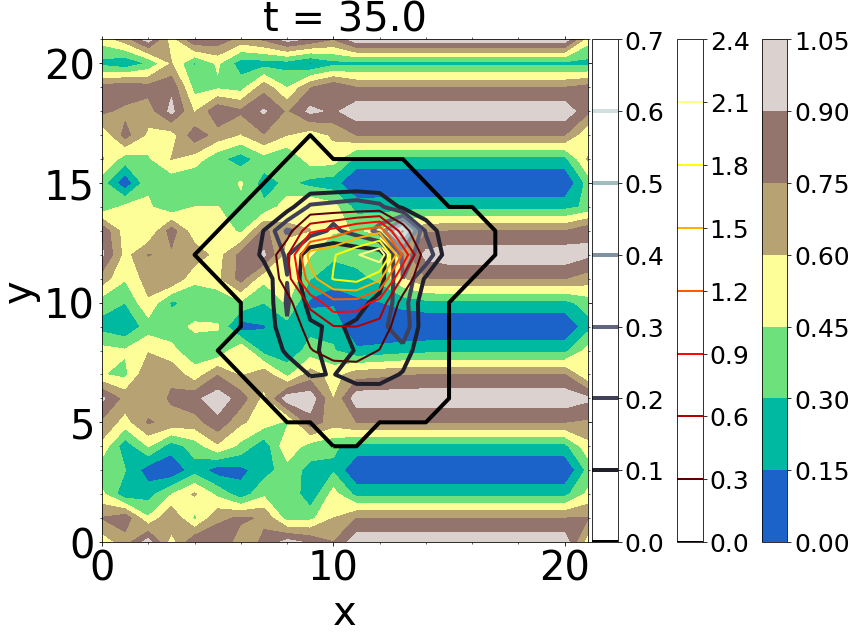}
	\includegraphics[trim = 0cm 0cm 0cm 0cm,clip,scale=.225]{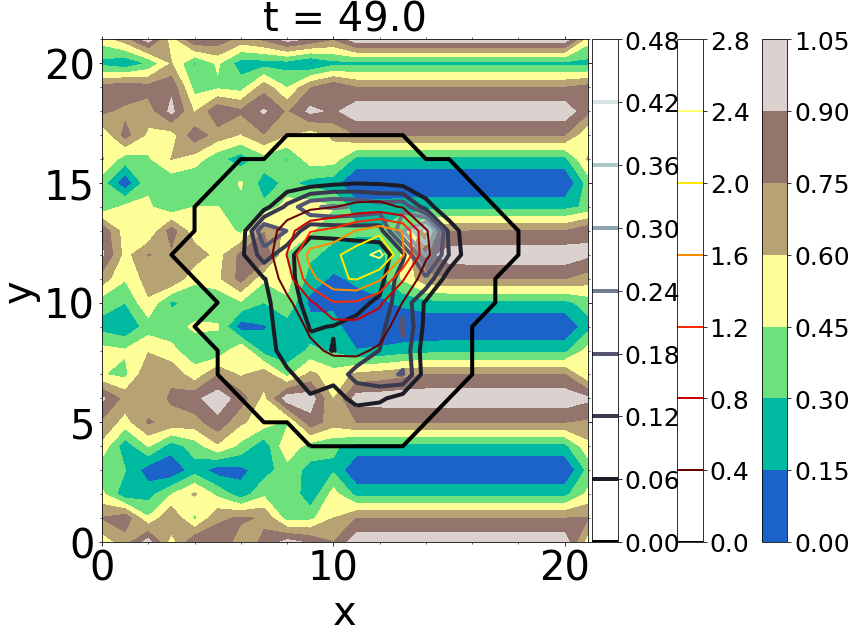}
\caption{Time snapshots of the numerical expectation.  \label{fig:sim2DEx}}
\end{figure}
\noindent
Figure \ref{fig:simFracVsStdDiff} depicts the time snapshots of solution to a simple standalone fractional and standard diffusion equation with periodic boundary conditions. For the fraction diffusion equation, the fraction $\alp$ is set to value $0.75$. Both the equations have the same initial condition and same diffusion coefficient of $0.025$. Based on the figure we see that the speed of the spread is much lower for the fractional diffusion, for the same diffusion coefficient. Thus for a crowded environment, such as cancerous tissue, it makes more sense to use a fractional diffusion operator. \\
Following this we now proceed with the simulation of the stochastic fractional diffusion model \eqref{eq:fracSpdeMod}. The numerically obtained results are shown in Figures \ref{fig:sim2DMC1}-\ref{fig:sim2DMC5}. They depict 5 different sample solutions, randomly chosen from a set of 500 sample solutions. The numerical and model parameters chosen for the simulations are given in Table \ref{tab:simparFracSPDE} and Table \ref{tab:modParFracSPDE}, respectively. Each of the figures consist of 2D-plots representing snapshots at four different time points. Each 2D-plot depicts the concentration distribution of $H$, $C$, and $N$ overlayed on each other. The concentration distribution of $N$ is displayed as a shaded contour plot which is seen as a patchy layer on the background. The concentration distributions of $H$ and $C$ are represented as contour lines,  where the former is represented by thin lines while the latter is represented as thick lines. The color bars on the right side of each plot show the values of the respective contour levels. The leftmost color bar indicates the contour levels of $N$, while the rightmost color bar indicates the contour levels of $H$. The center color bar indicates the contour levels of $C$. Based on these figures we observe the following: \\
{\bf 1.} Starting from a deterministic initial condition, we obtain various different patterns for the spread of cancer. These invasion patterns share the following similarities: \\
{\bf 1a.} In all sample paths, the majority of the tumor mass has the tendency to spread to the upper half of the tissue. This can be confirmed by observing the expected behavior of the sample solutions (as shown in Figure \ref{fig:sim2DEx}). This behavior is due to the haptotaxis term which has a negative velocity, meaning that the cancer tends to move away from high density regions of the tissue. The initial condition is chosen such that the upper part of the central region of the tissue is relatively less dense compared to that of the lower part of the central region of the tissue.\\
{\bf 1b.} A second factor that influences the direction of the invasion is the value of the proton index $H$. Recall that in our model (i.e. \eqref{eq:fracSpdeMod}) $H$ represents the ratio of extracellular and intracellular proton concentration, which provides a measure for the RpHG of the cell.  Since a higher value of the latter promotes cell movement, this feature is captured by the pH taxis term and serves as the primary source of direction for the movement of cells. The advection term in the $H$ equation facilitates propagation of the RpHG
along the tumor edge i.e. away from the dense core  and towards the periphery of the tumor. Thus the advection term in \eqref{eq:fracSpdeModH} and the pH-taxis term in \eqref{eq:fracSpdeModC} together form a  positive feedback which pushes the cancer wave.\\
{\bf 1c.} In the regions of relatively low or negligible proton index value, the spread of cancer is mainly governed by haptotaxis, fractional diffusion, and proliferation.\\
\noindent {\bf 2.} Apart from the above similarities, the stochasticity of the model induces the following interesting patterns:\\
{\bf 2a.} In Figure \ref{fig:sim2DMC2} and Figure \ref{fig:sim2DMC4} we observe that at time $t=49$, at the lower left part of
the tumor core, there is an island-like patch being formed, which is not seen in the expectation plot (Figure \ref{fig:sim2DEx}), hence indicates that this is a rare event. \\
{\bf 2b.} More interestingly, in Figure \ref{fig:sim2DMC1} and Figure \ref{fig:sim2DMC2}, for time $t=25.0$, we observe that the tumor has an opening on its lower side. However, looking at the expectation plot at time $t=25$ we see that the opening is not present. From this we can infer that this is a rare event, as well. \\
Altogether we get a diverse spread of cancer, in the sense that the invasion pattern is irregular, on the one hand due to the less regular fractional operator and on the other hand due to the randomness of the model and the heterogeneity of the tissue. 
\section{Summary}

In this article we have looked at the problem of modeling diversity in complex living systems, such as cancer, from a very fundamental level. We began with motivation from evolutionary biology which provided an intuitive motivation to interpret cancer as a selfish organism. Based on a novel axiomatic formulation of a very generic living system we highlighted the fundamental reason for living systems to diversify. We further justified this reasoning and demonstrated using an illustrative example (see Figures \ref{fig:divSystems} and  \ref{fig:randomRep}) that a diversifying living system has better chances of survival compared to non-diversifying systems. This is evident in competitive games such as football, cricket, volleyball, handball, etc. wherein a team having a set of diverse players has a better chance of winning the game. To be more precise a team having only right handed (legged) players has lower chances of winning when competing against a team that has good a mix of both right and left handed (legged) players. 
For a more biological example one can consider the cases of diseases and epidemies. It is clear that a population with diversified immune systems has better chances of surviving epidemic breakouts compared to a population having identical immune systems for all individuals. \\
Motivated by these observations and having considered cancer as a selfish organism that is playing a game of survival against the host,  we posed the question: how can one model diverse behavior of an organism as observed at the macroscopic scale starting from the microscopic level? This question was answered by Theorem \ref{thm:randACP} in Section \ref{sec:M2MACP}, where we used a functional analytic machinery (see Figure \ref{fig:randOpFrm}) to introduce a framework through which one is able to construct random equations at the macroscopic level. It should be noted here that the equation thus obtained at the macroscopic level is a random abstract Cauchy equation. As a result, in general, the obtained macroscopic equation need not have pointwise evaluation and even worse the operator need not have an explicit computable representation. This is one of the drawbacks of the approach and requires further investigation in order to make it more practical. Nonetheless, it is suitable for applications that involve periodic solutions at the macroscopic level, where one can solve the Cauchy equation in Schwartz space (i.e. on Fourier domain) via Fourier transform.\\ 
For the deduced abstract random Cauchy equation, we formulated Theorem \ref{thm:randACP} which was able to establish the existence of a unique solution which enabled us to go ahead with numerical simulations. Nonetheless, we believe that the current approach for the proving the existence of unique solution is not the most elegant approach. An alternative approach would be to either introduce an appropriate norm $\| \cdot \|_{N}$ on the set of continuous negative definite functions $N$ such that $(N, \| \cdot \|_N)$ is a Banach space or to introduce a metric on the set $N$ so that the notion of continuity on between the elements of $N$ could be established. Then the task would be to find a closed subspace $M \subset N$ which can generate sectorial operators. Once the existence of such a subspace $M$ is established, then the task of constructing random operators would involve specifying an appropriate $M$ valued stochastic process. \\
In Section \ref{sec:simulations} we provided simulations for the dynamics of cancer-invasion at two different levels namely the microscopic and macroscopic levels. In the microscopic simulations we modeled the intrinsic diversification signal by means of switching distribution of the noise term which consequentially conferred better odds for survival for cancer cells. This also provided the motivation for a direct macroscopic simulation where in the random fractional diffusion operator was used to model diverse invasive behavior of cancer cells. Since the micro to macro modeling framework is bidrectional, one has the flexibility to either start at the macro or the micro level.  Starting at a macroscopic equation one could in principle deduce the stochasitc process acting at the microscopic level, by making use of potential theory and the theory of Dirichlet forms [\citenum{MaRockner2012}]. Altogether the framework addresses the fundamental problem of diversification of living systems and introduces some new challenging mathematical problems. \\

\newpage



\bibliographystyle{abbrv}
\bibliography{lits}

\begin{thebibliography}{10}

\bibitem{Applebaum2004}
D.~Applebaum.
\newblock {\em L{\'e}vy Processes and Stochastic Calculus}.
\newblock Cambridge Studies in Advanced Mathematics. Cambridge University
  Press, 2004.

\bibitem{Aristov2019}
V.~Aristov.
\newblock Biological systems as nonequilibrium structures described by kinetic
  methods.
\newblock {\em Results in Physics}, 13:102232, 2019.

\bibitem{EVOCAN15}
A.~Arnal, B.~Ujvari, B.~Crespi, R.~Gatenby, T.~Tissot, M.~Vittecoq, P.~Ewald,
  A.~Casali, H.~Ducasse, C.~Jacqueline, D.~Missé, F.~Renaud, B.~Roche, and
  F.~Thomas.
\newblock Evolutionary perspective of cancer: myth, metaphors, and reality.
\newblock {\em Evolutionary Applications}, 8(6):541--544, 2015.

\bibitem{Bellomo2004}
N.~Bellomo and A.~Bellouquid.
\newblock From a class of kinetic models to the macroscopic equations for
  multicellular systems in biology.
\newblock {\em Discrete and Continuous Dynamical Systems - Series B},
  4(1):59--80, 2004.

\bibitem{Bellomo2008}
N.~Bellomo, A.~Bellouquid, and M.~Delitala.
\newblock From the mathematical kinetic theory of active particles to
  multiscale modelling of complex biological systems.
\newblock {\em Mathematical and Computer Modelling}, 47(7–8):687 -- 698,
  2008.
\newblock Mathematical Methods and Modelling of Biophysical Phenomena.

\bibitem{Bellomo2017}
N.~Bellomo, A.~Bellouquid, L.~Gibelli, and N.~Outada.
\newblock {\em A Quest Towards a Mathematical Theory of Living Systems}.
\newblock Springer International Publishing, 2017.

\bibitem{Bellomo2013}
N.~Bellomo, A.~Bellouquid, J.~Nieto, and J.~Soler.
\newblock Modeling chemotaxis from l–closure moments in kinetic theory of
  active particles.
\newblock 2013.

\bibitem{Bellomo2021}
N.~Bellomo, D.~Burini, G.~Dosi, and et~al.
\newblock What is life? a perspective of the mathematical kinetic theory of
  active particles.
\newblock {\em Mathematical Models and Methods in Applied Sciences},
  31(09):1821--1866, 2021.

\bibitem{Bellomo08}
N.~Bellomo and C.~Dogbe.
\newblock On the modelling crowd dynamics from scaling to hyperbolic
  macroscopic models.
\newblock {\em Mathematical Models and Methods in Applied Sciences},
  18(supp01):1317--1345, 2008.

\bibitem{Bellomo17}
N.~Bellomo and S.~Y. Ha.
\newblock A quest toward a mathematical theory of the dynamics of swarms.
\newblock {\em Mathematical Models and Methods in Applied Sciences},
  27(04):745--770, 2017.

\bibitem{Berg1993}
C.~Berg, K.~Boyadzhiev, and R.~Delaubenfels.
\newblock Generation of generators of holomorphic semigroups.
\newblock {\em Journal of the Australian Mathematical Society. Series A. Pure
  Mathematics and Statistics}, 55(2):246--269, 10 1993.

\bibitem{CasasSelves2011}
M.~Cas\'{a}s-Selves and J.~DeGregori.
\newblock How cancer shapes evolution, and how evolution shapes cancer.
\newblock {\em Evolution (N Y)}, 4(4):624--634, Dec 2011.
\newblock 23705033[pmid].

\bibitem{Chen2015}
H.~Chen, F.~Lin, K.~Xing, and X.~He.
\newblock The reverse evolution from multicellularity to unicellularity during
  carcinogenesis.
\newblock {\em Nature Communications}, 6:6367 EP, Mar 2015.
\newblock Article.

\bibitem{Davies2013}
P.~Davies.
\newblock Exposing cancer’s deep evolutionary roots.
\newblock {\em Physics world: Physics of cancer}, 26(7):37--40, 2013.

\bibitem{DeMasi1991}
A.~De~Masi and E.~Presutti.
\newblock {\em Mathematical Methods for Hydrodynamic Limits}, volume 1501.
\newblock Springer Berlin Heidelberg, Berlin, Heidelberg, 1991.

\bibitem{Delale1982}
C.~F. Delale.
\newblock The hilbert expansion to the boltzmann equation for steady flow.
\newblock {\em Journal of Statistical Physics}, 28(3):589--602, 1982.

\bibitem{DiPerna1989}
R.~J. DiPerna and P.~L. Lions.
\newblock On the cauchy problem for boltzmann equations: Global existence and
  weak stability.
\newblock {\em Annals of Mathematics}, 130(2):321--366, 1989.

\bibitem{Farkas2001}
W.~Farkas, N.~Jacob, and R.~Schilling.
\newblock {\em Function spaces related to continuous negative definite
  functions: $\psi$-Bessel potential spaces}.
\newblock 2001.

\bibitem{Friedl2003}
P.~Friedl and K.~Wolf.
\newblock Tumour-cell invasion and migration: diversity and escape mechanisms.
\newblock {\em Nat Rev Cancer}, 3(5):362--374, May 2003.

\bibitem{GCTFT07}
R.~Gatenby and R.~Gillies.
\newblock Glycolysis in cancer: A potential target for therapy.
\newblock {\em The International Journal of Biochemistry and Cell Biology},
  39(7–8):1358 -- 1366, 2007.

\bibitem{HAWEI11}
D.~Hanahan and R.~Weinberg.
\newblock Hallmarks of cancer: The next generation.
\newblock {\em Cell}, 144(5):646 -- 674, 2011.

\bibitem{Hilbert1916}
D.~Hilbert.
\newblock Begr\"undung der kinetischen gastheorie.
\newblock {\em Mathematische Annalen}, 72:562--577, 1916.

\bibitem{SMAMCI15}
S.~Hiremath and C.~Surulescu.
\newblock A stochastic multiscale model for acid mediated cancer invasion.
\newblock {\em Nonlinear Analysis: Real World Applications}, 22(0):176 -- 205,
  2015.

\bibitem{hs16b}
S.~A. Hiremath and C.~Surulescu.
\newblock Mathematical models for acid-mediated tumor invasion: deterministic
  and stochastic approaches.
\newblock In A.~Gerisch, R.~Penta, and J.~Lang, editors, {\em Multiscale Models
  in Mechano and Tumor Biology: Modeling, Homogenization, and Applications,
  Lecture Notes in Computational Science and Engineering}. Springer Verlag
  Heidelberg.
\newblock accepted.

\bibitem{SAIG16}
S.~A. Hiremath and C.~Surulescu.
\newblock A stochastic model featuring acid-induced gaps during tumor
  progression.
\newblock {\em Nonlinearity}, 29(3):851, 2016.

\bibitem{HiremathZhigun2016}
S.~A. Hiremath, A.~Zhigun, S.~Sonner, and C.~Surulescu.
\newblock On a coupled sde-pde system modeling acid-mediated tumor invasion.
\newblock DCDS-B (accepted), 2017.

\bibitem{Huang2014}
Y.~Huang and A.~Oberman.
\newblock Numerical methods for the fractional laplacian: A finite
  difference-quadrature approach.
\newblock {\em SIAM Journal on Numerical Analysis}, 52(6):3056--3084, 2014.

\bibitem{Iannini2016}
M.~L.~L. Iannini and R.~Dickman.
\newblock Kinetic theory of vehicular traffic.
\newblock {\em American Journal of Physics}, 84(2):135--145, 2016.

\bibitem{Illner1986}
R.~Illner and M.~Pulvirenti.
\newblock Global validity of the boltzmann equation for a two-dimensional rare
  gas in vacuum.
\newblock {\em Communications in Mathematical Physics}, 105(2):189--203, 1986.

\bibitem{Illner1984}
R.~Illner and M.~Shinbrot.
\newblock The boltzmann equation: Global existence for a rare gas in an
  infinite vacuum.
\newblock {\em Communications in Mathematical Physics}, 95(2):217--226, 1984.

\bibitem{Kipnis1999}
C.~Kipnis and C.~Landim.
\newblock {\em {Scaling limits of interacting particle systems}}.
\newblock Springer Verlag, 1999.

\bibitem{Kwasnicki2017}
M.~Kwasnicki.
\newblock Ten equivalent definitions of the fractional laplace operator.
\newblock In {\em Fractional Calculus and Applied Analysis}, volume~20, page~7.
  Degruyter, Feb 2017.
\newblock 1.

\bibitem{Lanford1975}
O.~E. Lanford.
\newblock Time evolution of large classical systems.
\newblock In J.~Moser, editor, {\em Dynamical Systems, Theory and Applications:
  Battelle Seattle 1974 Rencontres}, pages 1--111. Springer Berlin Heidelberg,
  Berlin, Heidelberg, 1975.

\bibitem{MaRockner2012}
Z.~Ma and M.~R{\"o}ckner.
\newblock {\em Introduction to the Theory of (Non-Symmetric) Dirichlet Forms}.
\newblock Universitext. Springer Berlin Heidelberg, 2012.

\bibitem{Merlo2006}
L.~M.~F. Merlo, J.~W. Pepper, B.~J. Reid, and C.~C. Maley.
\newblock Cancer as an evolutionary and ecological process.
\newblock {\em Nat Rev Cancer}, 6(12):924--935, Dec 2006.

\bibitem{MDPMO14}
O.~Miramontes, O.~DeSouza, L.~Paiva, A.~Marins, and S.~Orozco.
\newblock Lévy flights and self-similar exploratory behaviour of termite
  workers: Beyond model fitting.
\newblock {\em PLoS ONE}, 9(10):1--9, 10 2014.

\bibitem{Morrey54}
C.~B. Morrey.
\newblock On the derivation of the equations of hydrodynamics from statistical
  mechanics.
\newblock {\em Proceedings of the National Academy of Sciences of the United
  States of America}, 40(5):317--322, 1954.

\bibitem{Niklas2014}
K.~J. Niklas.
\newblock The evolutionary-developmental origins of multicellularity.
\newblock {\em American Journal of Botany}, 101(1):6--25, 2014.

\bibitem{MYECODEV15}
B.~Okamura, A.~Gruhl, and J.~Bartholomew.
\newblock {\em Myxozoan Evolution, Ecology and Development}.
\newblock Springer International Publishing, 2015.

\bibitem{Olla1993}
S.~Olla, S.~R.~S. Varadhan, and H.-T. Yau.
\newblock Hydrodynamical limit for a hamiltonian system with weak noise.
\newblock {\em Comm. Math. Phys.}, 155(3):523--560, 1993.

\bibitem{PAZY83}
A.~Pazy.
\newblock {\em Semigroups of Linear Operators and Applications to Partial
  Differential Equations}.
\newblock Applied Mathematical Sciences. Springer New York, 1983.

\bibitem{EVOLBIOCONAPP08}
P.~Pontarotti.
\newblock {\em Evolutionary Biology from Concept to Application}.
\newblock Springer Berlin Heidelberg, 2008.

\bibitem{Prevot2007}
C.~Pr{\'{e}}v{\^{o}}t and M.~R{\"{o}}ckner.
\newblock {\em {A concise course on stochastic partial differential
  equations}}, volume 1905 of {\em Lecture Notes in Mathematics}.
\newblock Springer, Berlin, 2007.

\bibitem{Pulvirenti1987}
M.~Pulvirenti.
\newblock Global validity of the boltzmann equation for a three-dimensional
  rare gas in vacuum.
\newblock {\em Communications in Mathematical Physics}, 113(1):79--85, 1987.

\bibitem{WDTHAG76}
E.~Racker.
\newblock Why do tumor cells have a high aerobic glycolysis?
\newblock {\em Journal of Cellular Physiology}, 89(4):697--700, 1976.

\bibitem{REY15}
A.~Reynolds.
\newblock Liberating {L\'evy} walk research from the shackles of optimal
  foraging.
\newblock {\em Physics of Life Reviews}, 14:59 -- 83, 2015.

\bibitem{Samko2001}
S.~Samko.
\newblock {\em Hypersingular Integrals and Their Applications}.
\newblock Analytical Methods and Special Functions. Taylor \& Francis, 2001.

\bibitem{Schilling2012}
R.~Schilling, R.~Song, and Z.~Vondracek.
\newblock {\em Bernstein Functions: Theory and Applications}.
\newblock De Gruyter Studies in Mathematics. De Gruyter, 2012.

\bibitem{Varadhan1993}
S.~R.~S. Varadhan.
\newblock Entropy methods in hydrodynamic scaling.
\newblock In C.~Cercignani and M.~Pulvirenti, editors, {\em Nonequilibrium
  Problems in Many-Particle Systems: Lectures given at the 3rd Session of the
  Centro Internazionale Matematico Estivo (C.I.M.E.) held in Montecatini,
  Italy, June 15--27, 1992}, pages 112--145. Springer Berlin Heidelberg,
  Berlin, Heidelberg, 1993.

\bibitem{VISRAPLU08}
G.~Viswanathan, E.~Raposo, and M.~da~Luz.
\newblock L\'evy flights and superdiffusion in the context of biological
  encounters and random searches.
\newblock {\em Physics of Life Reviews}, 5(3):133 -- 150, 2008.

\bibitem{VISAF96}
G.~M. Viswanathan, V.~Afanasyev, S.~V. Buldyrev, E.~J. Murphy, P.~A. Prince,
  and H.~E. Stanley.
\newblock L\'evy flight search patterns of wandering albatrosses.
\newblock {\em Nature}, 381:413--415, 1996.

\bibitem{YAGI09}
A.~Yagi.
\newblock {\em Abstract Parabolic Evolution Equations and Their Applications}.
\newblock Springer Monographs in Mathematics. Springer, 2009.

\bibitem{Yates2012}
L.~R. Yates and P.~J. Campbell.
\newblock Evolution of the cancer genome.
\newblock {\em Nat Rev Genet}, 13(11):795--806, Nov 2012.
\newblock 23044827[pmid].

\end{thebibliography}

\end{document}